\documentclass[final,onefignum,onetabnum]{siamart220329}

\usepackage{lipsum}
\usepackage{amsfonts}
\usepackage[export]{adjustbox}
\usepackage{graphicx}
\usepackage{epstopdf}
\usepackage{algorithmic}
\ifpdf
  \DeclareGraphicsExtensions{.eps,.pdf,.png,.jpg}
\else
  \DeclareGraphicsExtensions{.eps}
\fi

\newsiamremark{remark}{Remark}
\newsiamremark{hypothesis}{Hypothesis}
\crefname{hypothesis}{Hypothesis}{Hypotheses}
\newsiamthm{claim}{Claim}

\headers{First moments of a polyhedron clipped by a paraboloid}{F. Evrard, R. Chiodi, A. Han, B. van Wachem, O. Desjardins}

\title{First moments of a polyhedron clipped by a paraboloid}

\author{Fabien Evrard$^{\ddagger}$\thanks{Sibley School of Mechanical and Aerospace Engineering, Cornell University, {Ithaca, NY 14853, United States}}
\and Robert Chiodi\thanks{Computer, Computational, and Statistical Sciences Division, Los Alamos National Laboratory, {Los Alamos, NM 87545, United States}}  
\and Austin Han\footnotemark[1]
\and Berend van Wachem\thanks{Lehrstuhl f\"ur Mechanische Verfahrenstechnik, Otto-von-Guericke-Universit\"at Magdeburg, {Universit\"atsplatz 2, 39106 Magdeburg, Germany}}
\and Olivier Desjardins\footnotemark[1]}

\usepackage{amsopn}
\DeclareMathOperator{\diag}{diag}

\ifpdf
\hypersetup{
  pdftitle={Moments of the intersection between a paraboloid and a polyhedron},
  pdfauthor={F. Evrard, R. Chiodi, A. Han, B. van Wachem and O. Desjardins}
}
\fi

\usepackage{geometry}
\usepackage{amsmath,amssymb}
\usepackage[caption=false]{subfig}
\usepackage{multirow}
\usepackage{tikz,tikz-3dplot,tkz-euclide}
\usetikzlibrary{calc}
\usetikzlibrary{arrows}
\usetikzlibrary{shapes.geometric, arrows, intersections, through}
\usetikzlibrary{decorations.text, decorations.shapes, backgrounds}
\usetikzlibrary{arrows.meta}
\usetikzlibrary{colorbrewer}
\usetikzlibrary{patterns}
\usepackage{pgfplots} 
\usepackage{pgfplotstable}
\usepgfplotslibrary{colorbrewer}
\usepgfplotslibrary{colormaps}
\usepgfplotslibrary{patchplots}
\usetikzlibrary{colorbrewer}
\pgfplotsset{compat=newest}
\usepackage{ifthen}
\usepackage{setspace}
\usetikzlibrary{cd, arrows, matrix, intersections, math, calc}
\usepgfplotslibrary{external} 
\tikzset
  {midarrow/.style={decoration={markings,mark=at position 0.5 with
     {\arrow[xshift=1pt,scale=0.6]{angle 60[length=2pt]}}},postaction={decorate}}
  }
  \tikzset
  {revmidarrow/.style={decoration={markings,mark=at position 0.5 with
     {\arrowreversed[xshift=-1pt,scale=0.6]{angle 60[length=2pt]}}},postaction={decorate}}
  }
\usepackage[export]{adjustbox}
\usepackage{enumitem}
\usetikzlibrary{tikzmark}
\usetikzlibrary{decorations.pathreplacing,decorations.markings}

\tikzset{
  on each segment/.style={
    decorate,
    decoration={
      show path construction,
      moveto code={},
      lineto code={
        \path [#1]
        (\tikzinputsegmentfirst) -- (\tikzinputsegmentlast);
      },
      curveto code={
        \path [#1] (\tikzinputsegmentfirst)
        .. controls
        (\tikzinputsegmentsupporta) and (\tikzinputsegmentsupportb)
        ..
        (\tikzinputsegmentlast);
      },
      closepath code={
        \path [#1]
        (\tikzinputsegmentfirst) -- (\tikzinputsegmentlast);
      },
    },
  },
  mid arrow/.style={postaction={decorate,decoration={
        markings,
        mark=at position .5 with {\arrow[#1]{stealth}}
      }}},
}
\tikzset{
    set arrow inside/.code={\pgfqkeys{/tikz/arrow inside}{#1}},
    set arrow inside={end/.initial=>, opt/.initial=},
    /pgf/decoration/Mark/.style={
        mark/.expanded=at position #1 with
        {
            \noexpand\arrow[\pgfkeysvalueof{/tikz/arrow inside/opt}]{\pgfkeysvalueof{/tikz/arrow inside/end}}
        }
    },
    arrow inside/.style 2 args={
        set arrow inside={#1},
        postaction={
            decorate,decoration={
                markings,Mark/.list={#2}
            }
        }
    },
}

\makeatletter
\renewcommand*\env@matrix[1][\arraystretch]{%
  \edef\arraystretch{#1}%
  \hskip -\arraycolsep
  \let\@ifnextchar\new@ifnextchar
  \array{*\c@MaxMatrixCols c}}
\makeatother

\begin{document}

\overfullrule=0mm

\maketitle

\vspace{12mm}

\begin{abstract}
  We provide closed-form expressions for the first moments (i.e., the volume and volume-weighted centroid) of a polyhedron clipped by a paraboloid, that is, of a polyhedron intersected with the subset of the three-dimensional real space located on one side of a paraboloid. These closed-form expressions are derived following successive applications of the divergence theorem and the judicious parametrization of the intersection of the polyhedron's faces with the paraboloid. We provide means for identifying ambiguous discrete intersection topologies, and propose a corrective procedure for preventing their occurence. Finally, we put our proposed closed-form expressions and numerical approach to the test with millions of random and manually engineered polyhedron/paraboloid intersection configurations. The results of these tests show that we are able to provide robust machine-accurate estimates of the first moments at a computational cost that is within one order of magnitude of that of state-of-the-art half-space clipping algorithms.\\

  \noindent © 2024. This manuscript version is made available under the CC-BY-NC-ND 4.0 license.\\
  \noindent \url{http://creativecommons.org/licenses/by-nc-nd/4.0}
\end{abstract}



\vspace{4mm}

\section{Introduction}
\label{sec:introduction}
Many computational methods and applications, ranging from finite-element \cite{Sevilla2008,Burman2015}, cut-cell discontinuous Galerkin \cite{Engwer2020}, and immersed isometric analysis methods \cite{Hughes2005,Antolin2021}, 
to the initialization \cite{Bna2015,Jones2019,Kromer2019,Kromer2021a,Tolle2022} and transport \cite{Renardy2002} of interfaces for simulating gas-liquid flows, require estimating integrals over polyhedra that are clipped by curved surfaces. 
These applications have engendered multiple dedicated quadrature rules and integration strategies, most of which focusing on estimating the first few moments of these clipped polyhedra, thus considering polynomial integrands only.
The numerical approaches employed to estimate these moments vary greatly in terms of accuracy, computational cost, and robustness. Monte-Carlo methods \cite{Evans2000,Hahn2005} are extremely robust and straightforward to implement, however, they suffer from a poor convergence rate, hence their cost/accuracy ratio is significant. Approaches based on octree subdivision \cite{Abedian2013,Kudela2016,Divi2020} or surface triangulation/tesselation \cite{Tolle2022} exhibit better convergence rates, yet their computational cost remains prohibitive for numerical applications requiring \textit{on-the-fly} moment estimations. A number of recent approaches rely on successive applications of the divergence theorem, converting the first moments of a clipped polyhedron into two- and/or one-dimensional integrals. These integrals can then be numerically integrated at low computational cost or, for specific surface types, even be derived into closed-form expressions. 
Bn\`a~\textit{et al.}~\cite{Bna2015} estimate the volume (zeroth moment) of a cube clipped by an implicit surface represented with a level-set function through integrating the local height of the surface using a two-dimensional Gauss-Legendre quadrature rule. This work has been extended to the first moments of a clipped cuboid by Chierici~\textit{et al.}~\cite{Chierici2022}.
For the similar purpose of estimating the zeroth moment of a polyhedron clipped by an implicit surface, Jones~\textit{et al.}~\cite{Jones2019} decompose the clipped polyhedron into a set of simplices, themselves split into a reference polyhedron whose volume is computed analytically, and a set of fundamental curved domains whose volumes are estimated using a two-dimensional Gauss-Legendre quadrature rule. Chin and Sukumar~\cite{Chin2020} use a Duvanant quadrature rule \cite{Dunavant1985} for integrating over the faces of a polyhedron bounded by rational or non-rational B\'ezier and B-spline patches. For non-rational surface parametrizations, this yields exact integral estimations of polynomial integrands, provided that the order of the Duvanant rule is high enough. Kromer~and~Bothe~\cite{Kromer2019,Kromer2021a} estimate the zeroth moment of a polyhedron clipped by an implicit surface by locally approximating the implicit surface as a paraboloid and applying the divergence theorem twice, converting the clipped volume into a sum of one-dimensional integrals, which are then estimated with a Gauss-Legendre quadrature rule. Finally, using Berstein basis functions instead of monomial ones, Antolin~and~Hirschler~\cite{Antolin2021} recently showed that, following successive applications of the divergence theorem, polynomial integrands can be integrated in a straightforward and analytical manner over curved polyhedra bounded by non-rational B\'ezier or B-spline surfaces.\medskip

This manuscript is concerned with estimating the first moments of a specific type of curved polyhedra, that are planar non-convex polyhedra clipped by a paraboloid surface (as in~\cite{Kromer2019,Kromer2021a}). Moreover, we require this estimation to (a) reach machine accuracy, while (b) maintaining a computational expense that is low enough to enable its \textit{on-the-fly} execution in typical numerical applications (e.g., the simulation of two-phase flows with finite-volumes), and (c) being robust to singular configurations (e.g., paraboloid surfaces being parabolic cylinder or planes, and/or ambiguous discrete intersection topologies). These choices and requirements, although mainly motivated by the use of these moments for simulating two-phase flows, may also find applications in the numerical fields listed above. 
A main difficulty in clipping a polyhedron with a paraboloid lies in the fact that the faces of the clipped polyhedron cannot systematically be represented with non-rational or rational B\'ezier patches \cite{Teller1990,Lodha1990}. This prevents the use of recently proposed integration strategies designed for curved polyhedra bounded by B\'ezier or B-spline surfaces \cite{Chin2020,Antolin2021}. By successive applications of the divergence theorem, we show that the first moments of the clipped polyhedron can be expressed as a sum of one-dimensional integrals over straight line segments and conic section arcs. With a parametrization of the latter into rational B\'ezier arcs, we derive closed-form expressions for the first moments, rendering their numerical estimation exact within machine accuracy. Implemented within the half-edge data structure of the open-source \texttt{Interface Reconstruction Library} \cite{Chiodi2022}\footnote{The \texttt{Interface Reconstruction Library~(IRL)} source code is available under Mozilla Public License~2.0 (MPL-2.0) at \href{https://github.com/robert-chiodi/interface-reconstruction-library/tree/paraboloid_cutting}{\texttt{https://github.com/robert-chiodi/interface-reconstruction-library/tree/paraboloid\_cutting}}.}, the computational cost of these moment estimations is kept within an order of magnitude of that of clipping a polyhedron with a half-space. Finally, our choice of arc parametrization, in conjuction with the detection and treatment of ambiguous discrete topologies, allows for robust moment estimates even in degenerate configurations.\medskip

The remainder of this manuscript is organized as follows: \Cref{sec:problem} introduces the problem that we address and the notations employed throughout the manuscript. The closed-form expressions of the clipped polyhedron's first moments are derived in \Cref{sec:moments}. \Cref{sec:surface} touches upon the integration of quantities (e.g., the moments) of the clipped polyhedron's curved face(s). \Cref{sec:robustness} details the procedure employed for preventing ambiguous clipped polyhedron topologies. \Cref{sec:highorder} discusses the extension of our approach to higher-order moments. Finally, the accuracy, efficiency, and robustness of our proposed integration stategy are assessed in \Cref{sec:tests}, and we draw conclusions in \Cref{sec:conclusions}.



\section{Problem statement}
\label{sec:problem}

\begin{figure} \centering
    \subfloat[A polyhedron $\mathcal{P}\subset\mathbb{R}^3$. In this example, $\mathcal{P}$ is a regular dodecahedron with $n_{\!\mathcal{F}} = 12$ faces.]{\adjincludegraphics[width=0.45\textwidth,trim=0 {0.17\height} 0 0,clip=true]{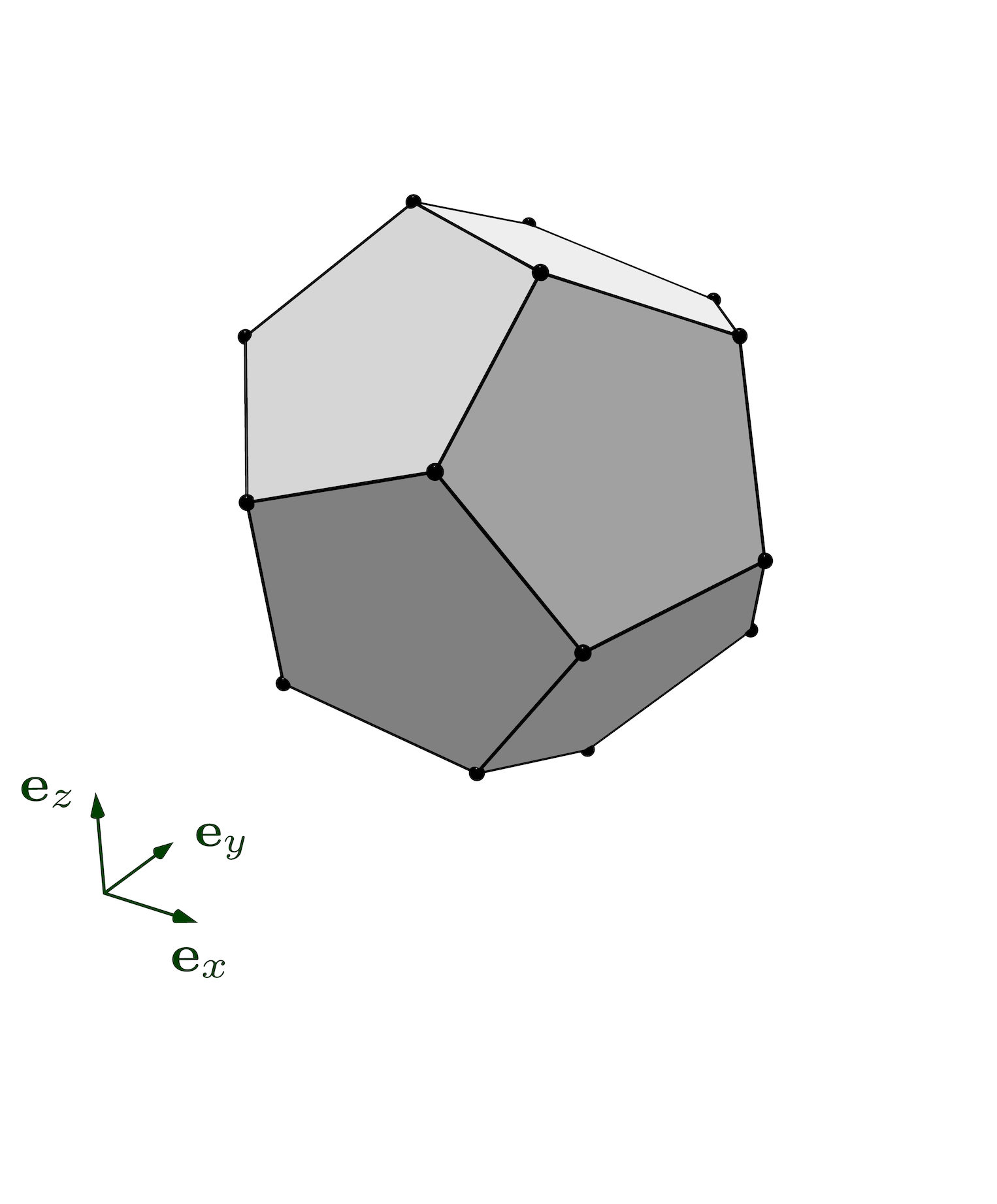} \label{fig:poly}}\quad
    \subfloat[The faces of $\smash{\partial\mathcal{P} = \cup_{i}\mathcal{F}_i}$ and their outward-pointing normal vectors, $\smash{\mathbf{n}_i}$, $\smash{i \in \{1,\ldots,n_{\!\mathcal{F}}\}}$.]{\adjincludegraphics[width=0.45\textwidth,trim=0 {0.17\height} 0 0,clip=true]{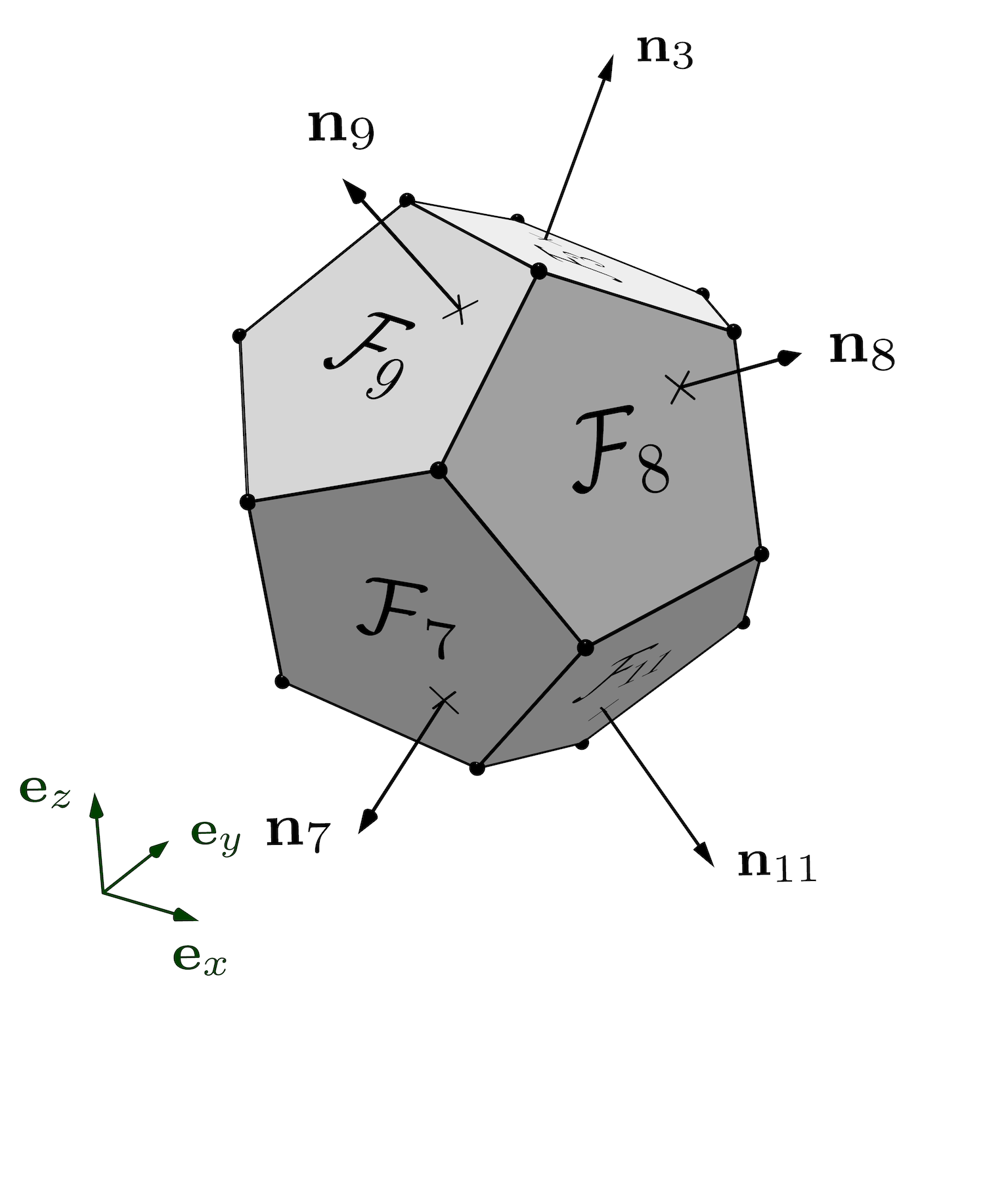} \label{fig:poly_withfaces}}\\
    \subfloat[A paraboloid surface $\mathcal{S}$ intersecting $\mathcal{P}$. The domain $\mathcal{Q}$ is the subset of $\mathbb{R}^3$ located below $\mathcal{S}$, with respect to $\mathbf{e}_z$.]{\adjincludegraphics[width=0.45\textwidth,trim=0 {0.17\height} 0 {0.07\height},clip=true]{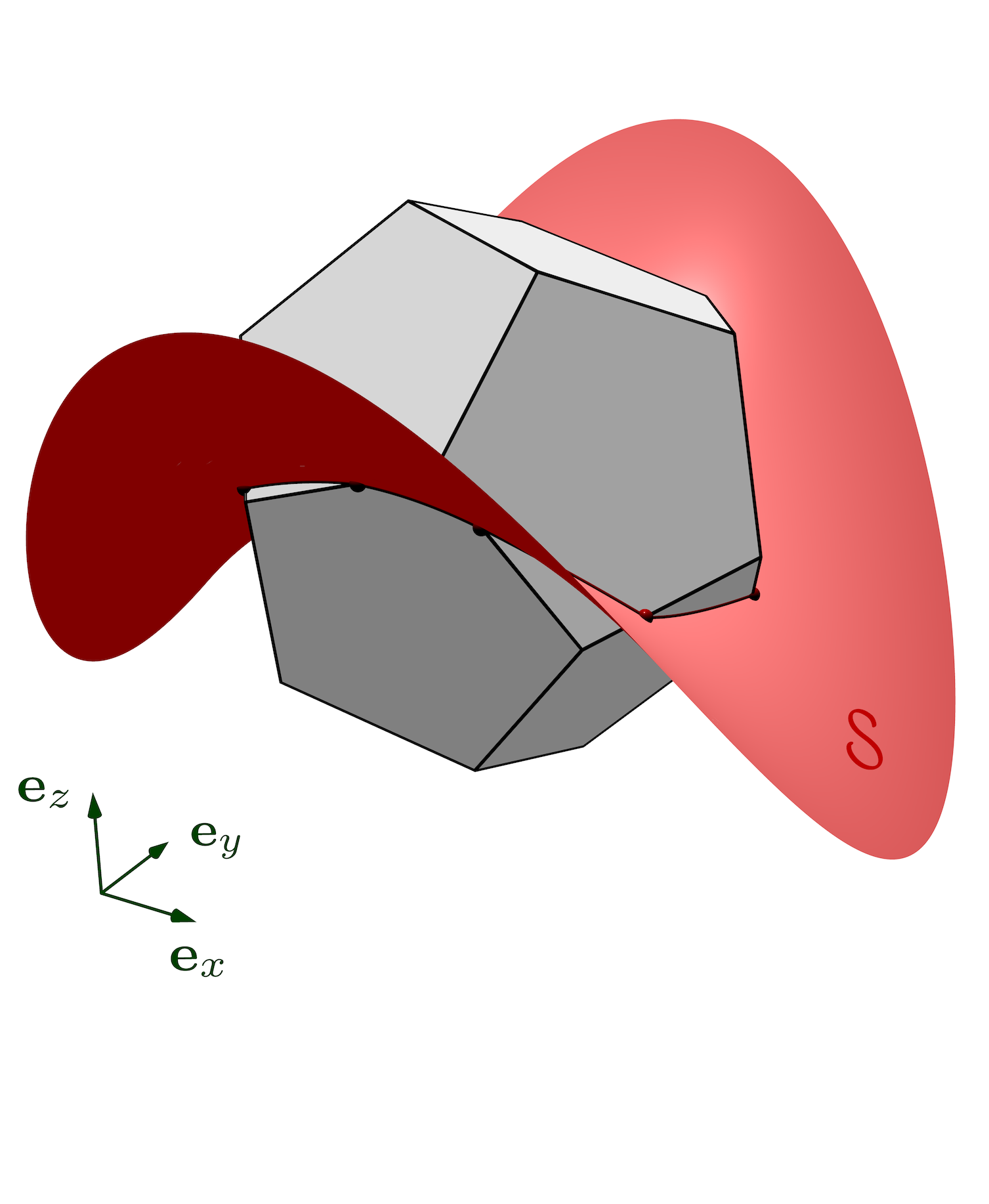}  \label{fig:poly_and_para}}\quad
    \subfloat[The intersection of the paraboloid $\mathcal{S}$ with the polyhedron $\mathcal{P}$, $\smash{\tilde{\mathcal{S}} = \mathcal{S} \cap \mathcal{P}}$.]{\adjincludegraphics[width=0.45\textwidth,trim=0 {0.17\height} 0 {0.07\height},clip=true]{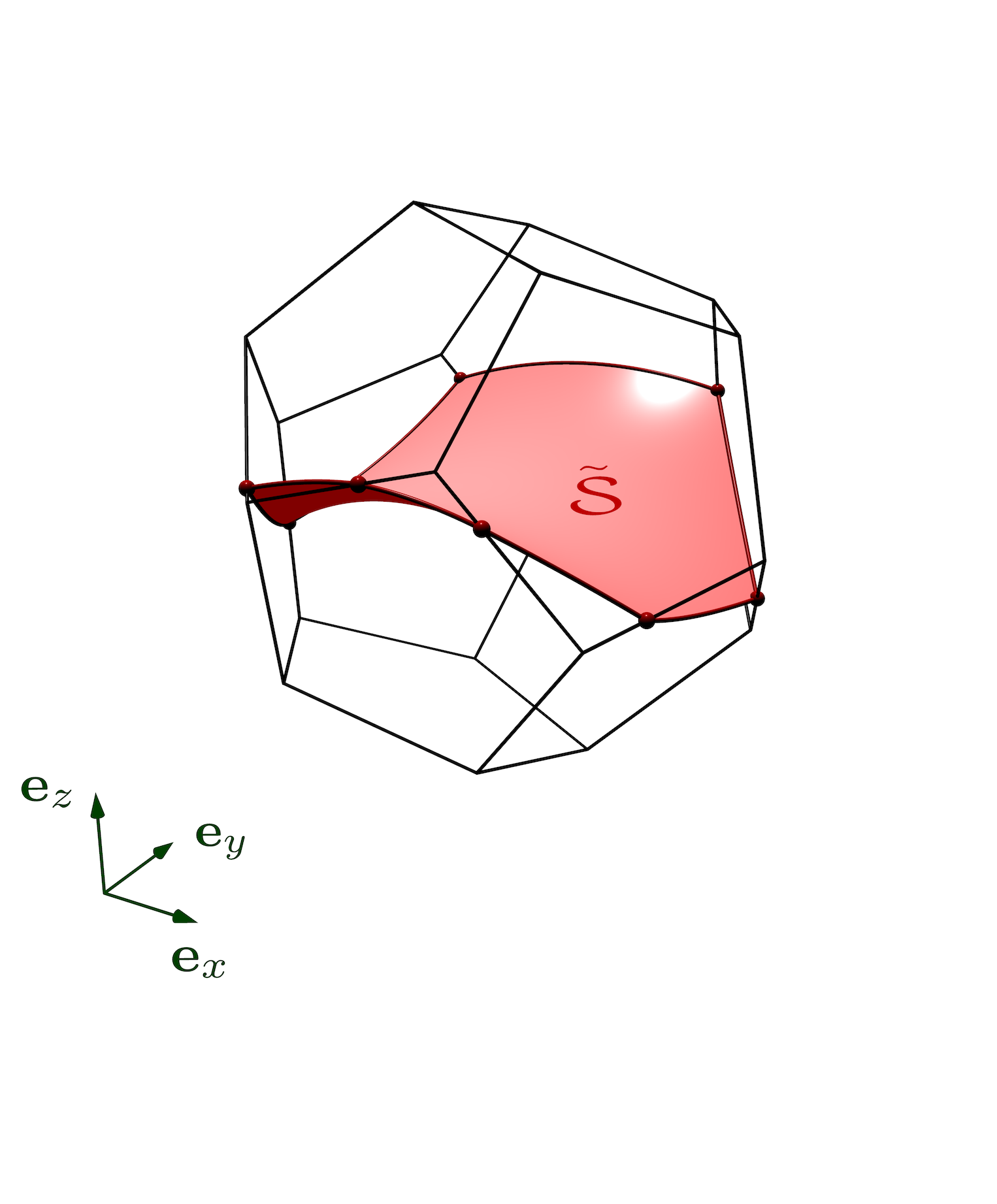}  \label{fig:poly_inter_surf}}\\
    \subfloat[The clipped polyhedron $\smash{\hat{\mathcal{P}} = \mathcal{P} \cap \mathcal{Q}}$, where $\mathcal{Q}$ is the subset of $\mathbb{R}^3$ below $\mathcal{S}$.]{\adjincludegraphics[width=0.45\textwidth,trim=0 {0.17\height} 0 {0.25\height},clip=true]{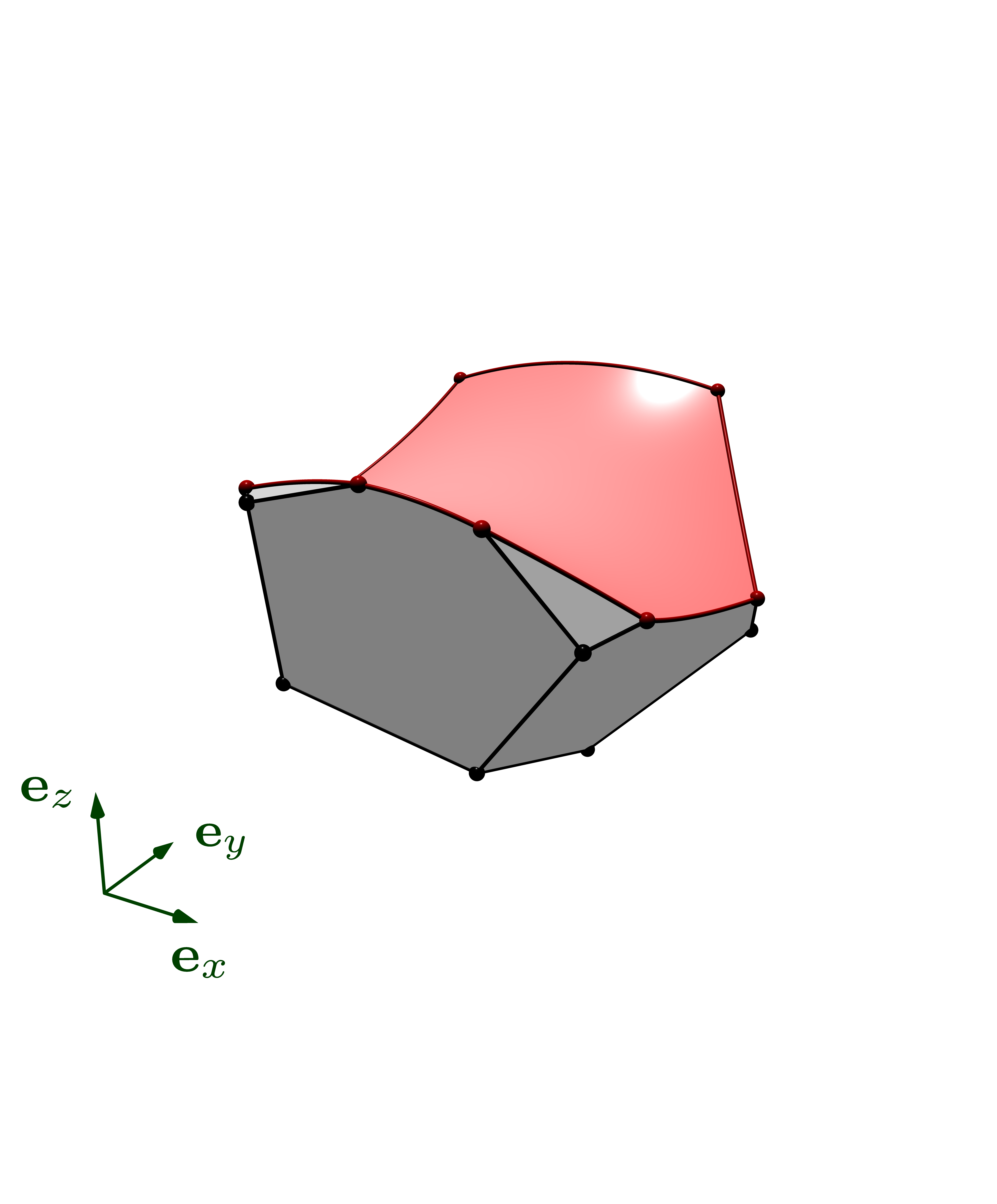}  \label{fig:clippedpoly}}\quad
    \subfloat[The faces of $\smash{\partial\hat{\mathcal{P}} =\cup_{i} \hat{\mathcal{F}}_i \cup \tilde{\mathcal{S}}}$.]{\adjincludegraphics[width=0.45\textwidth,trim=0 {0.17\height} 0 {0.25\height},clip=true]{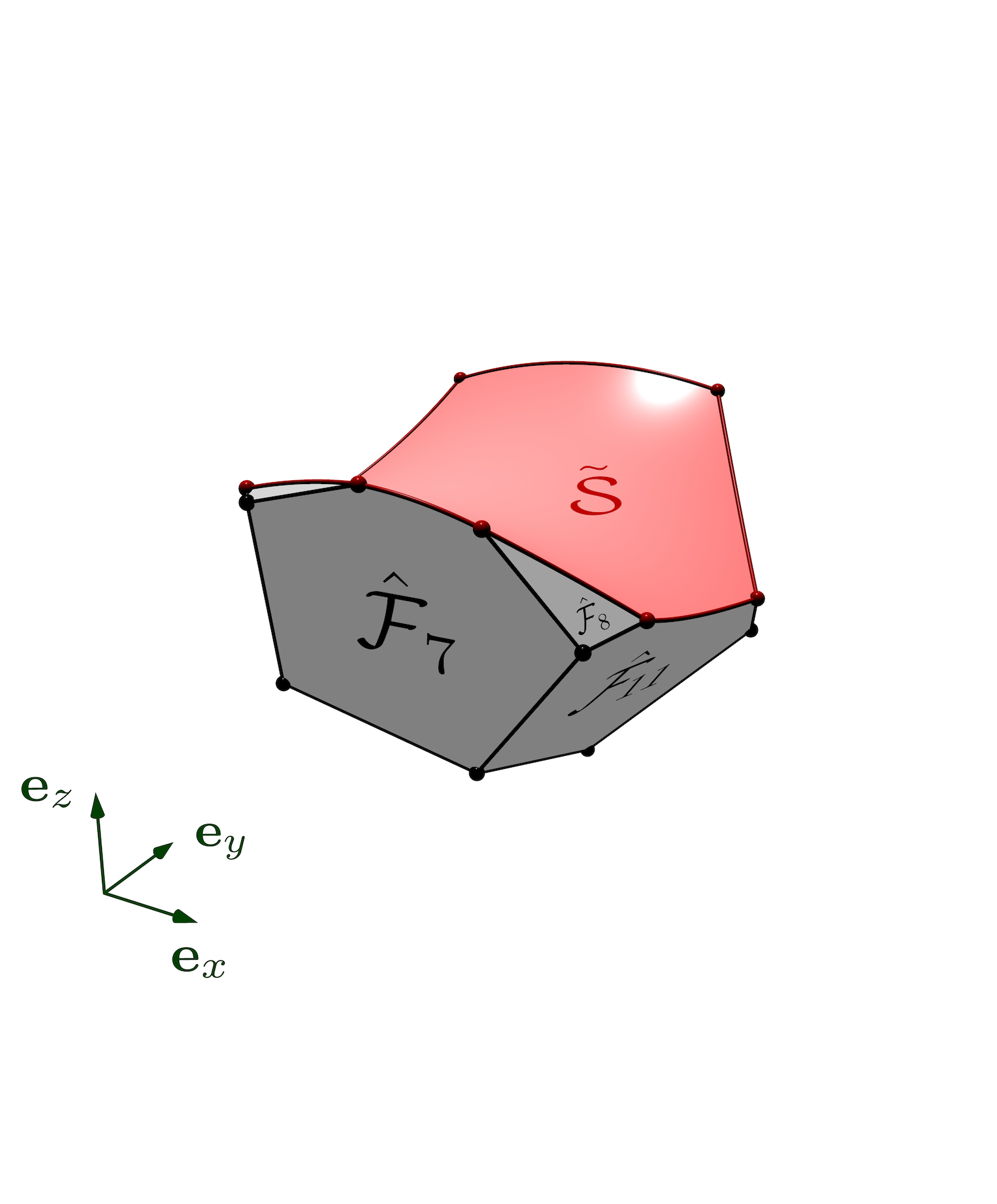}  \label{fig:clippedpoly_withfaces}}
    \caption{A polyhedron $\smash{\mathcal{P}\subset\mathbb{R}^3}$ intersected by the clipping region $\smash{\mathcal{Q}\subset\mathbb{R}^3}$ located below a paraboloid surface $\smash{\mathcal{S}}$.}
    \label{fig:intro_poly}
\end{figure}

Consider the two following subsets of $\smash{\mathbb{R}^3}$: 
\begin{itemize}\setlength\itemsep{0.2em}
    \item A polyhedron~$\smash{\mathcal{P}}$ delimited by $\smash{n_{\!\mathcal{F}}}$ planar polygonal faces $\smash{\mathcal{F}_i, i \in \{1,\ldots,n_{\!\mathcal{F}}\}}$, each equipped with an outward-pointing normal vector $\smash{\mathbf{n}_i}$ (e.g., see \cref{fig:poly,fig:poly_withfaces}). Neither $\smash{\mathcal{P}}$ nor its faces $\smash{\mathcal{F}_i}$ are required to be convex.
    \item The region $\smash{\mathcal{Q}}$, located on one side of a paraboloid~$\smash{\mathcal{S}}$ (e.g., see \cref{fig:poly_and_para}).
\end{itemize}
Without loss of generality, we assume to be working in a Cartesian coordinate system  
equipped with the orthonormal basis $\smash{(\mathbf{e}_x, \mathbf{e}_y,\mathbf{e}_z)}$, within which the position vector reads $\smash{\mathbf{x} = \begin{bmatrix}x & y & z\end{bmatrix}^\intercal}$, and where $\smash{\mathcal{Q}}$ and $\smash{\mathcal{S}}$ are implicitly defined as
\begin{align}
    \mathcal{Q} & = \{ \mathbf{x} \in \mathbb{R}^3 : \phi(\mathbf{x}) \le 0 \} \, , \label{def:halfspace} \\
    \mathcal{S} & = \{ \mathbf{x} \in \mathbb{R}^3 : \phi(\mathbf{x}) = 0 \} \, , \label{def:paraboloid}
\end{align}
with
\begin{equation}
    \phi : \left\{\begin{array}{lcl} 
        \mathbb{R}^3 & \!\to & \mathbb{R} \\
        \mathbf{x} & \!\mapsto & \alpha x^2 + \beta y^2 + z \,, 
    \end{array}\right. \quad (\alpha,\beta) \in \mathbb{R}^2 \, . \label{def:psi} 
\end{equation}
These assumptions do not restrict $\smash{\mathcal{Q}}$ and $\smash{\mathcal{S}}$ since, for any paraboloid-bounded clipping region in $\smash{\mathbb{R}^3}$, there exists a combination of rotations and translations of the canonical coordinate system resulting in such implicit definitions of $\smash{\mathcal{Q}}$ and $\smash{\mathcal{S}}$. 
For the sake of clarity and conciseness, we introduce the following notations:
\begin{itemize}\setlength\itemsep{0.2em}
    \item The subscript $\smash{\square_i}$ refers to a topological element or quantity related to the $i$th face of the polyhedron $\smash{\mathcal{P}}$.
    \item The superscript $\smash{\hat{\square}}$ implies an intersection with the clipping region $\smash{\mathcal{Q}}$, e.g., $\smash{\hat{\mathcal{P}} \equiv \mathcal{P}\cap\mathcal{Q}}$ or $\smash{\hat{\mathcal{F}}_i \equiv \mathcal{F}_i\cap\mathcal{Q}}$.
    \item The superscript $\smash{\tilde{\square}}$ implies an intersection with the polyhedron $\smash{\mathcal{P}}$, e.g., $\smash{\tilde{\mathcal{S}} \equiv \mathcal{S}\cap\mathcal{P}}$. This means that $\smash{\hat{\mathcal{P}} \equiv \tilde{\mathcal{Q}}}$.
\end{itemize}
As mentioned in \Cref{sec:introduction}, we are interested in calculating the zeroth and first moments of $\smash{\hat{\mathcal{P}}=\mathcal{P}\cap\mathcal{Q}}$ (e.g., see \cref{fig:clippedpoly}), i.e., its volume and volume-weighted centroid, given as
\begin{align}
	\mathrm{M}_0^{\hat{\mathcal{P}}} & = \int_{\hat{\mathcal{P}}} 1 \ \mathrm{d}\mathbf{x} \, , \quad \quad \text{and} \quad \quad \mathbf{M}^{\hat{\mathcal{P}}}_1 = \int_{\hat{\mathcal{P}}} \mathbf{x} \ \mathrm{d}\mathbf{x} \, .
\end{align}
In the remainder of this work, we shall refer to these quantities as ``the first moments'' or ``the moments'' of~$\smash{\hat{\mathcal{P}}}$, which we group into the vector
\begin{align}
	\boldsymbol{\mathcal{M}}^{\hat{\mathcal{P}}} & = \begin{bmatrix} \mathrm{M}_0^{\hat{\mathcal{P}}} \\ \mathbf{M}^{\hat{\mathcal{P}}}_1 \end{bmatrix} = \int_{\hat{\mathcal{P}}}  \boldsymbol{\Upsilon}(\mathbf{x}) \ \mathrm{d}\mathbf{x} \, , \label{eq:moments}
\end{align}
where
\begin{equation}
    \boldsymbol{\Upsilon} : \left\{\begin{array}{lcl} 
        \mathbb{R}^3 & \!\to & \mathbb{R}^4 \\
        \mathbf{x} & \!\mapsto & \begin{bmatrix}  1 & x & y & z \end{bmatrix}^\intercal \end{array}\right.  \, . \label{eq:defmonomials}
\end{equation}

\section{Moments derivation}
\label{sec:moments}
\begin{figure} \centering
    \subfloat[Projection of $\smash{\partial\hat{\mathcal{P}}}$ (where $\smash{\hat{\mathcal{P}}}$ is shown in \cref{fig:clippedpoly}) onto the $xy$-plane.]{\adjincludegraphics[width=0.45\textwidth,trim={0.03\width} {0.05\height} {0.2\width} {0.1\height},clip=true]{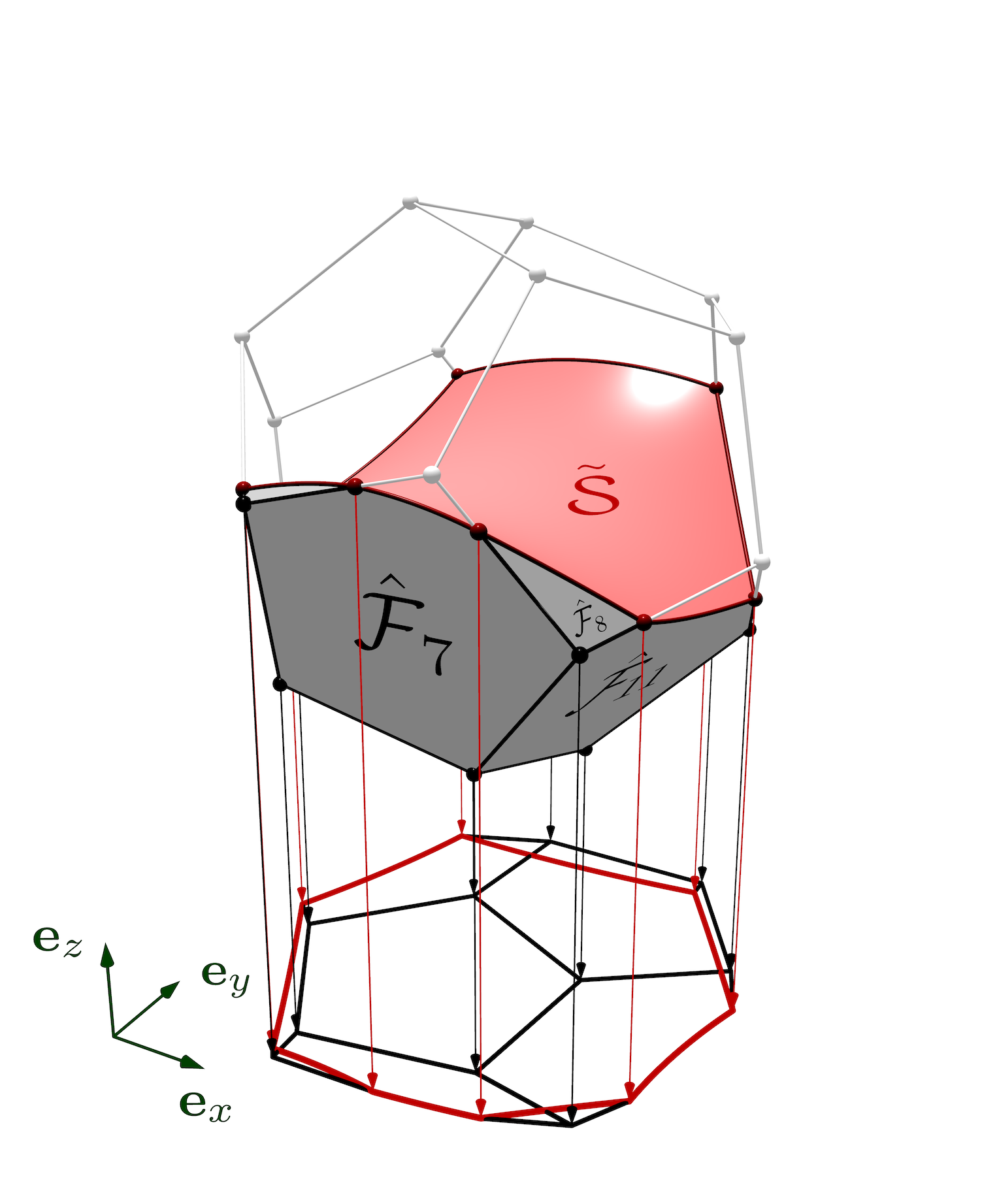} \label{fig:projection_3d}} \quad
    \subfloat[Projected clipped faces in the $xy$-plane.]{\includegraphics{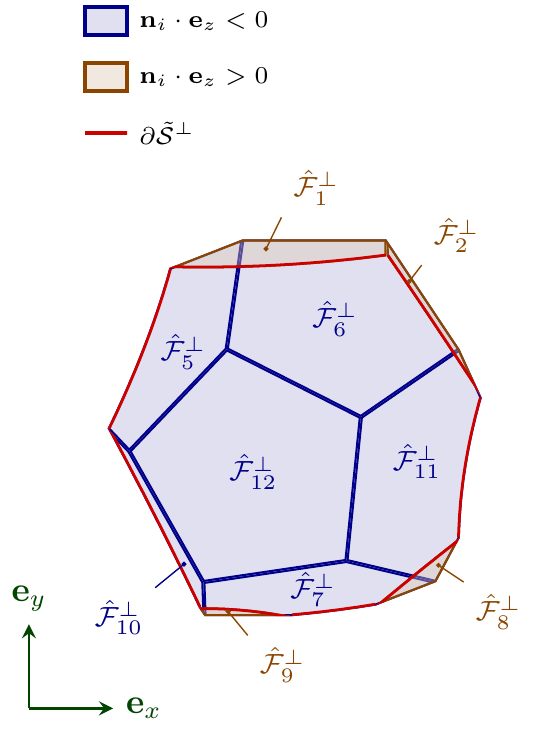} \label{fig:proj_viewxy}}\\
    \subfloat[Projected clipped face $\smash{\hat{\mathcal{F}}_{11}^\perp}$ and the splitting of its boundary into $\smash{n_{\partial\hat{\mathcal{F}}_{11}} = 6}$ parametrized arcs.]{\includegraphics{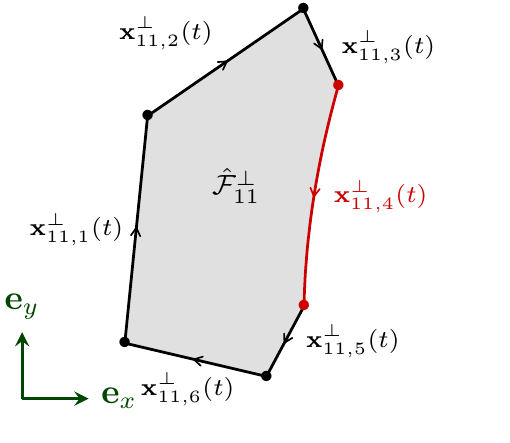} \label{fig:proj_f11}}  \quad\quad\quad
    \subfloat[Oriented closed curve on which $\smash{\boldsymbol{\mathcal{M}}^{\hat{\mathcal{P}}_1}_{11}}$ is integrated.]{\includegraphics{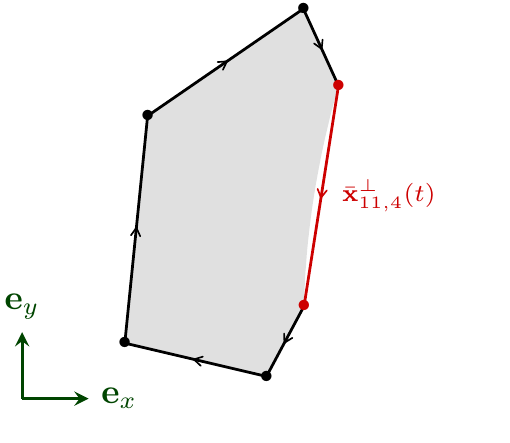} \label{fig:proj_f11_m1}}\\
    \subfloat[Oriented closed curve on which $\smash{\boldsymbol{\mathcal{M}}^{\hat{\mathcal{P}}_2}_{11,\bigtriangleup}}$ is integrated.]{\includegraphics{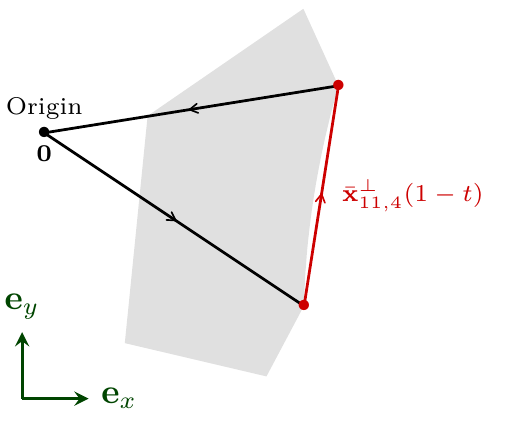} \label{fig:proj_f11_m2}}  \quad\quad\quad
    \subfloat[Oriented closed curve on which $\smash{\boldsymbol{\mathcal{M}}^{\hat{\mathcal{P}}_3}_{11}}$ is integrated.]{\includegraphics{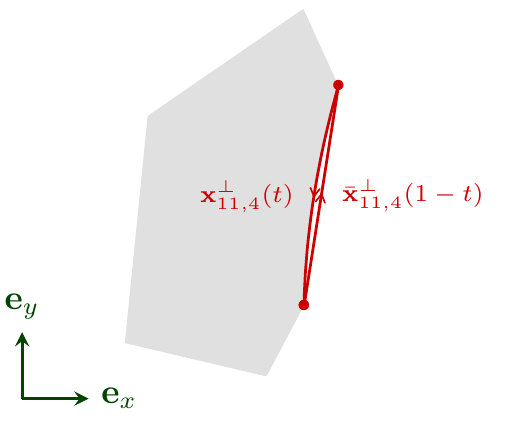} \label{fig:proj_f11_m3}}
    \caption{Illustration of the one-dimensional integration domains used for integrating the first, second, and third contributions to the moments.} \label{fig:projection}
\end{figure}

Using the divergence theorem, Eq.~\cref{eq:moments} can be rewritten as
\begin{equation}
	\boldsymbol{\mathcal{M}}^{\hat{\mathcal{P}}} = \int_{\hat{\mathcal{P}}} \nabla \cdot \left( \boldsymbol{\Phi}(\mathbf{x}) \otimes \mathbf{e}_z \right) \ \mathrm{d} \mathbf{x} 
    = \int_{\partial\hat{\mathcal{P}}} \boldsymbol{\Phi}(\mathbf{x}) (\mathbf{n}_{\partial\hat{\mathcal{P}}} \cdot \mathbf{e}_z) \ \mathrm{d} \mathbf{a}  \, , \label{eq:vol_after_div1}
\end{equation}
where $\smash{\mathrm{d}\mathbf{a}}$ is an infinitesimal surface element on $\smash{\partial\hat{\mathcal{P}} = \cup_{i} \hat{\mathcal{F}}_i \cup \tilde{\mathcal{S}}}$, the boundary of $\smash{\hat{\mathcal{P}}}$, $\smash{\mathbf{n}_{\partial\hat{\mathcal{P}}}}$ is the normal to $\smash{\partial\hat{\mathcal{P}}}$ pointing towards the outside of $\smash{\hat{\mathcal{P}}}$, and $\smash{\boldsymbol{\Phi}}$ is defined as
\begin{equation}
    \boldsymbol{\Phi} : \left\{\begin{array}{lcl} 
        \mathbb{R}^3 & \!\to & \mathbb{R}^4 \\
        \mathbf{x} & \!\mapsto & {\displaystyle \int_0^z} \boldsymbol{\Upsilon}(\mathbf{x}) \ \mathrm{d}z \ \ \ \left(= \begin{bmatrix} z & xz & yz & \frac{1}{2}z^2 \end{bmatrix}^\intercal \right) \end{array}\right.  \, .
\end{equation}
Eq.~\cref{eq:vol_after_div1} can be split into the following sum of integrals,
\begin{equation}
	\boldsymbol{\mathcal{M}}^{\hat{\mathcal{P}}} = \int_{\tilde{\mathcal{S}}}\boldsymbol{\Phi}(\mathbf{x}) (\mathbf{n}_{\tilde{\mathcal{S}}} \cdot \mathbf{e}_z) \ \mathrm{d}\mathbf{a} + \sum\limits_{i=1}^{n_\mathcal{F}} \int_{\hat{\mathcal{F}}_i} \boldsymbol{\Phi}(\mathbf{x}) (\mathbf{n}_{i} \cdot \mathbf{e}_z) \ \mathrm{d}\mathbf{a} \, , \label{eq:vol_after_div2}
\end{equation}
where $\smash{\tilde{\mathcal{S}} = \mathcal{S} \cap \mathcal{P}}$ is the portion of the paraboloid $\smash{\mathcal{S}}$ inside the polyhedron $\smash{\mathcal{P}}$ (e.g., see \cref{fig:poly_inter_surf}), $\smash{\mathbf{n}_{\tilde{\mathcal{S}}}}$ is the normal to $\smash{\tilde{\mathcal{S}}}$ pointing outwards of $\smash{\mathcal{Q}}$ (i.e., $\smash{\mathbf{n}_{\tilde{\mathcal{S}}} \cdot \mathbf{e}_z \ge 0}$), and $\smash{\hat{\mathcal{F}}_i = \mathcal{F}_i \cap\mathcal{Q}}$ is the portion of the face $\smash{\mathcal{F}_i}$ inside the clipping region $\smash{\mathcal{Q}}$ (e.g., see \cref{fig:clippedpoly_withfaces}).
Owing to the definitions of $\smash{\mathcal{Q}}$ and $\smash{\mathcal{S}}$, as given in Eqs.~\cref{def:paraboloid,def:halfspace}, the normal to $\smash{\mathcal{S}}$ pointing outwards of $\smash{\mathcal{Q}}$ reads as
\begin{equation}
    \mathbf{n}_{{\mathcal{S}}} = \frac{\nabla \phi}{\left\|\nabla \phi\right\|} \, ,
\end{equation}
yielding
\begin{equation}
	\boldsymbol{\mathcal{M}}^{\hat{\mathcal{P}}} = \int_{\tilde{\mathcal{S}}} \boldsymbol{\Phi}(\mathbf{x}) \left(\frac{\nabla \phi(\mathbf{x})\cdot \mathbf{e}_z}{\left\|\nabla \phi(\mathbf{x})\right\|}\right)  \ \mathrm{d}\mathbf{a} + \sum\limits_{i=1}^{n_\mathcal{F}} \int_{\hat{\mathcal{F}}_i} \boldsymbol{\Phi}(\mathbf{x}) (\mathbf{n}_{i} \cdot \mathbf{e}_z) \ \mathrm{d}\mathbf{a} \, . \label{eq:vol_after_div3}
\end{equation}
The surface $\smash{\tilde{\mathcal{S}}}$ can be expressed in the parametric form
\begin{equation}
    \tilde{\mathcal{S}} = \left\{ \begin{bmatrix} x \\ y \\ - \alpha x^2 - \beta y^2\end{bmatrix} , (x,y) \in \tilde{\mathcal{S}}^{\perp} \right\} \, , \label{def:paramS}
\end{equation}
with $\alpha$ and $\beta$ the coefficients introduced in Eq.~\cref{def:psi}, and with $\smash{\tilde{\mathcal{S}}^{\perp}}$ the projection of $\smash{\tilde{\mathcal{S}}}$ onto the $xy$-plane. Under the assumption that $\smash{\mathbf{n}_{i} \cdot \mathbf{e}_z \ne 0}$, each face $\smash{\hat{\mathcal{F}}_i}$, $\smash{i \in \{1,\ldots,n_{\!\mathcal{F}}\}}$, can also be expressed in the parametric form
\begin{equation}
    \hat{\mathcal{F}}_i = \left\{ \begin{bmatrix} x \\ y \\ \delta_i - \lambda_i x - \tau_i y \end{bmatrix}, (x,y) \in \hat{\mathcal{F}}_{i}^{\perp} \right\} \, ,  \label{def:paramFk}
\end{equation}
with 
\begin{align}
    \delta_i & = \frac{\mathbf{n}_i\cdot\mathbf{x}_{\mathcal{F}_i}}{\mathbf{n}_i \cdot \mathbf{e}_z} , \quad \text{for any } \mathbf{x}_{\mathcal{F}_i} \in \mathcal{F}_i \, , \label{eq:distance} \\
    \lambda_i & = \frac{\mathbf{n}_i \cdot \mathbf{e}_x}{\mathbf{n}_i \cdot \mathbf{e}_z} \, , \label{eq:lambda} \\
    \tau_i & = \frac{\mathbf{n}_i \cdot \mathbf{e}_y}{\mathbf{n}_i \cdot \mathbf{e}_z} \, ,\label{eq:tau} 
\end{align}
and with $\smash{\hat{\mathcal{F}}_{i}^{\perp}}$ the projection of $\smash{\hat{\mathcal{F}}_{i}}$ onto the $xy$-plane. These explicit parametrizations yield
\begin{equation}
	\boldsymbol{\mathcal{M}}^{\hat{\mathcal{P}}} = \int_{\tilde{\mathcal{S}}^{\perp}}  \boldsymbol{\Phi}_{\!\mathcal{S}}(x,y) \ \mathrm{d}\mathbf{a}^\perp + \sum\limits_{i=1}^{n_\mathcal{F}} \text{sign}(\mathbf{n}_{i} \cdot \mathbf{e}_z)\int_{\hat{\mathcal{F}}_{i}^{\perp}} \boldsymbol{\Phi}_{\!\mathcal{F}_i}(x,y) \ \mathrm{d}\mathbf{a}^\perp \, , \label{eq:vol_after_div4}
\end{equation}
where $\smash{\boldsymbol{\Phi}_{\!\mathcal{S}}}$ and $\smash{\boldsymbol{\Phi}_{\!\mathcal{F}_i}}$ are the function vectors
\begin{equation}
    \boldsymbol{\Phi}_{\!\mathcal{S}}(x,y)  = \begin{bmatrix} - \alpha x^2 - \beta y^2 \\  - \alpha x^3 - \beta x y^2 \\  - \alpha y x^2 - \beta y^3 \\  \frac{1}{2}\left(\alpha x^2 + \beta y^2\right)^2 \end{bmatrix} \, , \quad \quad
    \boldsymbol{\Phi}_{\!\mathcal{F}_i}(x,y) = \begin{bmatrix} \delta_i - \lambda_i x - \tau_i y \\ \delta_i x - \lambda_i x^2 - \tau_i x y \\ \delta_i y - \lambda_i x y - \tau_i y^2 \\ \frac{1}{2} \left(\delta_i - \lambda_i x - \tau_i y \right)^2 \end{bmatrix} \, , \label{eq:primitives}
\end{equation}
and $\smash{\mathrm{d}\mathbf{a}^\perp = \mathrm{d}x\mathrm{d}y}$ is an infinitesimal surface element on the $xy$-plane.
Note that in order to simplify Eq.~\cref{eq:vol_after_div3} into \cref{eq:vol_after_div4}, we have used the fact that $\smash{\left\|\nabla \phi(\mathbf{x})\right\|}$ is the determinant of the parametrization \cref{def:paramS} of $\smash{\tilde{\mathcal{S}}}$, and that $\smash{\left|\mathbf{n}_{i} \cdot \mathbf{e}_z\right|^{-1}}$ is the determinant of the parametrization \cref{def:paramFk} of $\smash{\hat{\mathcal{F}}_i}$. The projected integration domains $\smash{\tilde{\mathcal{S}}^\perp}$ and $\smash{\hat{\mathcal{F}}_i^\perp}$ corresponding to the configuration of \cref{fig:intro_poly} are illustrated in \cref{fig:projection}.\medskip

Eq.~\cref{eq:vol_after_div4} can be rewritten as
\begin{align}
	\boldsymbol{\mathcal{M}}^{\hat{\mathcal{P}}} & = \int_{\tilde{\mathcal{S}}^{\perp}} \nabla \cdot \left( \boldsymbol{\Psi}_{\mathcal{S}}(x,y) \otimes \mathbf{e}_x\right) \ \mathrm{d}\mathbf{a}^\perp \label{eq:divsecondtime}\\
    & \quad \quad \quad \quad \quad \quad + \sum\limits_{i=1}^{n_\mathcal{F}} \text{sign}(\mathbf{n}_{i} \cdot \mathbf{e}_z)\int_{\hat{\mathcal{F}}_{i}^{\perp}} \nabla \cdot \left(\boldsymbol{\Psi}_{\mathcal{F}_i}(x,y) \otimes \mathbf{e}_x\right)\ \mathrm{d}\mathbf{a}^\perp \, , \nonumber
\end{align}
where $\smash{\boldsymbol{\Psi}_{\mathcal{S}}}$ and $\smash{\boldsymbol{\Psi}_{\mathcal{F}_i}}$ are defined as
\begin{equation}
    \boldsymbol{\Psi}_{\mathcal{S}} : \left\{\begin{array}{lcl} 
        \mathbb{R}^2 & \!\to & \mathbb{R}^4 \\
        (x,y) & \!\mapsto & {\displaystyle \int_0^x} \boldsymbol{\Phi}_{\!\mathcal{S}}(x,y) \ \mathrm{d}x \end{array}\right.  \, . \label{eq:defpsiS}
\end{equation}
and $\smash{\forall i \in \{1,\ldots,n_{\!\mathcal{F}}\}}$,
\begin{equation}
    \boldsymbol{\Psi}_{\mathcal{F}_i} : \left\{\begin{array}{lcl} 
        \mathbb{R}^2 & \!\to & \mathbb{R}^4 \\
        (x,y) & \!\mapsto & {\displaystyle \int_0^x} \boldsymbol{\Phi}_{\!\mathcal{F}_i}(x,y) \ \mathrm{d}x \end{array}\right. \, . \label{eq:defpsiF}
\end{equation}
Note that the choice made here of integrating $\smash{\boldsymbol{\Phi}_{\!\mathcal{S}}}$ and $\smash{\boldsymbol{\Phi}_{\!\mathcal{F}_i}}$ with respect to $x$ is arbitrary, and that we could have equivalently integrated them with respect to $y$, requiring to replace $\smash{\mathbf{e}_x}$ by $\smash{\mathbf{e}_y}$ in Eq.~\cref{eq:divsecondtime}.
Using the divergence theorem once again, this gives
\begin{equation}
	\boldsymbol{\mathcal{M}}^{\hat{\mathcal{P}}} = \int_{\partial\tilde{\mathcal{S}}^{\perp}} \boldsymbol{\Psi}_{\mathcal{S}}(x,y) \left(\mathbf{n}_{\partial\tilde{\mathcal{S}}^{\perp}} \cdot \mathbf{e}_x\right) \ \mathrm{d}l + \sum\limits_{i=1}^{n_\mathcal{F}} \int_{\partial\hat{\mathcal{F}}_{i}^{\perp}} \boldsymbol{\Psi}_{\mathcal{F}_i}(x,y) (\mathbf{n}_{\partial\hat{\mathcal{F}}_{i}^{\perp}} \cdot \mathbf{e}_x) \ \mathrm{d}l \, , \label{eq:vol_after_div6}
\end{equation}
where $\smash{\mathrm{d}l}$ is an infinitesimal line element on the integration domains $\smash{\partial\tilde{\mathcal{S}}^{\perp}}$ and $\smash{\partial\hat{\mathcal{F}}_{i}^{\perp}}$, which are the boundaries of the projections of the faces of $\smash{\hat{\mathcal{P}}}$ onto the $xy$-plane. As such, they consist of closed curves in the $xy$-plane, that are successions of conic section arcs and/or line segments (e.g., see \cref{fig:proj_viewxy,fig:proj_f11}). Note that the term ``$\smash{\text{sign}(\mathbf{n}_{i} \cdot \mathbf{e}_z)}$'' is now implicitly accounted for, as the closed curves $\smash{\partial\hat{\mathcal{F}}_i}$ (and therefore their projection onto the $xy$-plane, $\smash{\partial\hat{\mathcal{F}}_i^\perp}$) are oriented so as to produce a normal vector pointing towards the outside of $\smash{\hat{\mathcal{P}}}$. It should also be noted that the integration domains $\smash{\partial\tilde{\mathcal{S}}^{\perp}}$ and $\smash{\partial\hat{\mathcal{F}}_{i}^{\perp}}$ do not necessarily consist of one closed curve each -- they may each be the union of several non-intersecting oriented closed curves.\medskip

Let us assume that a parametrization
\begin{equation}
    \mathbf{x}_{i,j}(t) = x_{i,j}(t) \mathbf{e}_x + y_{i,j}(t) \mathbf{e}_y + z_{i,j}(t) \mathbf{e}_z, \quad t \in [0,1] , \quad j \in \{1,\ldots,n_{\partial\hat{\mathcal{F}_i}}\} \, , \label{eq:param_arc}
\end{equation}
is known for each of the $\smash{n_{\partial\hat{\mathcal{F}_i}}}$ arcs of $\smash{\partial\hat{\mathcal{F}}_{i}}$, where the functions $\smash{x_{i,j}}$, $\smash{y_{i,j}}$, and $\smash{z_{i,j}}$ belong to $\mathcal{C}^1([0,1])$. Moreover, let us note that each parametrized conic section arc belonging to $\smash{\partial\tilde{\mathcal{S}}}$ is necessarily present in one and only one of the clipped face boundaries $\smash{\partial\hat{\mathcal{F}}_i}$, where it is traversed in the opposite direction for integrating $\smash{\boldsymbol{\Psi}_{\!\mathcal{F}_i}}$. Eq.~\cref{eq:vol_after_div6} can then be written as
\begin{align}
	\boldsymbol{\mathcal{M}}^{\hat{\mathcal{P}}} = \sum\limits_{i=1}^{n_{\!\mathcal{F_{\phantom{k}\!\!}}}} \sum\limits_{j = 1}^{n_{\partial\hat{\mathcal{F}}_i}}  \int_{0}^{1} \left(\boldsymbol{\Psi}_{\mathcal{F}_i}(x_{i,j}(t),y_{i,j}(t)) - 1^{\partial\tilde{\mathcal{S}}}_{i,j} \, \boldsymbol{\Psi}_{\mathcal{S}}(x_{i,j}(t),y_{i,j}(t)) \right) y^\prime_{i,j}(t) \ \mathrm{d}t \, , \label{eq:vol_after_div7}
\end{align}
where
\begin{align}
    1^{\partial\tilde{\mathcal{S}}}_{i,j} & = \left\{ \begin{array}{ll} 1 & \text{if the $j$th arc of $\partial\hat{\mathcal{F}}_i$ also belongs to $\partial\tilde{\mathcal{S}}$} \\ 0 & \text{otherwise} \end{array}\right. \, . 
\end{align}
Note that in Eq.~\cref{eq:vol_after_div7} and in the remainder of this manuscript, the superscript $\smash{\square^\prime}$ indicates that a function has been differentiated with respect to its unique variable.
A closed-form expression can be derived for the integral in Eq.~\cref{eq:vol_after_div7}, however, its use for numerically calculating the moments is undesirable for two main reasons: 
\begin{enumerate} 
    \item the expression contains many terms, rendering its numerical calculation expensive;
    \item the expression depends on $\smash{\delta_i}$, $\smash{\lambda_i}$, and $\smash{\tau_i}$, which all tend towards infinity as $\mathbf{n}_i \cdot \mathbf{e}_z$ tends towards zero, leading to large round-off errors in the context of floating-point arithmetics.
\end{enumerate}
Instead, the introduction of a twin parametrization of the arcs of $\smash{\partial\partial\hat{\mathcal{P}}}$ and the judicious splitting of the integral in Eq.~\cref{eq:vol_after_div7} can both reduce the complexity of its closed-form expression and remove its direct dependency on the potentially singular coefficients $\smash{\delta_i}$, $\smash{\lambda_i}$, and $\smash{\tau_i}$. Let us then introduce the parametrization
\begin{equation}
    \bar{\mathbf{x}}_{i,j}(t) = \bar{x}_{i,j}(t) \mathbf{e}_x + \bar{y}_{i,j}(t) \mathbf{e}_y + \bar{z}_{i,j}(t) \mathbf{e}_z, \quad t \in [0,1] , \quad j \in \{1,\ldots,n_{\partial\hat{\mathcal{F}_i}}\} \, , \label{eq:param_line}
\end{equation}
which links $\smash{{\mathbf{x}}_{i,j}(0)}$ to $\smash{{\mathbf{x}}_{i,j}(1)}$ by a straight line. For the sake of conciseness, we shall refer to these two points as $\smash{{\mathbf{x}}_{i,j}(0) = {\mathbf{x}}_{i,j,0}}$ and $\smash{{\mathbf{x}}_{i,j}(1) = {\mathbf{x}}_{i,j,1}}$ in the remainder of this work. This twin parametrization of each arc is simply given as
\begin{equation}
    \bar{\mathbf{x}}_{i,j}(t) = (1-t) {\mathbf{x}}_{i,j,0}  + t {\mathbf{x}}_{i,j,1}, \quad t \in [0,1] , \quad j \in \{1,\ldots,n_{\partial\hat{\mathcal{F}_i}}\} \, . \label{eq:param_arc_straight}
\end{equation}
If the $j$th arc of $\smash{\partial\hat{\mathcal{F}}_i}$ does not belong to $\smash{\tilde{\mathcal{S}}}$, then $\smash{x_{i,j} \equiv \bar{x}_{i,j}}$, $\smash{y_{i,j} \equiv \bar{y}_{i,j}}$ and $\smash{z_{i,j} \equiv \bar{z}_{i,j}}$. We can then re-organize Eq.~\cref{eq:vol_after_div7} as
\begin{normalsize}\begin{align}
    \boldsymbol{\mathcal{M}}^{\hat{\mathcal{P}}} & = \sum\limits_{i=1}^{n_{\!\mathcal{F_{\phantom{k}\!\!}}}} \sum\limits_{j = 1}^{n_{\partial\hat{\mathcal{F}}_i}} \left[ \; \int_{0}^{1}  \boldsymbol{\Psi}_{\mathcal{F}_i}(\bar{x}_{i,j}(t),\bar{y}_{i,j}(t))\bar{y}^\prime_{i,j}(t)  \ \mathrm{d}t  \right. \label{eq:vol_after_div8} \\
    & \quad\quad\quad\quad\quad\quad - 1^{\partial\tilde{\mathcal{S}}}_{i,j} \int_{0}^{1} \boldsymbol{\Psi}_{\mathcal{S}}(\bar{x}_{i,j}(t),\bar{y}_{i,j}(t))\bar{y}^\prime_{i,j}(t) \ \mathrm{d}t  \nonumber\\
    & \quad\quad\quad\quad\quad\quad - 1^{\partial\tilde{\mathcal{S}}}_{i,j} \int_{0}^{1} \left( \phantom{\frac{1}{1}\!\!\!\!}\boldsymbol{\Psi}_{\mathcal{F}_i}(\bar{x}_{i,j}(t),\bar{y}_{i,j}(t)) - \boldsymbol{\Psi}_{\mathcal{S}}(\bar{x}_{i,j}(t),\bar{y}_{i,j}(t)) \right) \bar{y}^\prime_{i,j}(t) \nonumber \\
    & \quad\quad\quad\quad\quad\quad  \left. \phantom{\int_{0}^{1}} + \left(\phantom{\frac{1}{1}\!\!\!\!}\boldsymbol{\Psi}_{\mathcal{S}}(x_{i,j}(t),y_{i,j}(t)) - \boldsymbol{\Psi}_{\mathcal{F}_i}(x_{i,j}(t),y_{i,j}(t)) \right)y^\prime_{i,j}(t) \ \mathrm{d}t \; \right]\, , \nonumber
\end{align}\end{normalsize}The moments can thus be described as the sum of three distinct contributions, i.e.,
 \begin{equation}
    \boldsymbol{\mathcal{M}}^{\hat{\mathcal{P}}} = \boldsymbol{\mathcal{M}}^{\hat{\mathcal{P}}_1} + \boldsymbol{\mathcal{M}}^{\hat{\mathcal{P}}_2} + \boldsymbol{\mathcal{M}}^{\hat{\mathcal{P}}_3} = \sum\limits_{i=1}^{n_{\!\mathcal{F_{\phantom{k}\!\!}}}} \boldsymbol{\mathcal{M}}^{\hat{\mathcal{P}}_1}_{i} + \sum\limits_{i=1}^{n_{\!\mathcal{F_{\phantom{k}\!\!}}}} \boldsymbol{\mathcal{M}}^{\hat{\mathcal{P}}_2}_{i} + \sum\limits_{i=1}^{n_{\!\mathcal{F_{\phantom{k}\!\!}}}} \boldsymbol{\mathcal{M}}^{\hat{\mathcal{P}}_3}_{i} \, ,
\end{equation}
where
\begin{align}
    \boldsymbol{\mathcal{M}}^{\hat{\mathcal{P}}_1}_{i} &= \int_{0}^{1} \sum\limits_{j = 1}^{n_{\partial\hat{\mathcal{F}}_i}} \boldsymbol{\Psi}_{\mathcal{F}_i}(\bar{x}_{i,j}(t),\bar{y}_{i,j}(t))\bar{y}^\prime_{i,j}(t) \ \mathrm{d}t \, , \label{eq:defM1}\\
    \boldsymbol{\mathcal{M}}^{\hat{\mathcal{P}}_2}_{i} & = \int_{0}^{1} \sum\limits_{j = 1}^{n_{\partial\hat{\mathcal{F}}_i}} -1^{\partial\tilde{\mathcal{S}}}_{i,j} \boldsymbol{\Psi}_{\mathcal{S}}(\bar{x}_{i,j}(t),\bar{y}_{i,j}(t))\bar{y}^\prime_{i,j}(t) \ \mathrm{d}t \, , \label{def:m2k} \\
    \boldsymbol{\mathcal{M}}^{\hat{\mathcal{P}}_3}_{i} &= \int_{0}^{1} \sum\limits_{j = 1}^{n_{\partial\hat{\mathcal{F}}_i}} 1^{\partial\tilde{\mathcal{S}}}_{i,j} \left(\phantom{\frac{1}{1}\!\!\!\!\!\!}\left( \boldsymbol{\Psi}_{\mathcal{S}}(\bar{x}_{i,j}(t),\bar{y}_{i,j}(t)) - \boldsymbol{\Psi}_{\mathcal{F}_i}(\bar{x}_{i,j}(t),\bar{y}_{i,j}(t)) \right)\bar{y}^\prime_{i,j}(t) \right. \label{eq:defM3}\\
    & \quad\quad\quad\quad\quad\quad \left. + \left(\boldsymbol{\Psi}_{\mathcal{F}_i}(x_{i,j}(t),y_{i,j}(t)) - \boldsymbol{\Psi}_{\mathcal{S}}(x_{i,j}(t),y_{i,j}(t)) \right)y^\prime_{i,j}(t) \phantom{\frac{1}{1}\!\!\!\!\!\!}\right) \ \mathrm{d}t \, . \nonumber 
\end{align}

The contributions $\smash{\boldsymbol{\mathcal{M}}^{\hat{\mathcal{P}}_1}}$ and $\smash{\boldsymbol{\mathcal{M}}^{\hat{\mathcal{P}}_2}}$ require the integration of the paraboloid and plane primitives over straight lines only (e.g., see \cref{fig:proj_f11_m1,fig:proj_f11_m2}), hence are straightforward to derive. The contribution $\smash{\boldsymbol{\mathcal{M}}^{\hat{\mathcal{P}}_3}}$, on the other hand, requires the parametrization of the conic section arcs in $\smash{\partial\tilde{\mathcal{S}}}$ (e.g., see \cref{fig:proj_f11_m3}). It should also be noted that each arc of the clipped faces $\smash{\partial\hat{\mathcal{F}}_i}$ contributes to $\smash{\boldsymbol{\mathcal{M}}^{\hat{\mathcal{P}}_1}}$, whereas only the conic section arcs in those faces (originating from the intersection of $\smash{\partial\mathcal{P}}$ with $\smash{\mathcal{S}}$) contribute to $\smash{\boldsymbol{\mathcal{M}}^{\hat{\mathcal{P}}_2}}$ and $\smash{\boldsymbol{\mathcal{M}}^{\hat{\mathcal{P}}_3}}$, owing to the presence of the coefficient~$\smash{1^{\partial\tilde{\mathcal{S}}}_{i,j}}$.

\subsection{First term: $\boldsymbol{\mathcal{M}}^{\hat{\mathcal{P}}_1}$}\label{sec:firstcontri}
Let us be reminded that the boundary of each clipped face, $\smash{\partial\hat{\mathcal{F}}_i}$, is a succession of $\smash{n_{\partial\hat{\mathcal{F}}_i}}$ conic section arcs and/or straight line segments that each link a start point $\smash{{\mathbf{x}}_{i,j,0}}$ to an end point $\smash{{\mathbf{x}}_{i,j,1}}$. Now recall that we aim to derive expressions that are free of the coefficents $\smash{\delta_i}$, $\smash{\lambda_i}$, and $\smash{\tau_i}$, so as to avoid round-off errors in the numerical calculation of the moments. To do so, we assign to each face~$\smash{\hat{\mathcal{F}}_i}$ a reference point $\smash{\mathbf{x}_{i,\text{ref}}}$ whose only requirement is to belong to the plane containing~$\smash{\hat{\mathcal{F}}_i}$, e.g., $\smash{\mathbf{x}_{i,\text{ref}} = {\mathbf{x}}_{i,1,0}}$. For each arc of each clipped face, rather than integrating on the straight line linking $\smash{{\mathbf{x}}_{i,j,0}}$ to $\smash{{\mathbf{x}}_{i,j,1}}$, we integrate instead on the oriented triangle $\smash{T^{(1)}_{i,j} = \left({\mathbf{x}}_{i,j,0},{\mathbf{x}}_{i,j,1},{\mathbf{x}}_{i,\text{ref}}\right)}$. Since $\smash{\partial\hat{\mathcal{F}}_i}$ is the union of closed curves, the start point of each of its constituting arcs is necessarily the end point of another arc, and the sum of all these triangle integrals is equal to the sum of the straight arc integrals. In other terms, substituting Eq.~\cref{eq:param_arc_straight} into Eq.~\cref{eq:defM1}, the latter can be rewritten as
\begin{small}\begin{align}
    \boldsymbol{\mathcal{M}}^{\hat{\mathcal{P}}_1}_{i} & = {\int_{0}^{1}} \sum\limits_{j = 1}^{n_{\partial\hat{\mathcal{F}}_i}}  \boldsymbol{\Psi}_{\mathcal{F}_i}((1-t){x}_{i,j,0}+t{x}_{i,j,1},(1-t){y}_{i,j,0}+t{y}_{i,j,1})({y}_{i,j,1}-{y}_{i,j,0}) \ \mathrm{d}t \, , \\
      & = {\int_{0}^{1}} \sum\limits_{j = 1}^{n_{\partial\hat{\mathcal{F}}_i}} \left(\!\!\begin{array}{lll}
        \boldsymbol{\Psi}_{\mathcal{F}_i}((1-t){x}_{i,j,0}+t{x}_{i,j,1},(1-t){y}_{i,j,0}+t{y}_{i,j,1})({y}_{i,j,1}-{y}_{i,j,0}) \\
     \ + \boldsymbol{\Psi}_{\mathcal{F}_i}((1-t){x}_{i,j,1}+t{x}_{i,\mathrm{ref}},(1-t){y}_{i,j,1}+t{y}_{i,\mathrm{ref}})({y}_{i,\mathrm{ref}}-{y}_{i,j,1}) \\
     \ \ \ + \boldsymbol{\Psi}_{\mathcal{F}_i}((1-t){x}_{i,\mathrm{ref}}+t{x}_{i,j,0},(1-t){y}_{i,\mathrm{ref}}+t{y}_{i,j,0})({y}_{i,j,0}-{y}_{i,\mathrm{ref}})\end{array}\!\!\right) \mathrm{d}t \, . \nonumber 
\end{align}\end{small}Making use of Eqs.~\cref{eq:primitives,eq:defpsiF}, and of the fact that
\begin{align}
    {z}_{i,\text{ref}} & =  \delta_i - \lambda_i {x}_{i,\text{ref}} - \tau_i {y}_{i,\text{ref}}  \\
    {z}_{i,j} & = \delta_i - \lambda_i {x}_{i,j} - \tau_i {y}_{i,j} \ , \quad \forall j \in \{1,\ldots,n_{\partial\hat{\mathcal{F}_i}}\}
\end{align}
it follows that
\begin{align}
    \boldsymbol{\mathcal{M}}^{\hat{\mathcal{P}}_1}_i & = \sum\limits_{j = 1}^{n_{\partial\hat{\mathcal{F}}_i}} \mathcal{A}\left({\mathbf{x}}_{i,j,0},{\mathbf{x}}_{i,j,1},{\mathbf{x}}_{i,\text{ref}}\right) \boldsymbol{\mathcal{B}}^{(1)}\left({\mathbf{x}}_{i,j,0},{\mathbf{x}}_{i,j,1},{\mathbf{x}}_{i,\text{ref}}\right) \, ,
\end{align}
where $\mathcal{A}$ is the operator for calculating the signed projected area of a triangle from the knowledge of its three corners, i.e.,
\begin{equation}
    \mathcal{A} : \left\{ \begin{array}{ccl} \mathbb{R}^3\times\mathbb{R}^3\times\mathbb{R}^3 & \to & \mathbb{R} \\ \left(\mathbf{x}_{a},\mathbf{x}_{b},\mathbf{x}_{c}\right) & \mapsto & \frac{1}{2} \left(x_a(y_b-y_c) + x_b(y_c-y_a) + x_c(y_a-y_b) \right) \end{array}\right. \, , \label{eq:triangle_area}
\end{equation}
and $\smash{\boldsymbol{\mathcal{B}}^{(1)} : \mathbb{R}^3\times\mathbb{R}^3\times\mathbb{R}^3 \to \mathbb{R}^4}$ reads as
\begin{equation}
    \boldsymbol{\mathcal{B}}^{(1)}\left(\mathbf{x}_{a},\mathbf{x}_{b},\mathbf{x}_{c}\right) = \frac{1}{12}
\begin{bmatrix} 4 \left( z_{a} + z_{b} + z_{c} \right) \\  
    \left( z_{a} + z_{b} + z_{c} \right) \left( x_{a} + x_{b} + x_{c} \right) + x_{a} z_{a} + x_{b} z_{b} + x_{c} z_{c}
     \\
     \left( z_{a} + z_{b} + z_{c} \right) \left( y_{a} + y_{b} + y_{c} \right) + y_{a} z_{a} + y_{b} z_{b} + y_{c} z_{c}
    \\ z_{a}^2 + z_{b}^2 +
  z_{c}^2 + z_{b} z_{c} +
  z_{a} z_{b} + z_{a} z_{c}\end{bmatrix} \, .
\end{equation}

\subsection{Second term: $\boldsymbol{\mathcal{M}}^{\hat{\mathcal{P}}_2}$}
For computing $\smash{\boldsymbol{\mathcal{M}}^{\hat{\mathcal{P}}_2}}$, similarly as in \Cref{sec:firstcontri}, we choose a reference point $\smash{\mathbf{x}_{\mathcal{S}}}$ belonging to $\smash{\mathcal{S}}$. An obvious choice for this reference point is the origin of our coordinate system, i.e., $\smash{\mathbf{x}_{\mathcal{S}}=\boldsymbol{0}}$. For each conic section arc of each clipped face, rather than integrating on the straight line linking $\smash{{\mathbf{x}}_{i,j,0}}$ to $\smash{{\mathbf{x}}_{i,j,1}}$, we integrate instead on the oriented triangle $\smash{T^{(2)}_{i,j} = \left({\mathbf{x}}_{i,j,0},{\mathbf{x}}_{i,j,1},\mathbf{x}_{\mathcal{S}}\right)}$.
For each $i$th face of $\smash{\mathcal{P}}$, this yields the moment contribution\begin{small}\begin{align}
    \boldsymbol{\mathcal{M}}^{\hat{\mathcal{P}}_2}_{i,\bigtriangleup} & = 
        {\int_{0}^{1}} \sum\limits_{j = 1}^{n_{\partial\hat{\mathcal{F}}_i}} -1^{\partial\tilde{\mathcal{S}}}_{i,j}\left(\!\!\begin{array}{lll}
            \boldsymbol{\Psi}_{\mathcal{S}}((1-t){x}_{i,j,0}+t{x}_{i,j,1},(1-t){y}_{i,j,0}+t{y}_{i,j,1})({y}_{i,j,1}-{y}_{i,j,0}) \\
         \ + \boldsymbol{\Psi}_{\mathcal{S}}((1-t){x}_{i,j,1}+t{x}_{\mathcal{S}},(1-t){y}_{i,j,1}+t{y}_{\mathcal{S}})({y}_{\mathcal{S}}-{y}_{i,j,1}) \\
         \ \ \ + \boldsymbol{\Psi}_{\mathcal{S}}((1-t){x}_{\mathcal{S}}+t{x}_{i,j,0},(1-t){y}_{\mathcal{S}}+t{y}_{i,j,0})({y}_{i,j,0}-{y}_{\mathcal{S}})\end{array}\!\!\right) \mathrm{d}t \, ,
\end{align}\end{small}Using Eqs.~\cref{eq:primitives,eq:defpsiS}, substituting $\smash{\mathbf{x}_{\mathcal{S}}=\boldsymbol{0}}$, and making use of the fact that
\begin{align}
    1^{\partial\tilde{\mathcal{S}}}_{i,j} \ne 0 & \Leftrightarrow \left\{\begin{array}{l} {z}_{i,j,0} = - \alpha {x}_{i,j,0}^2 - \beta {y}_{i,j,0}^2 \\  {z}_{i,j,1} = - \alpha {x}_{i,j,1}^2 - \beta {y}_{i,j,1}^2 \end{array}\right.  \ , \quad \forall j \in \{1,\ldots,n_{\partial\hat{\mathcal{F}_i}}\}
\end{align}
it follows that
\begin{align}
    \boldsymbol{\mathcal{M}}^{\hat{\mathcal{P}}_2}_{i,\bigtriangleup} & = \sum\limits_{j = 1}^{n_{\partial\hat{\mathcal{F}}_i}} 1^{\partial\tilde{\mathcal{S}}}_{i,j} \mathcal{A}\left({\mathbf{x}}_{i,j,0},{\mathbf{x}}_{i,j,1},\boldsymbol{0}\right) \boldsymbol{\mathcal{B}}^{(2)}\left({\mathbf{x}}_{i,j,0},{\mathbf{x}}_{i,j,1}\right) \, ,
\end{align}
where $\smash{\mathcal{A}}$ is the signed projected triangle area operator defined in Eq.~\cref{eq:triangle_area} and the operator $\smash{\boldsymbol{\mathcal{B}}^{(2)} : \mathbb{R}^3\times\mathbb{R}^3 \to \mathbb{R}^4}$ reads as
\begin{equation}
    \boldsymbol{\mathcal{B}}^{(2)}\left(\mathbf{x}_{a},\mathbf{x}_{b}\right) = {\small\frac{1}{90}
\begin{bmatrix}[1.1] 15 \left( \alpha x_{a} x_{b} + \beta y_{a} y_{b} - z_{a} - z_{b} \right) \\  
    6 \beta (y_{a} - y_{b})
           (x_{b} y_{a} - x_{a} y_{b}) +
       9 (x_{a} + x_{b}) (z_{a} + z_{b})
     \\
     6 \alpha (x_{a} - x_{b})
           (y_{b} x_{a} - y_{a} x_{b}) +
       9 (y_{a} + y_{b}) (z_{a} + z_{b})
    \\ 2 \alpha \beta
    (x_{b} y_{a} - x_{a} y_{b})^2 +
   3 (z_{a} + z_{b}) (\alpha x_{a} x_{b} + \beta y_{a} y_{b} - z_{a} - z_{b})
   + 3 z_{a}  z_{b}
    \end{bmatrix}} \, .
\end{equation}
Note that $\smash{\boldsymbol{\mathcal{M}}^{\hat{\mathcal{P}}_2}_{i,\bigtriangleup} \ne \boldsymbol{\mathcal{M}}^{\hat{\mathcal{P}}_2}_{i}}$, owing to our choice of integrating over triangles rather than the individual arcs, however their sum over all faces is equal, yielding
\begin{equation}
    \boldsymbol{\mathcal{M}}^{\hat{\mathcal{P}}_2} = \sum\limits_{i=1}^{n_{\!\mathcal{F_{\phantom{k}\!\!}}}} \boldsymbol{\mathcal{M}}^{\hat{\mathcal{P}}_2}_{i,\bigtriangleup} \, .
\end{equation}

\subsection{Third term: $\boldsymbol{\mathcal{M}}^{\hat{\mathcal{P}}_3}$}
\begin{figure}[h] \centering
\includegraphics{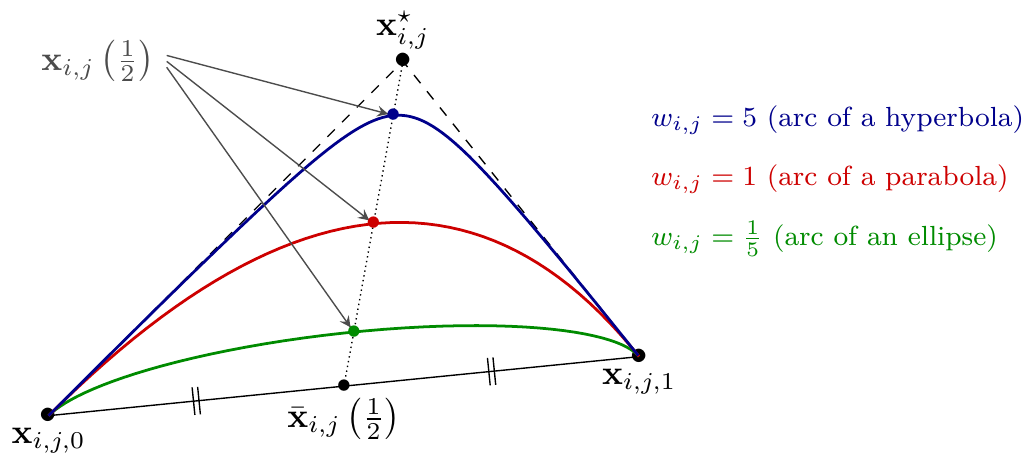}
\caption{Three rational B\'ezier arcs linking a start point $\smash{\mathbf{x}_{i,j,0}}$ to an end point $\smash{\mathbf{x}_{i,j,1}}$. Also shown are the control point $\smash{\mathbf{x}^{\star}_{i,j}}$ and weight $\smash{w_{i,j}}$, along with the points $\smash{\bar{\mathbf{x}}_{i,j}\left(\tfrac{1}{2}\right)}$ and $\smash{{\mathbf{x}}_{i,j}\left(\tfrac{1}{2}\right)}$ that can be used for determining $\smash{w_{i,j}}$.}
\label{fig:bezierarc}
\end{figure}

To derive a closed-form expression for $\smash{\boldsymbol{\mathcal{M}}^{\hat{\mathcal{P}}_3}}$, a parametrization of the conic section arcs in $\smash{\partial\tilde{\mathcal{S}}}$ must be provided. For the elliptic and hyperbolic cases, traditional parametrizations using trigonometric functions are obvious choices, however they can yield significant round-off errors due to very large values of their constitutive parameters (e.g., the semi-major and semi-minor axes). 
\begin{figure} \centering
\includegraphics{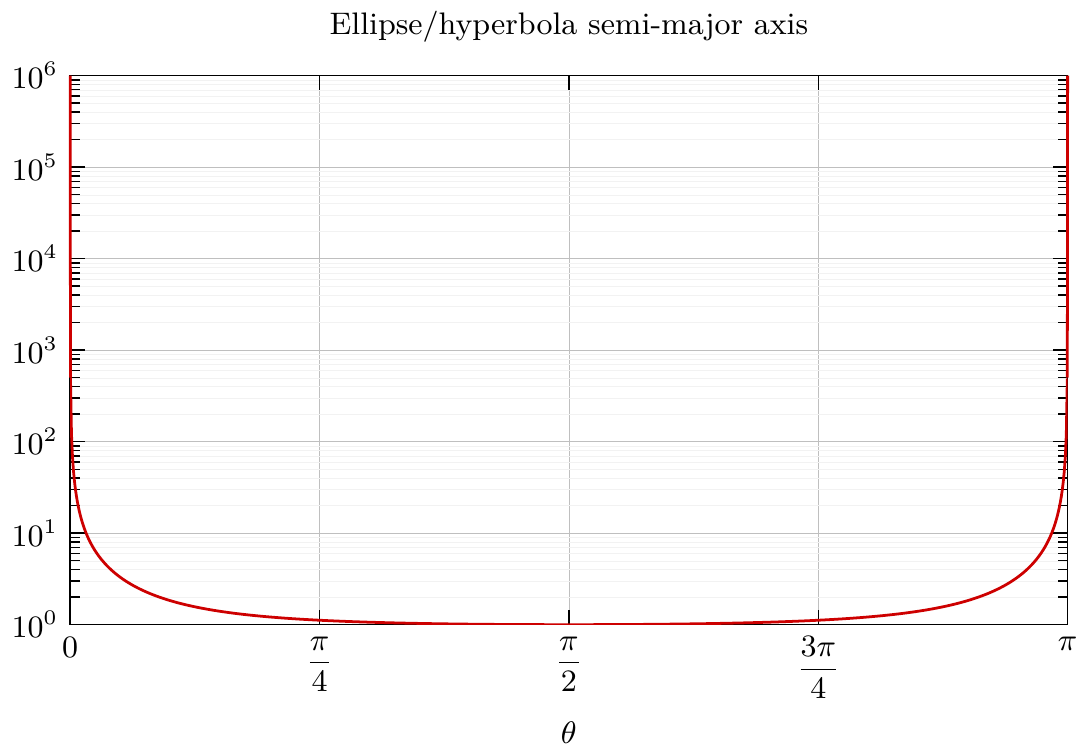}
\caption{Semi-major axis of the conic section generated by the intersection of the paraboloid $\mathcal{S}$, as defined in Eq.~\cref{def:paraboloid}, when $\smash{\alpha=|\beta|=1}$, with the plane implicitely defined by $\smash{\mathbf{n}\cdot\mathbf{x} = -\mathbf{n}\cdot\mathbf{e}_z}$, with $\smash{\mathbf{n} = \begin{bmatrix} \cos(\theta) & 0 & \sin(\theta)\end{bmatrix}^\intercal}$, as a function of the angle~$\theta$. The semi-major axis of the conic section tends to infinity when $\theta$ approaches a multiple of~$\pi$, meaning that a parametrization of any arc of this conic section with trigonometric functions becomes singular if $\smash{\mathbf{n}\cdot\mathbf{e}_z \to 0}$.}
\label{fig:singular_param}
\end{figure}
This is illustrated in \cref{fig:singular_param} where the semi-major axis of the conic section generated by the intersection of a paraboloid with a plane is plotted as this plane is rotated about the $\mathbf{e}_y$ basis vector. To avoid singular arc parametrizations, we express each conic section arc as a rational B\'ezier curve \cite{Farin2001}. This provides a general parametrization that is valid over all conic section cases (i.e., elliptic, hyperbolic, and parabolic) and allows a seamless and smooth transition between cases. A conic section arc linking a start point $\smash{{\mathbf{x}}_{i,j,0}}$ to an end point $\smash{{\mathbf{x}}_{i,j,1}}$ can be exactly represented by the rational quadratic B\'ezier curve parametrically defined as
\begin{equation}
	\mathbf{x}_{i,j}(t) = \dfrac{B_0(t) \, {\mathbf{x}}_{i,j,0} + w_{i,j} \, B_1(t) \, {\mathbf{x}}_{i,j}^{\star} + B_2(t) \, {\mathbf{x}}_{i,j,1}}{B_0(t) + w_{i,j} \, B_1(t) + B_2(t)} \, , \quad t \in [0,1] \, ,\label{eq:rat_bezier}
\end{equation}
where
\begin{align}
	B_0(t) & = (1-t)^2 \, ,\\
	B_1(t) & = 2(1-t)t \, ,\\
	B_2(t) & = t^2 \, ,
\end{align}
are the Bernstein polynomials of degree $2$, and $\smash{w_{i,j}}$ is a weight associated with the control point $\smash{{\mathbf{x}}_{i,j}^{\star}}$. This control point is located at the intersection of the tangents to the conic section at the start and end points. In the case of the intersection of a planar face $\smash{\mathcal{F}_i}$ with a paraboloid surface $\smash{\mathcal{S}}$, these tangents are obtained by the intersection of the planes tangent to $\smash{\mathcal{S}}$ at the start and end points, with the plane containing $\smash{\mathcal{F}_i}$. The point $\smash{\mathbf{x}_{i,j}\left(\tfrac{1}{2}\right)}$ is, by definition, located at the intersection of $\smash{\mathcal{S}}$ with the segment linking $\smash{\bar{\mathbf{x}}_{i,j}\left(\tfrac{1}{2}\right)}$ to $\smash{{\mathbf{x}}_{i,j}^{\star}}$. Substituting $t$ by $\smash{\tfrac{1}{2}}$ in Eq.~\cref{eq:rat_bezier}, it follows that
\begin{equation}
	w_{i,j} \left(\mathbf{x}_{i,j}\left(\tfrac{1}{2}\right) - {\mathbf{x}}_{i,j}^{\star}\right) =  \bar{\mathbf{x}}_{i,j}\left(\tfrac{1}{2}\right) - \mathbf{x}_{i,j}\left(\tfrac{1}{2}\right) \, ,
\end{equation}
from which $\smash{w_{i,j}}$ can be deducted. If $\smash{\left|w_{i,j}\right|<1}$, the rational B\'ezier curve is an arc of an ellipse, if $\smash{\left|w_{i,j}\right|=1}$, the rational B\'ezier curve is an arc of a parabola, and if $\smash{\left|w_{i,j}\right|>1}$, the rational B\'ezier curve is an arc of a hyperbola (e.g., see \cref{fig:bezierarc}).

Alternatively, the weight $\smash{w_{i,j}}$ (when positive) relates to the absolute curvature $\kappa$ of the conic section arc at its extremities $\smash{{\mathbf{x}}_{i,j,0}}$ and $\smash{{\mathbf{x}}_{i,j,1}}$ following~\cite{Farin1992}
\begin{equation}
    w_{i,j} = 
    \sqrt{\frac{\left|\mathcal{A}\left({\mathbf{x}}_{i,j,0},{\mathbf{x}}_{i,j,1},\mathbf{x}_{i,j}^{\star}\right)\right|}{\kappa({\mathbf{x}}_{i,j,0}) \|\mathbf{x}_{i,j}^{\star}{\mathbf{x}}_{i,j,0}\|^3}}
    = \sqrt{\frac{\left|\mathcal{A}\left({\mathbf{x}}_{i,j,0},{\mathbf{x}}_{i,j,1},\mathbf{x}_{i,j}^{\star}\right)\right|}{\kappa({\mathbf{x}}_{i,j,1}) \|\mathbf{x}_{i,j}^{\star}{\mathbf{x}}_{i,j,1}\|^3}} \, ,
\end{equation}
where $\kappa$ at any point $\mathbf{x}$ on the conic section is given by~\cite{Hartmann1996}
\begin{equation}
    \kappa (\mathbf{x}) = \frac{\left|\left(\mathbf{n}_i \times \nabla \phi(\mathbf{x})\right) \cdot \left(\mathbf{H}_\phi \left(\mathbf{n}_i \times \nabla \phi(\mathbf{x})\right)\right)\right|}{\| \mathbf{n}_i \times \nabla \phi(\mathbf{x})\|} \, ,
\end{equation}
with $\smash{\mathbf{H}_\phi}$ the Hessian matrix
\begin{equation}
    \mathbf{H}_\phi = 2 \begin{bmatrix} \alpha & 0 & 0 \\ 0 & \beta & 0 \\ 0 & 0 & 0 \end{bmatrix} \, .
\end{equation}

Note that, in order to prevent round-off errors in the numerical calculation of $\smash{\boldsymbol{\mathcal{M}}^{\hat{\mathcal{P}}_3}}$, we limit our implementation to cases where $\smash{w_{i,j} \ge 0}$. As a consequence, conic section arcs that would result in negative rational B\'ezier weights are recursively split until positive weights are found. With such a parametrization of the conic section arcs, it can be shown that
\begin{align}
    \boldsymbol{\mathcal{M}}^{\hat{\mathcal{P}}_3}_{i} & = \sum\limits_{j = 1}^{n_{\partial\hat{\mathcal{F}}_i}} 1^{\partial\tilde{\mathcal{S}}}_{i,j} \mathcal{A}\left({\mathbf{x}}_{i,j,0},{\mathbf{x}}_{i,j,1},\mathbf{x}_{i,j}^{\star}\right) \boldsymbol{\mathcal{B}}^{(3)}\left(w_{i,j},{\mathbf{x}}_{i,j,0},{\mathbf{x}}_{i,j,1},\mathbf{x}_{i,j}^{\star}\right) \, , \label{eq:m3k}
\end{align}
where $\smash{\mathcal{A}}$ is the signed projected triangle area operator defined in Eq.~\cref{eq:triangle_area} and the operator $\smash{\boldsymbol{\mathcal{B}}^{(3)}: \mathbb{R}^{+} \times \mathbb{R}^3 \times \mathbb{R}^3 \times \mathbb{R}^3 \to \mathbb{R}^4}$ is given in \cref{apdx:M3}.

\subsection{A special elliptic case}
When $\smash{\alpha\beta > 0}$ (i.e., $\smash{\mathcal{S}}$ is an elliptic paraboloid) and the normal to $\smash{\mathcal{F}_i}$ is such that $\smash{\mathbf{n}_i \cdot \mathbf{e}_z \ne 0}$, the intersection of the plane containing the face $\smash{\mathcal{F}_i}$ with the surface $\smash{\mathcal{S}}$ is an ellipse. The intersection of $\smash{\mathcal{F}_i}$ with $\smash{\mathcal{S}}$ can then be: empty, a collection of arcs of this ellipse, or the entire ellipse. In the latter case, the sum of the contributions $\smash{\boldsymbol{\mathcal{M}}^{\hat{\mathcal{P}}_2}_{i}}$ and $\smash{\boldsymbol{\mathcal{M}}^{\hat{\mathcal{P}}_3}_{i}}$ can be directly calculated by integrating $\smash{\boldsymbol{\Psi}_\mathcal{S}}$ and $\smash{\boldsymbol{\Psi}_{\mathcal{F}_i}}$ over the full ellipse, which yields the more concise expression 
\begin{equation}
    \boldsymbol{\mathcal{M}}^{\hat{\mathcal{P}}_2}_{i}+\boldsymbol{\mathcal{M}}^{\hat{\mathcal{P}}_3}_{i} =  -\mathrm{sign}\left(\mathbf{n}_i \cdot \mathbf{e}_z\right)\pi\frac{(\tau_i^2 \alpha + \lambda_i^2 \beta - 4 \alpha \beta \delta_i)^2}{32\left(\alpha\beta\right)^{5/2}} 
    \begin{bmatrix}[1.2]
        1 \\
        \frac{\lambda_i }{2\alpha} \\
        \frac{\tau_i }{2\beta} \\
        \frac{5\tau_i^2}{12\beta} + \frac{5\lambda_i^2}{12\alpha} - \frac{2\delta_i}{3}
    \end{bmatrix} \, ,
\end{equation}
where $\smash{\delta_i}$, $\smash{\lambda_i}$, and $\smash{\tau_i}$ have been defined in Eqs.~\cref{eq:distance,eq:lambda,eq:tau}. Note that this case only occurs for $\smash{\mathbf{n}_i \cdot \mathbf{e}_z \ne 0}$, hence these coefficients are here non-singular.

\section{Integrating on $\tilde{\mathcal{S}}$}
\label{sec:surface}
Although this manuscript is concerned with estimating the first moments of $\smash{\hat{\mathcal{P}}}$, the integration domains and parametrizations introduced in \Cref{sec:moments} can also be used for integrating quantities associated with the clipped surface $\smash{\tilde{\mathcal{S}}}$, e.g., its moments. The area of $\tilde{\mathcal{S}}$, for instance, given as
\begin{align}
	\mathrm{M}_0^{\tilde{\mathcal{S}}} & = \int_{\tilde{\mathcal{S}}} 1 \ \mathrm{d}\mathbf{x} \, ,
\end{align}
also reads after application of the divergence theorem as
\begin{align}
	\mathrm{M}_0^{\tilde{\mathcal{S}}} & =\sum\limits_{i=1}^{n_{\!\mathcal{F_{\phantom{k}\!\!}}}} \sum\limits_{j = 1}^{n_{\partial\hat{\mathcal{F}}_i}} -1^{\partial\tilde{\mathcal{S}}}_{i,j} \int_0^1 \gamma(x_{i,j}(t),y_{i,j}(t)) y^\prime_{i,j}(t) \ \mathrm{d}t \, , \label{eq:surface_area}
\end{align}
where
\begin{equation}
    \gamma : \left\{\begin{array}{lcl} 
        \mathbb{R}^2 & \!\to & \mathbb{R} \\
        (x,y) & \!\mapsto & {\displaystyle \int_0^x} \sqrt{1+4\alpha^2x^2 + 4\beta^2y^2} \ \mathrm{d}x \end{array} \right.  \, .
\end{equation}
The first moments of $\smash{\tilde{\mathcal{S}}}$, as well as the average normal vector, average Gaussian curvature and average mean curvature of $\smash{\tilde{\mathcal{S}}}$, for example, can be expressed similarly as sums of one-dimensional integrals over the parametric arcs of $\smash{\partial\tilde{\mathcal{S}}^\perp}$. Contrary to the first moments of $\smash{\hat{\mathcal{P}}}$, however, they need to be estimated using a numerical quadrature rule, as closed-form expressions cannot be derived.

\section{On floating-point arithmetics and robustness}
\label{sec:robustness}
In the context of floating-point arithmetics, there exist cases for which the computed topology of the clipped faces $\smash{\hat{\mathcal{F}}_i}$ may be ill-posed, preventing the accurate calculation of the moments of $\smash{\hat{\mathcal{P}}}$. This occurs when, in the discrete sense:
\begin{enumerate}
    \item The surface $\smash{\mathcal{S}}$ is tangent to one or more edges of the polyhedron $\smash{\mathcal{P}}$;
    \item At least one corner or vertex of the polyhedron $\smash{\mathcal{P}}$ belongs to the surface $\smash{\mathcal{S}}$.
\end{enumerate}
For any given face of $\smash{\mathcal{P}}$, the former case is numerically detected by computing the absolute value of the dot product between the normalized tangent to $\smash{\mathcal{S}}$ and the normalized edge from which the tangent originates, and checking whether it lies within $\epsilon_\text{tangent}$ of unity. The latter case is detected by checking whether the intersection of an edge with $\smash{\mathcal{S}}$ lies within $\epsilon_\text{corner}$ of a corner or vertex of the polyhedron $\mathcal{P}$. When any of these configurations is detected, the polyhedron $\smash{\mathcal{P}}$ is randomly translated and rotated about its centroid by a distance and an angle equal to $\smash{\epsilon_\text{nudge}}$, and the clipped face discrete topologies are re-computed. Moreover, we also switch to a higher-accuracy floating-point format, e.g., from a $64$-bit to a $128$-bit format.
In the current work, where we aim to produce results in ``double precision'', we find the values $\smash{\epsilon_\text{nudge} = 10^{10} \times \epsilon_{\mathrm{128}}}$, $\smash{\epsilon_\text{corner} = 10^2 \times \epsilon_{\mathrm{64}}}$, and $\smash{\epsilon_\text{tangent} =  10^6 \times \epsilon_{\mathrm{64}}}$, with $\smash{\epsilon_{\mathrm{64}} = 2^{-52}}$ and $\smash{\epsilon_{\mathrm{128}} = 2^{-112}}$ the upper bounds of the relative approximation error in $64$-bit and $128$-bit floating-point arithmetics, respectively, to prevent the computation of any ill-posed topologies in all tests presented in \Cref{sec:tests} (for which more than $\smash{5\times10^7}$ occurences of the current ``nudging'' procedure are forced to occur). Choosing lower values for these tolerances may result in the generation of non-valid discrete topologies and/or erroneous moments. It should be noted that, for a polyhedron $\smash{\mathcal{P}}$ with volume $\smash{\mathrm{M}^\mathcal{P}_0 = \mathcal{O}(1)}$, the rate of occurence of the cases triggering this correction procedure is extremely small (that is, for non-engineered intersection configurations).

\section{Higher-order moments}
\label{sec:highorder}
Second- and higher-order moments of the clipped polyhedron $\smash{\hat{\mathcal{P}}}$ can be derived using the same procedure as presented in \Cref{sec:moments} for the calculation of the zeroth and first moments. This merely entails appending higher-order monomials to the function vector $\smash{\boldsymbol{\Upsilon}}$ defined in Eq.~\cref{eq:defmonomials}. The arc parametrizations introduced in \Cref{sec:moments} and used for calculating the three moment contributions $\smash{\boldsymbol{\mathcal{M}}^{\hat{\mathcal{P}}_1}}$, $\smash{\boldsymbol{\mathcal{M}}^{\hat{\mathcal{P}}_2}}$, and $\smash{\boldsymbol{\mathcal{M}}^{\hat{\mathcal{P}}_3}}$ would remain unchanged, while the operators $\smash{\boldsymbol{\mathcal{B}}^{(1)}}$, $\smash{\boldsymbol{\mathcal{B}}^{(2)}}$, and $\smash{\boldsymbol{\mathcal{B}}^{(3)}}$ would then contain additional components respectively corresponding to the mononials appended to $\smash{\boldsymbol{\Upsilon}}$. We hypothesize that closed-form expressions similar to those given in \Cref{sec:moments} can be derived for moments of $\smash{\hat{\mathcal{P}}}$ of arbitrary order. However, these would involve an ever increasing amount of high-order monomials of $\alpha$, $\beta$, and of the components of the vertices of $\smash{\hat{\mathcal{P}}}$ and of the control points of the conic section arcs in $\smash{\partial\partial\hat{\mathcal{P}}}$, rendering the computation of the moments more sensitive to round-off errors.

\section{Verification}
\label{sec:tests}
In this section, the closed-form expressions derived in \Cref{sec:moments}, the approach of \Cref{sec:surface} for integrating on the clipped surface, and the corrective procedure of \Cref{sec:robustness} are tested on a wide variety of engineered and random intersection configurations. When analytical expressions of the exact moments are not available, we recursively split the faces of the polyhedron $\mathcal{P}$ so as to approximate $\smash{\partial{\mathcal{P}}\cap\mathcal{Q}}$ and $\smash{\partial\mathcal{P}\cap(\mathbb{R}^3\setminus\mathcal{Q})}$ by collections of oriented triangles. We refer to this procedure as the adaptive mesh refinement (AMR) of the faces of the polyhedron $\mathcal{P}$. We then exactly integrate $\smash{\boldsymbol{\Phi}_{\mathcal{S}}}$ on the triangulated approximation of $\smash{\partial\mathcal{P}\cap(\mathbb{R}^3\setminus\mathcal{Q})}$ and $\smash{\boldsymbol{\Phi}_{\mathcal{F}_i}}$ on the triangulated approximation of each clipped face $\smash{\hat{\mathcal{F}}_i}$, effectively approximating Eq.~\cref{eq:vol_after_div4}. For each case, we ensure that enough levels of recursive refinement are employed in order to reach machine-zero. Accumulated errors due to the summation of the contributions of all triangles are avoided by the use of compensated summation, also known as Kahan summation~\cite{Kahan1965}.

\subsection{Unit cube translating along $\mathbf{e}_z$}

\begin{figure}[h]
    \centering
    \begin{tabular}{ccc}
    \subfloat[$k=1/2$]{\adjincludegraphics[width=0.3\textwidth,trim=0 0 0 0,clip=true]{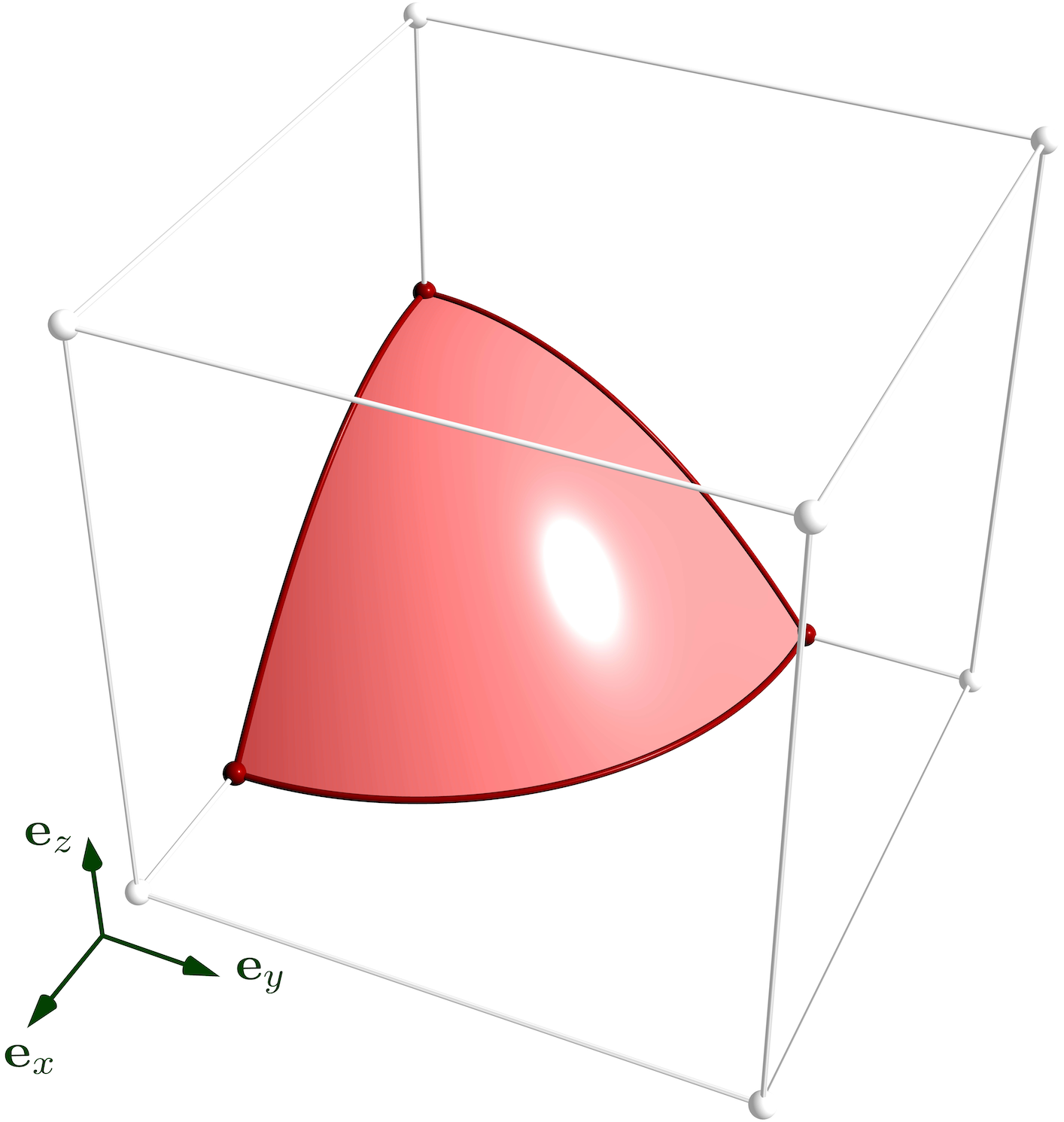}} \quad
    \subfloat[$k=3/2$]{\adjincludegraphics[width=0.3\textwidth,trim=0 0 0 0,clip=true]{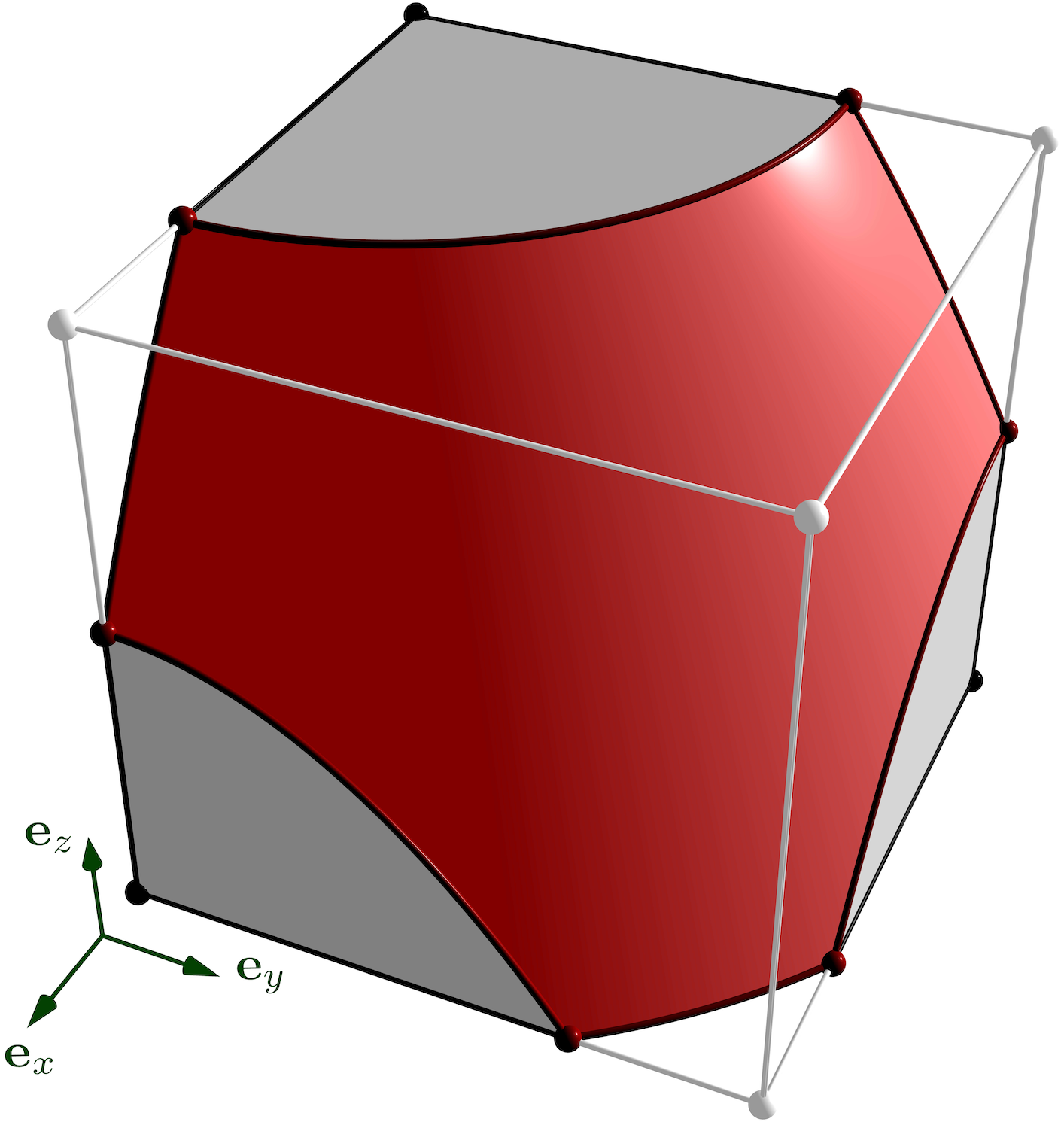}} \quad
    \subfloat[$k=5/2$]{\adjincludegraphics[width=0.3\textwidth,trim=0 0 0 0,clip=true]{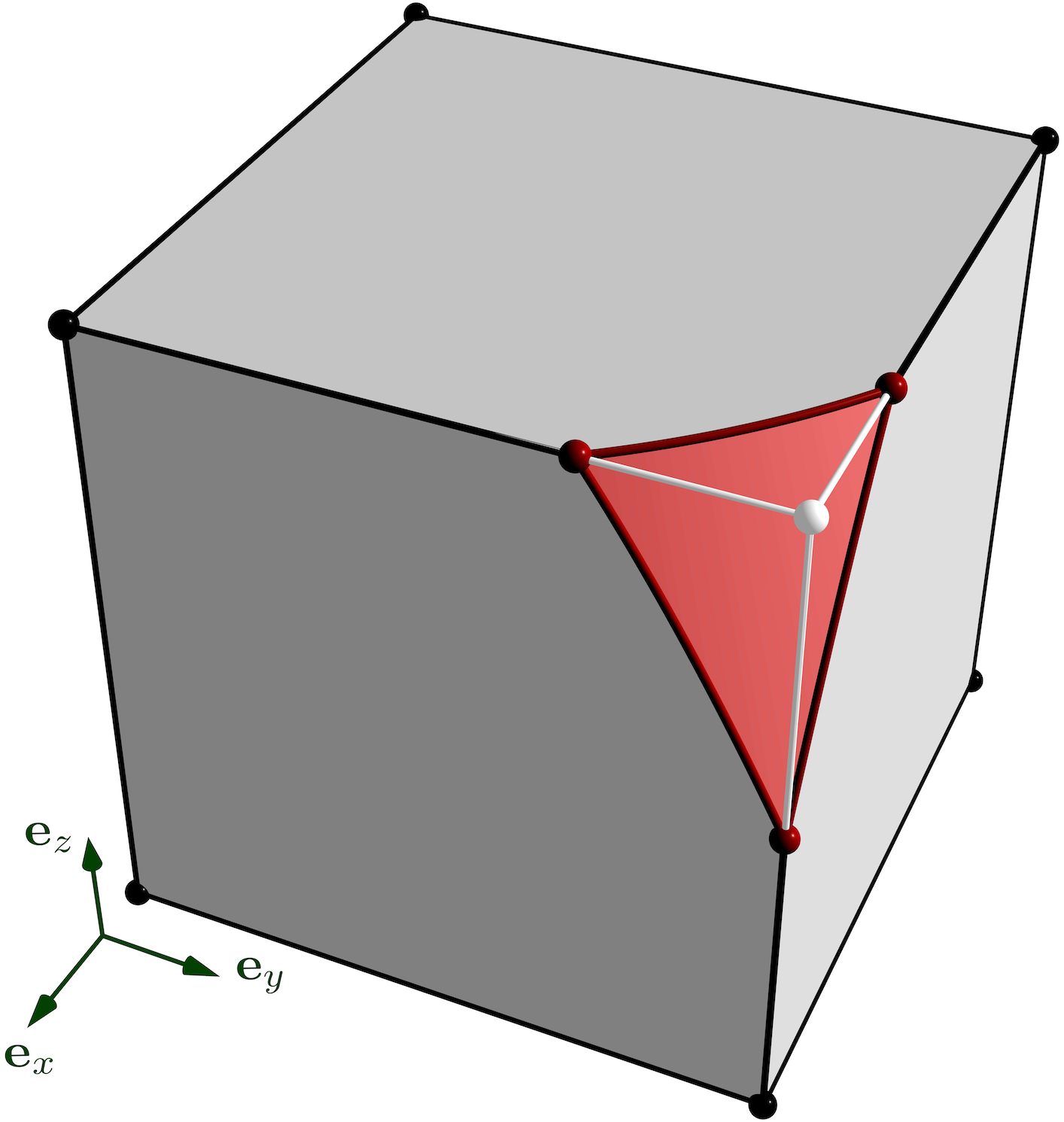}}
    \end{tabular}
    \caption{Unit cube centred at $\begin{bmatrix} 1/2 & 1/2 & 1/2 - k \end{bmatrix}^\intercal$ clipped by the elliptic paraboloid parametrically defined as $z = -x^2 - y^2$.}
    \label{fig:testcube}
\end{figure}

In a first test, we consider the elliptic paraboloid defined by Eq.~\cref{def:paraboloid} with $\alpha = \beta = 1$, intersecting with the unit cube centered at $\smash{\mathbf{x}_c = \begin{bmatrix} 1/2 &1/2 & 1/2 - k \end{bmatrix}^\intercal}$ (as illustrated in \cref{fig:testcube}). For this case, the zeroth and first moments of $\smash{\hat{\mathcal{P}}}$, as well as the zeroth moment of $\smash{\tilde{\mathcal{S}}}$, can be derived as analytical functions of $k$\footnote{A \texttt{C++} implementation of these functions can be found at the beginning of \href{https://github.com/robert-chiodi/interface-reconstruction-library/blob/paraboloid_cutting/tests/src/paraboloid_intersection_test.cpp}{\texttt{/tests/src/paraboloid\_intersection\_test.cpp}} in the open-source \href{https://github.com/robert-chiodi/interface-reconstruction-library/tree/paraboloid_cutting}{\texttt{Interface Reconstruction Library}}.}. We compare these exact moments against those computed using the closed-form expressions derived in \Cref{sec:moments} for estimating the moments of $\smash{\hat{\mathcal{P}}}$, and by integrating Eq.~\cref{eq:surface_area} numerically with an adaptive Gauss-Legendre quadrature rule for estimating the zeroth moment of $\smash{\tilde{\mathcal{S}}}$. The parameter $k$ is regularly sampled on $[0,3]$ with a uniform spacing $\smash{\Delta k = 10^{-3}}$.

\begin{figure}
    \centering
    \includegraphics{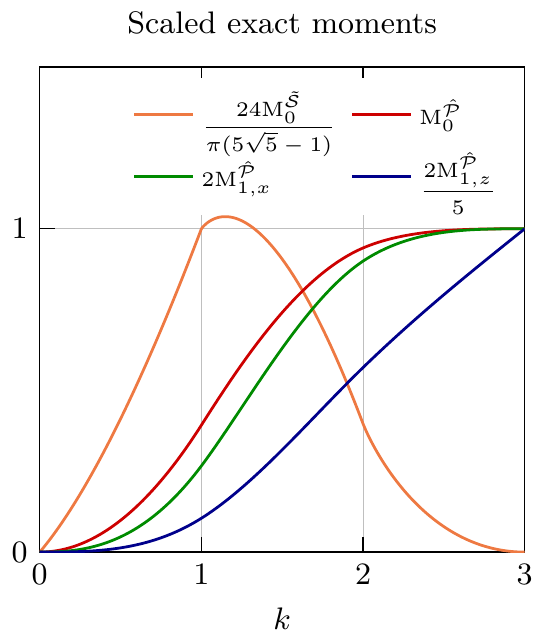}
    \includegraphics{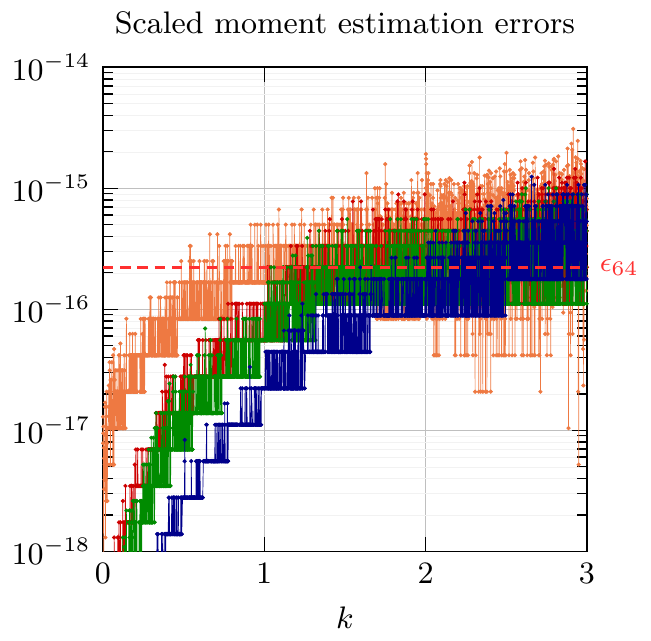}
    \caption{Moments of the unit cube centred at $\begin{bmatrix} 1/2 & 1/2 & 1/2 - k \end{bmatrix}^\intercal$ clipped by the elliptic paraboloid parametrically defined as $z = -x^2 - y^2$, and the error in their estimation with $64$-bit floating point arithmetics. The volume moments and their estimation errors are scaled with respect to the moments at $k=3$. The surface area and its estimation error is scaled with respect to the surface area at $k=1$. The volume moments are computed from the analytical expressions derived in~\Cref{sec:moments}, whereas the surface area is computed from Eq.~\cref{eq:surface_area} using an adaptive Gauss-Legendre quadrature rule. The $64$-bit machine-epsilon $\epsilon_{64} = 2^{-52}$ is shown as the dashed red line.}
    \label{fig:resultscube}
\end{figure}

The left-hand graph of \cref{fig:resultscube} shows the exact moments of $\smash{\hat{\mathcal{P}}}$, scaled by their maximum values, as well as the exact zeroth moment of $\smash{\tilde{\mathcal{S}}}$, scaled by its value at $k=1$. The right-hand side of \cref{fig:resultscube} shows the errors associated with their estimation, scaled similarly. These errors are all contained within an order of magnitude of $\smash{\epsilon_{64} = 2^{-52}}$, the upper bound of the relative approximation error in $64$-bit floating-point arithmetics.

\subsection{Parameter sweep for several geometries} \label{sec:paramsweep}

\newcolumntype{R}{>{$}r<{$}} %
\newcolumntype{V}[1]{>{[\;}*{#1}{R@{\;\;}}R<{\;]}} %

\begin{table}[tbhp]
    \renewcommand\arraystretch{1.4}
    \scriptsize
    \captionsetup{position=top} 
    \caption{Five polyhedra are considered: a regular tetrahedron, a cube, a regular dodecahedron, a hollow cube, and the triangulated Stanford bunny \cite{StanfordBunny}. The first three polyhedra are convex, whereas the last two are not. The hollow cube contains non-convex faces.}\label{tab:geometries}
    \begin{center}
    \begin{tabular}{|r|r|r|r|r|r|} \hline
                \multicolumn{1}{|c|}{\textbf{Geometry}} & Tetrahedron & Cube & Dodecahedron & Hollow cube & Stanford bunny \\ \hline
                \multicolumn{1}{|c|}{\textbf{Number of}} & \multirow{2}{*}{$4$} & \multirow{2}{*}{$8$} & \multirow{2}{*}{$20$} & \multirow{2}{*}{$16$} &  \multirow{2}{*}{$167,891$}\\[-4pt] 
                \multicolumn{1}{|c|}{\textbf{vertices}} & & & & & \\ \hline
                \multicolumn{1}{|c|}{\textbf{Number of}} & \multirow{2}{*}{$4$} & \multirow{2}{*}{$6$} & \multirow{2}{*}{$12$} & \multirow{2}{*}{$12$} &  \multirow{2}{*}{$335,778$}\\[-4pt]  
                \multicolumn{1}{|c|}{\textbf{faces}} & & & & & \\ \hline
                \multicolumn{1}{|c|}{\textbf{Are all faces}} & \multirow{2}{*}{Yes} & \multirow{2}{*}{Yes} & \multirow{2}{*}{Yes} & \multirow{2}{*}{No} & \multirow{2}{*}{Yes} \\[-4pt] 
                \multicolumn{1}{|c|}{\textbf{convex?}} & & & & & \\ \hline
                \multicolumn{1}{|c|}{\textbf{Is polyhedron}} & \multirow{2}{*}{Yes} & \multirow{2}{*}{Yes} & \multirow{2}{*}{Yes} & \multirow{2}{*}{No} & \multirow{2}{*}{No} \\[-4pt]                  
                \multicolumn{1}{|c|}{\textbf{convex?}} & & & & & \\ \hline
                \multicolumn{1}{|c|}{\textbf{Snapshot}} 
                &  \multicolumn{1}{c|}{\begin{minipage}{0.12\textwidth}\centering\adjincludegraphics[width=0.99\textwidth,trim=0 0 0 0,clip=true]{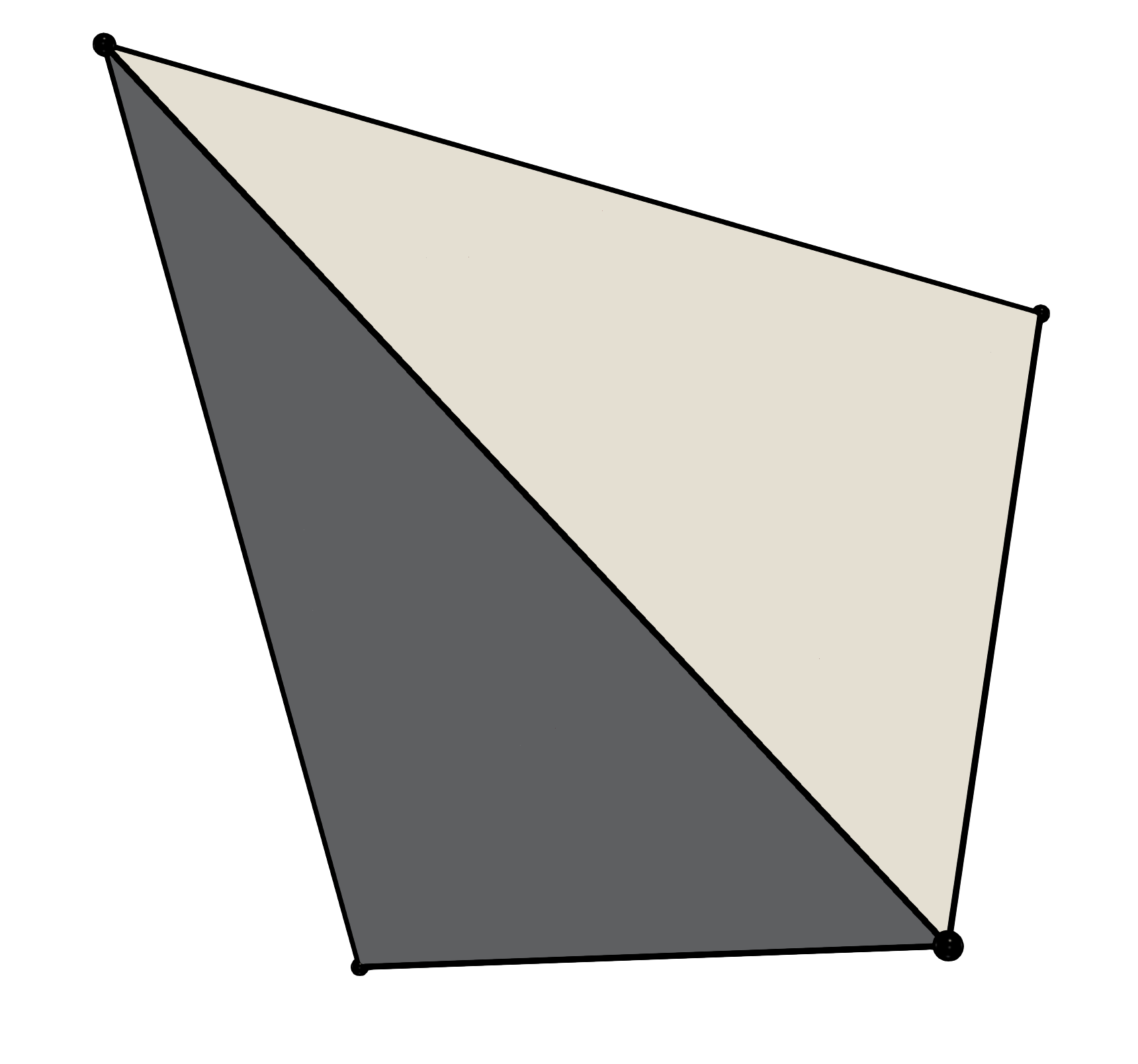}\end{minipage}} 
                &  \multicolumn{1}{c|}{\begin{minipage}{0.12\textwidth}\centering\adjincludegraphics[width=0.99\textwidth,trim=0 0 0 0,clip=true]{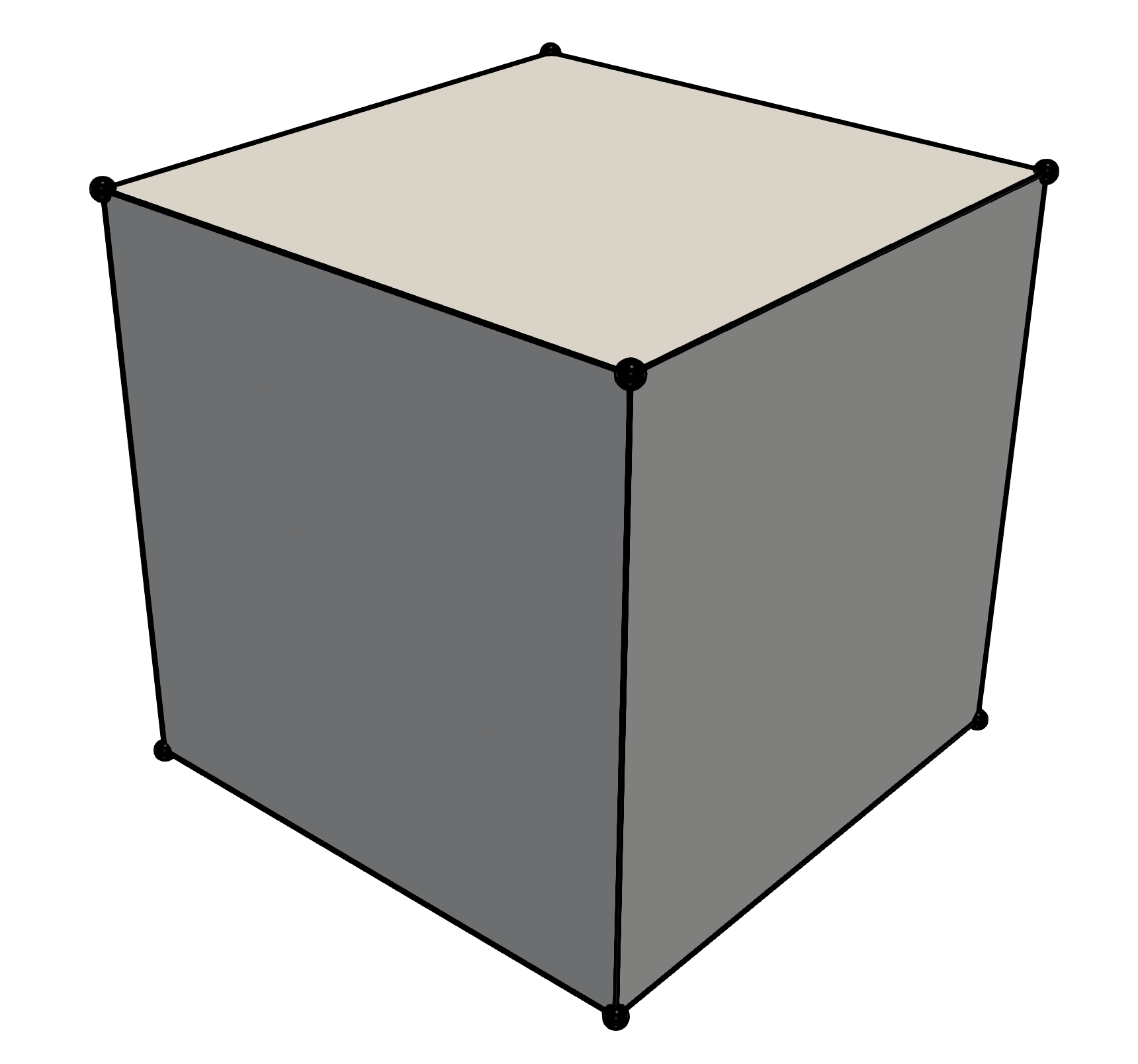}\end{minipage}} 
                &  \multicolumn{1}{c|}{\begin{minipage}{0.12\textwidth}\centering\adjincludegraphics[width=0.99\textwidth,trim=0 0 0 0,clip=true]{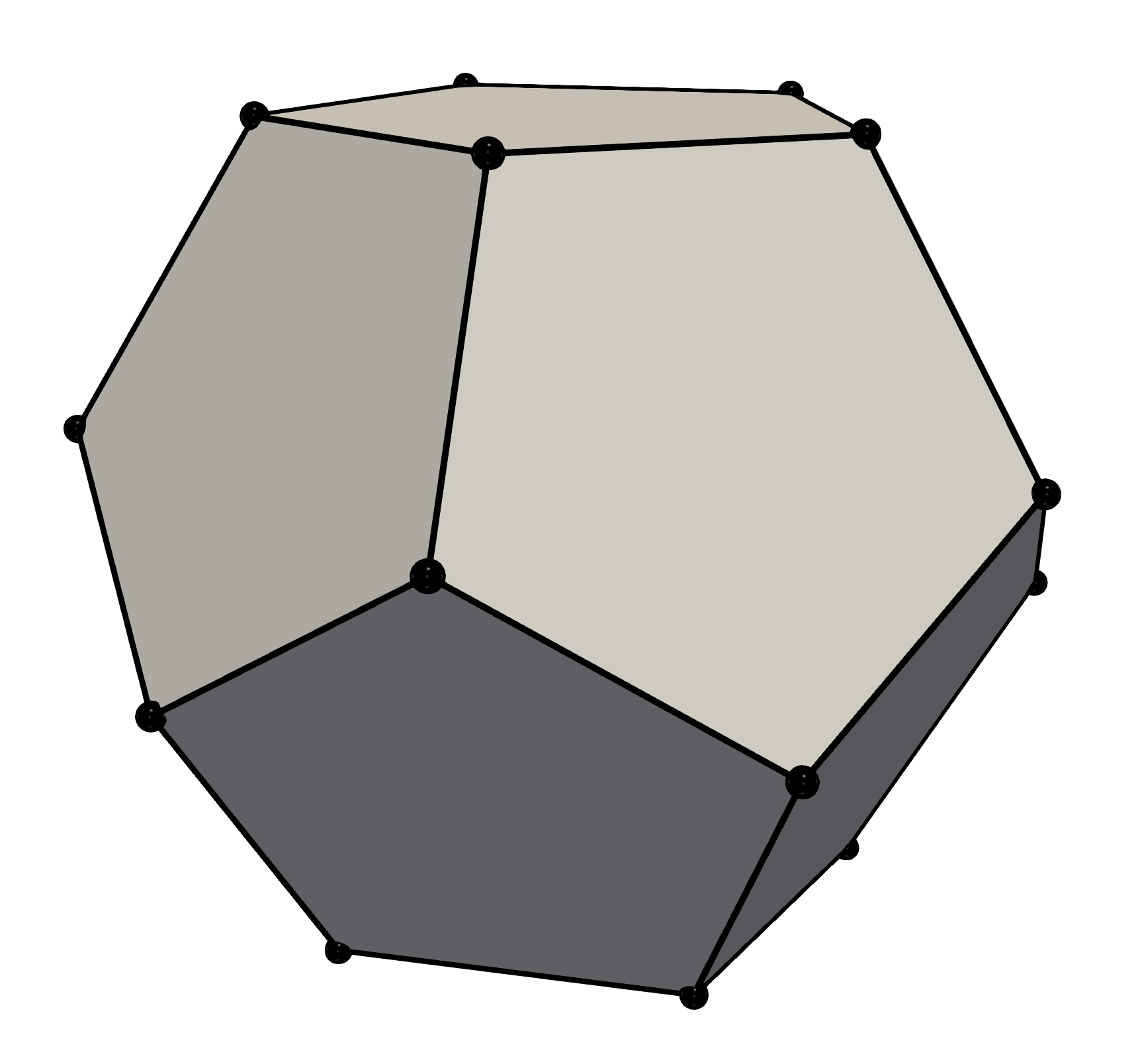}\end{minipage}} 
                &  \multicolumn{1}{c|}{\begin{minipage}{0.12\textwidth}\centering\adjincludegraphics[width=0.99\textwidth,trim=0 0 0 0,clip=true]{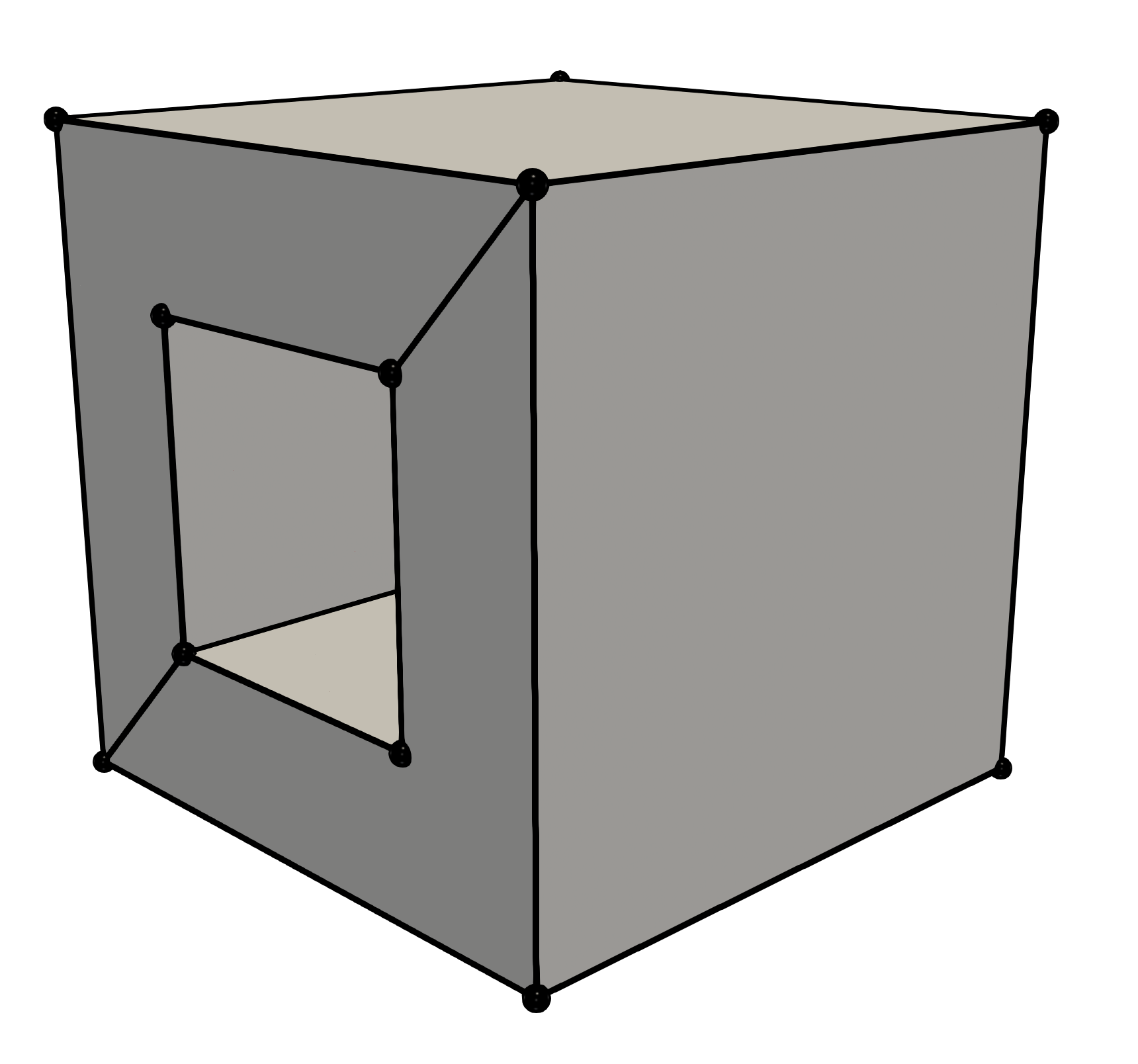}\end{minipage}} 
                &  \multicolumn{1}{c|}{\begin{minipage}{0.12\textwidth}\centering\vspace{0.5mm}\adjincludegraphics[width=0.9\textwidth,trim=0 0 0 0,clip=true]{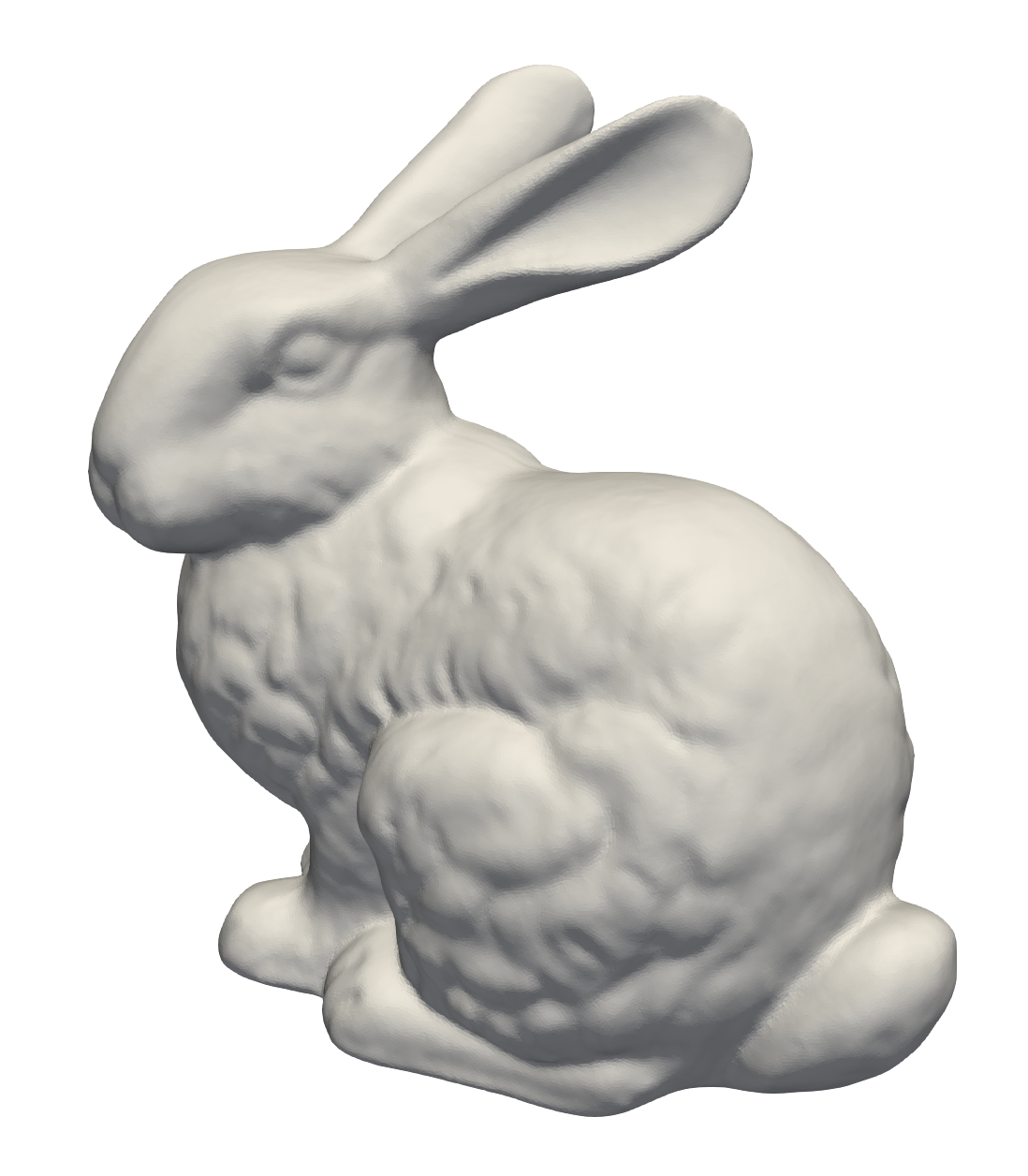}\end{minipage}} 
                \\ \hline
    \end{tabular}
    \end{center}
\end{table}

In a second test, we consider a selection of convex polyhedra (a regular tetrahedron, a cube, and a regular dodecahedron) and non-convex polyhedra (a hollow cube and the triangulated Stanford bunny \cite{StanfordBunny}). These polyhedra, whose properties are summarized in \cref{tab:geometries}, are scaled so as to have a unit volume (i.e., $\smash{\mathrm{M}_0^\mathcal{P} = 1}$), and a centroid initially located at the origin (i.e., $\smash{\mathbf{M}_1^\mathcal{P} = \mathbf{0}}$). They are then translated along $\mathbf{t} = \begin{bmatrix} t_x & t_y & t_z \end{bmatrix}^\intercal$, and rotated about the three basis vectors $\smash{(\mathbf{e}_x, \mathbf{e}_y,\mathbf{e}_z)}$ with the angles $\theta_x$, $\theta_y$, and $\theta_z$, respectively. Throughout this section, $(t_x,t_y,t_z)$ are varied in $[-\tfrac{1}{2},\tfrac{1}{2}]^3$, $(\theta_x,\theta_y,\theta_z)$ are varied in $[-\pi,\pi]^3$, and the paraboloid coefficients $(\alpha,\beta)$ are varied in $[-5,5]^2$.\medskip

\begin{figure}[tbhp] \centering
    \begin{tabular}{c}
        \subfloat[Tetrahedron]{\begin{minipage}{0.30\textwidth}\centering\adjincludegraphics[width=3cm,trim=0 0 0 0,clip=true]{figure7a.png}\label{fig:tet}\end{minipage}}
        \subfloat[Clipped tetrahedron]{\begin{minipage}{0.30\textwidth}\centering\adjincludegraphics[width=3cm,trim=0 0 0 0,clip=true]{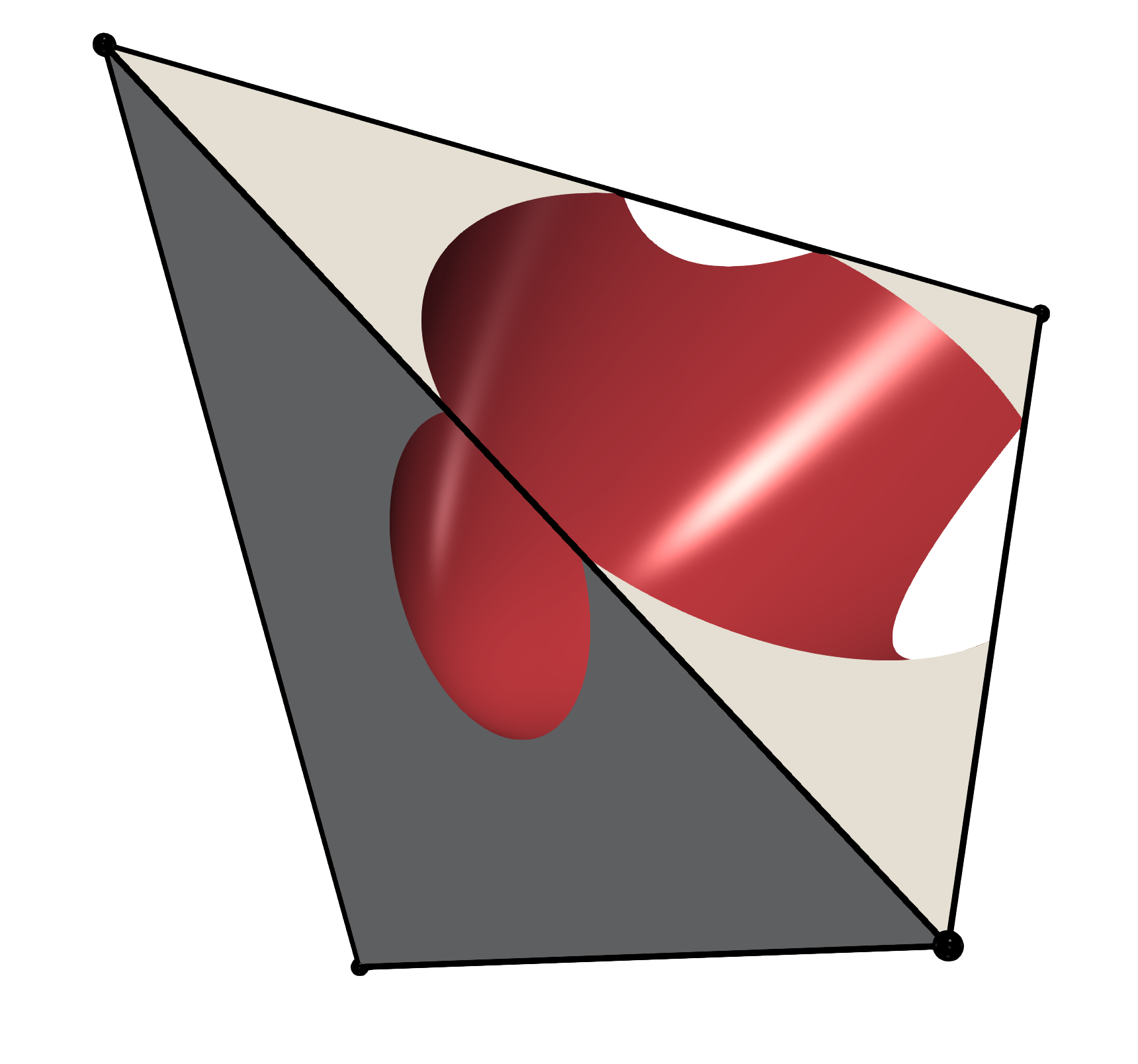}\end{minipage}}
        \subfloat[Clipped tetrahedron (AMR ref.)]{\begin{minipage}{0.39\textwidth}\centering\adjincludegraphics[width=3cm,trim=0 0 0 0,clip=true]{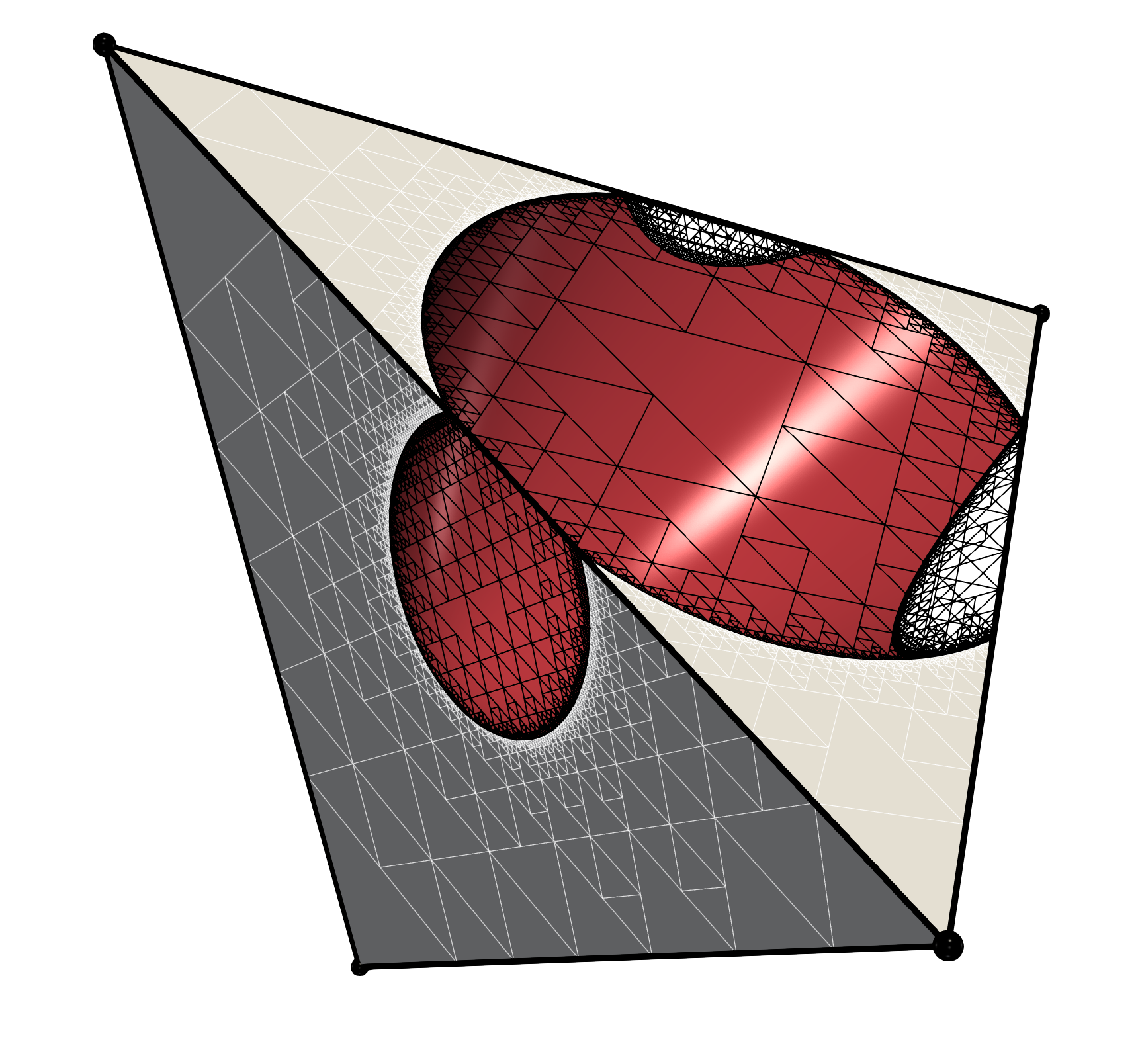}\end{minipage}}\\
        \subfloat[Cube]{\begin{minipage}{0.30\textwidth}\centering\adjincludegraphics[width=3cm,trim=0 0 0 0,clip=true]{figure7d.png}\label{fig:cube}\end{minipage}}
        \subfloat[Clipped cube]{\begin{minipage}{0.30\textwidth}\centering\adjincludegraphics[width=3cm,trim=0 0 0 0,clip=true]{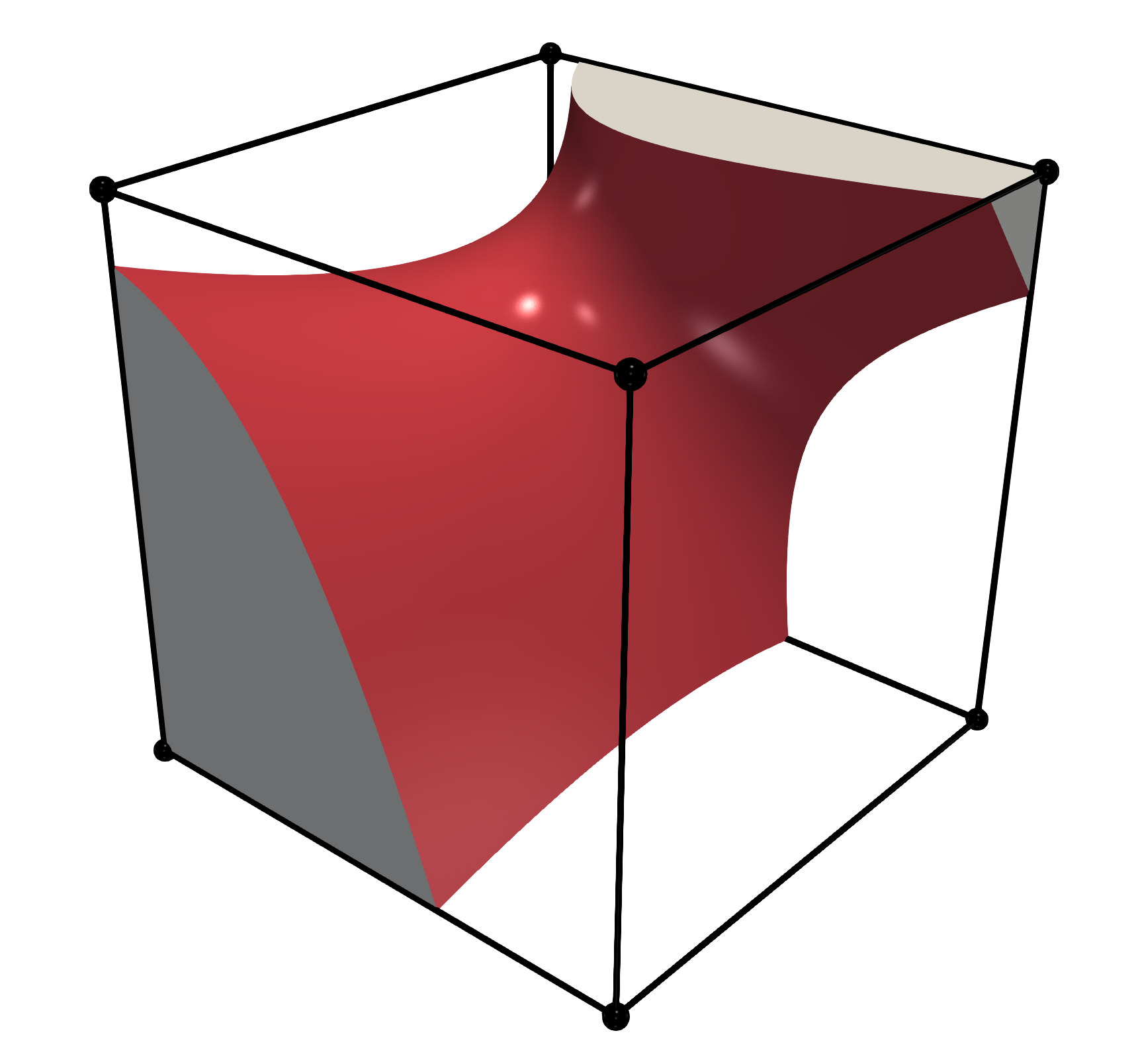}\end{minipage}}
        \subfloat[Clipped cube (AMR ref.)]{\begin{minipage}{0.39\textwidth}\centering\adjincludegraphics[width=3cm,trim=0 0 0 0,clip=true]{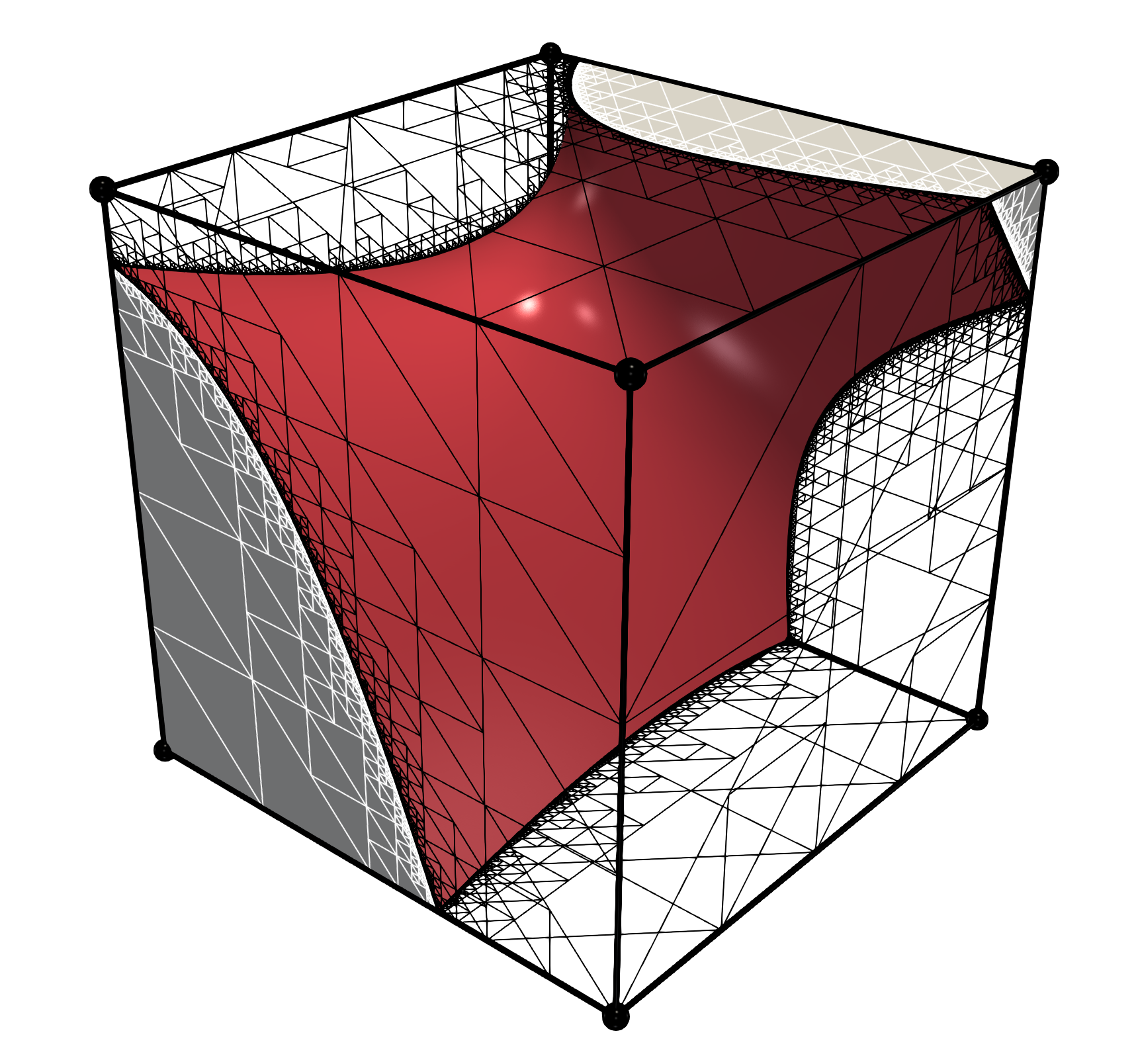}\end{minipage}}\\
        \subfloat[Dodecahedron]{\begin{minipage}{0.30\textwidth}\centering\adjincludegraphics[width=3cm,trim=0 0 0 0,clip=true]{figure7g.png}\label{fig:dodeca}\end{minipage}}
        \subfloat[Clipped dodecahedron]{\begin{minipage}{0.30\textwidth}\centering\adjincludegraphics[width=3cm,trim=0 0 0 0,clip=true]{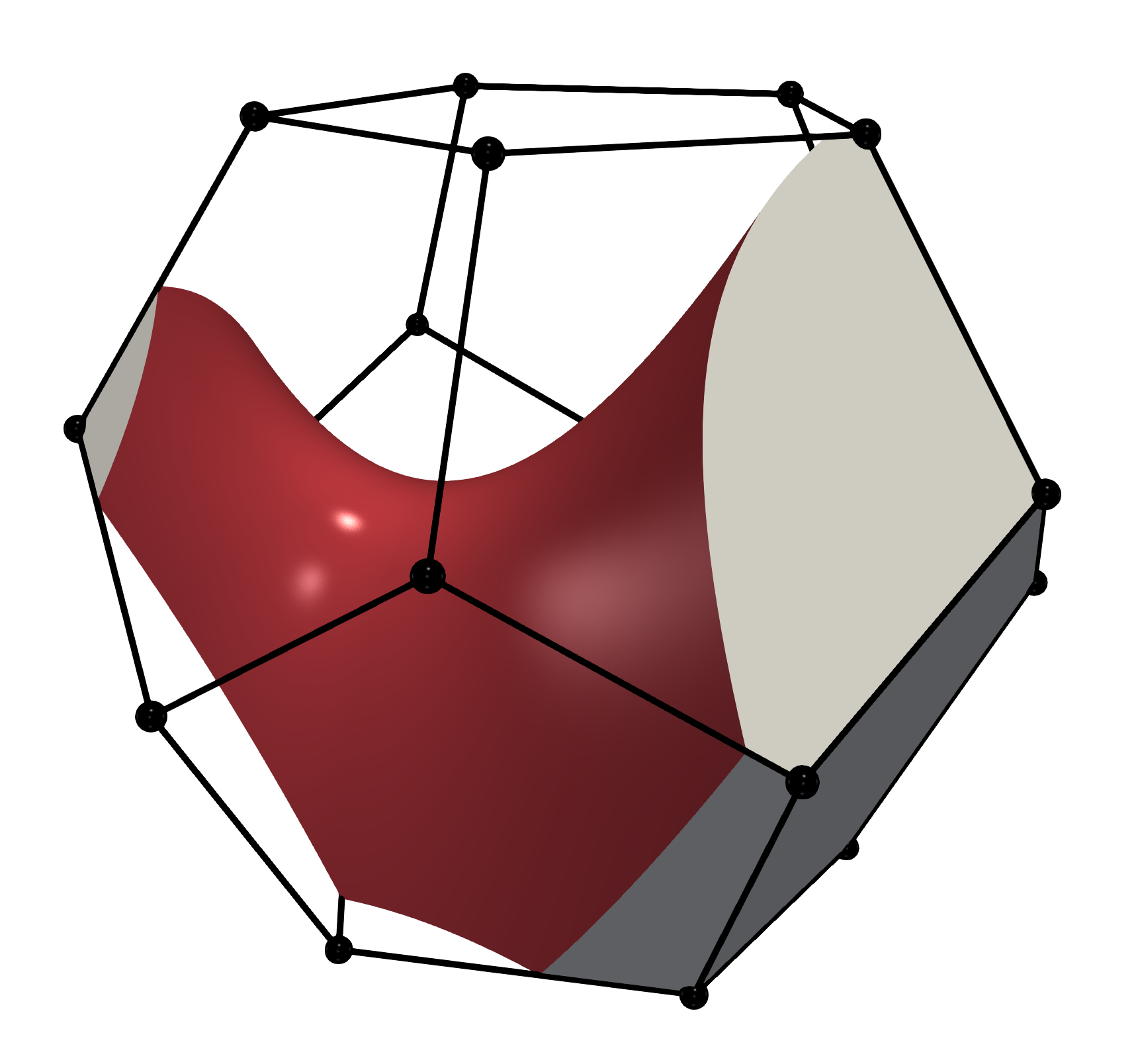}\end{minipage}}
        \subfloat[Clipped dodecahedron (AMR ref.)]{\begin{minipage}{0.39\textwidth}\centering\adjincludegraphics[width=3cm,trim=0 0 0 0,clip=true]{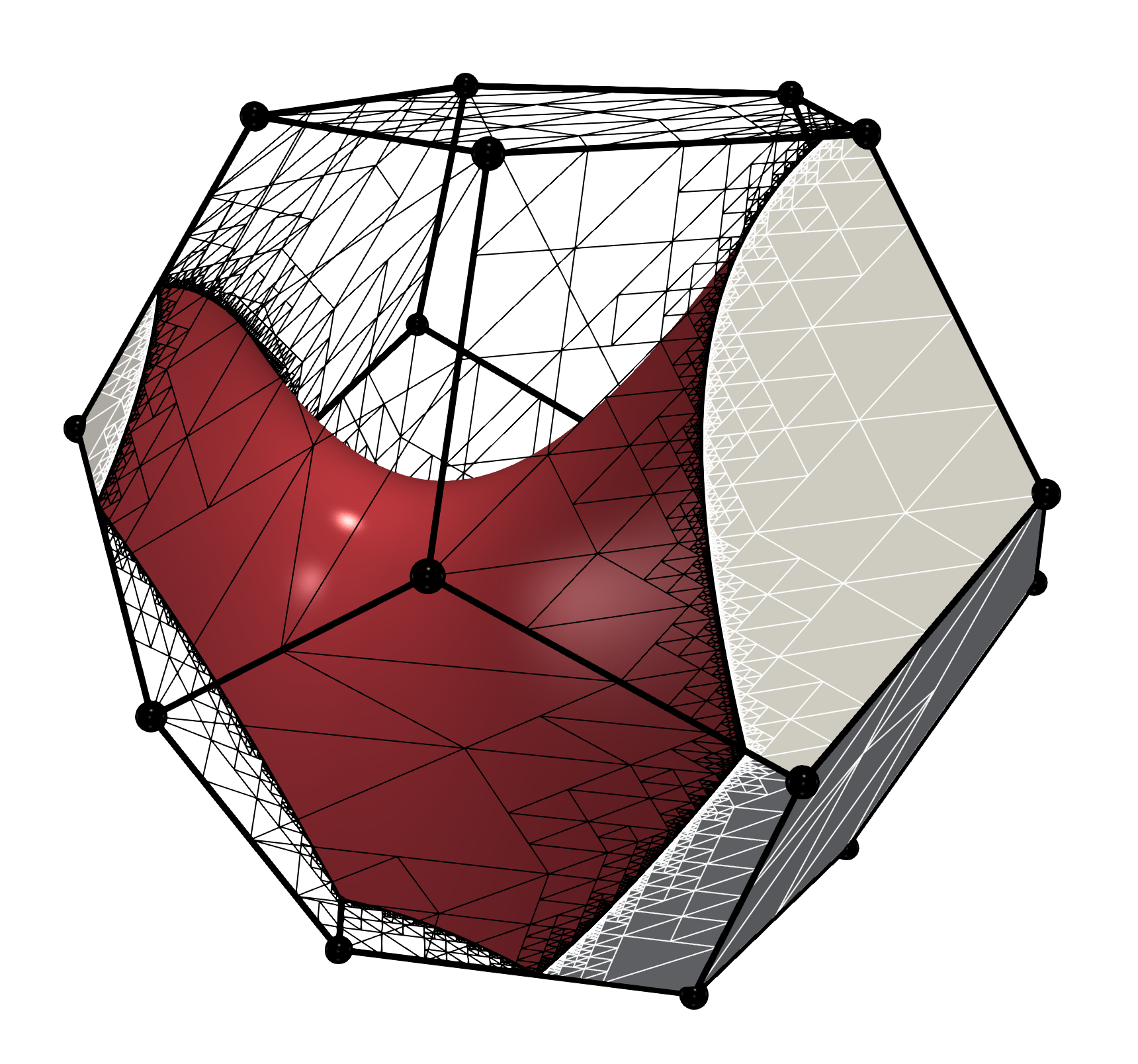}\end{minipage}}\\
        \subfloat[Hollow cube]{\begin{minipage}{0.30\textwidth}\centering\adjincludegraphics[width=3cm,trim=0 0 0 0,clip=true]{figure7j.png}\label{fig:cubehole}\end{minipage}}
        \subfloat[Clipped hollow cube]{\begin{minipage}{0.30\textwidth}\centering\adjincludegraphics[width=3cm,trim=0 0 0 0,clip=true]{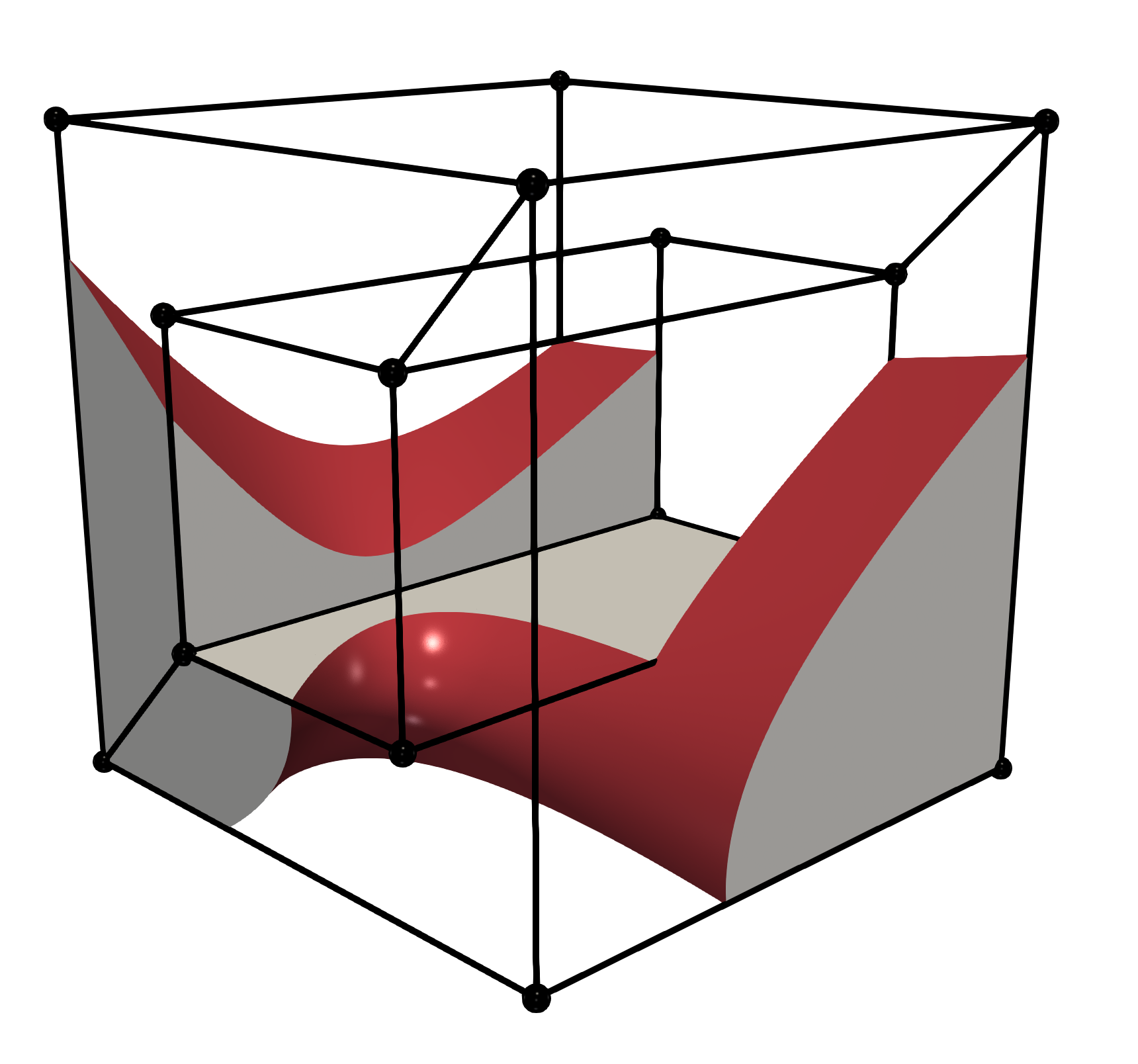}\end{minipage}}
        \subfloat[Clipped hollow cube (AMR ref.)]{\begin{minipage}{0.39\textwidth}\centering\adjincludegraphics[width=3cm,trim=0 0 0 0,clip=true]{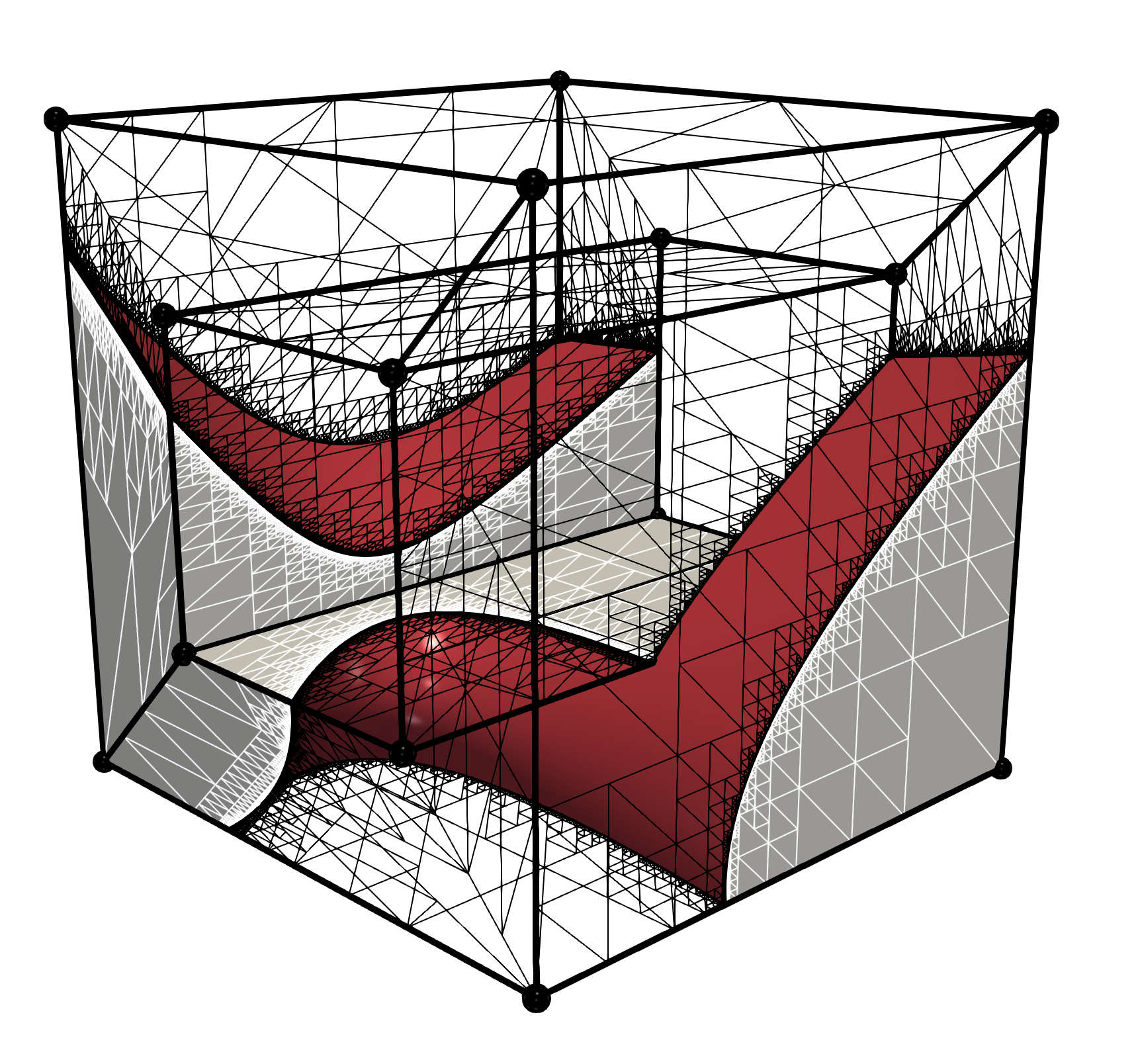}\end{minipage}}\\
        \subfloat[Stanford bunny \cite{StanfordBunny}]{\begin{minipage}{0.30\textwidth}\centering\adjincludegraphics[width=3cm,trim=0 0 0 0,clip=true]{figure7m.png}\label{fig:bunny}\end{minipage}}
        \subfloat[Clipped bunny]{\begin{minipage}{0.30\textwidth}\centering\adjincludegraphics[width=3cm,trim=0 0 0 0,clip=true]{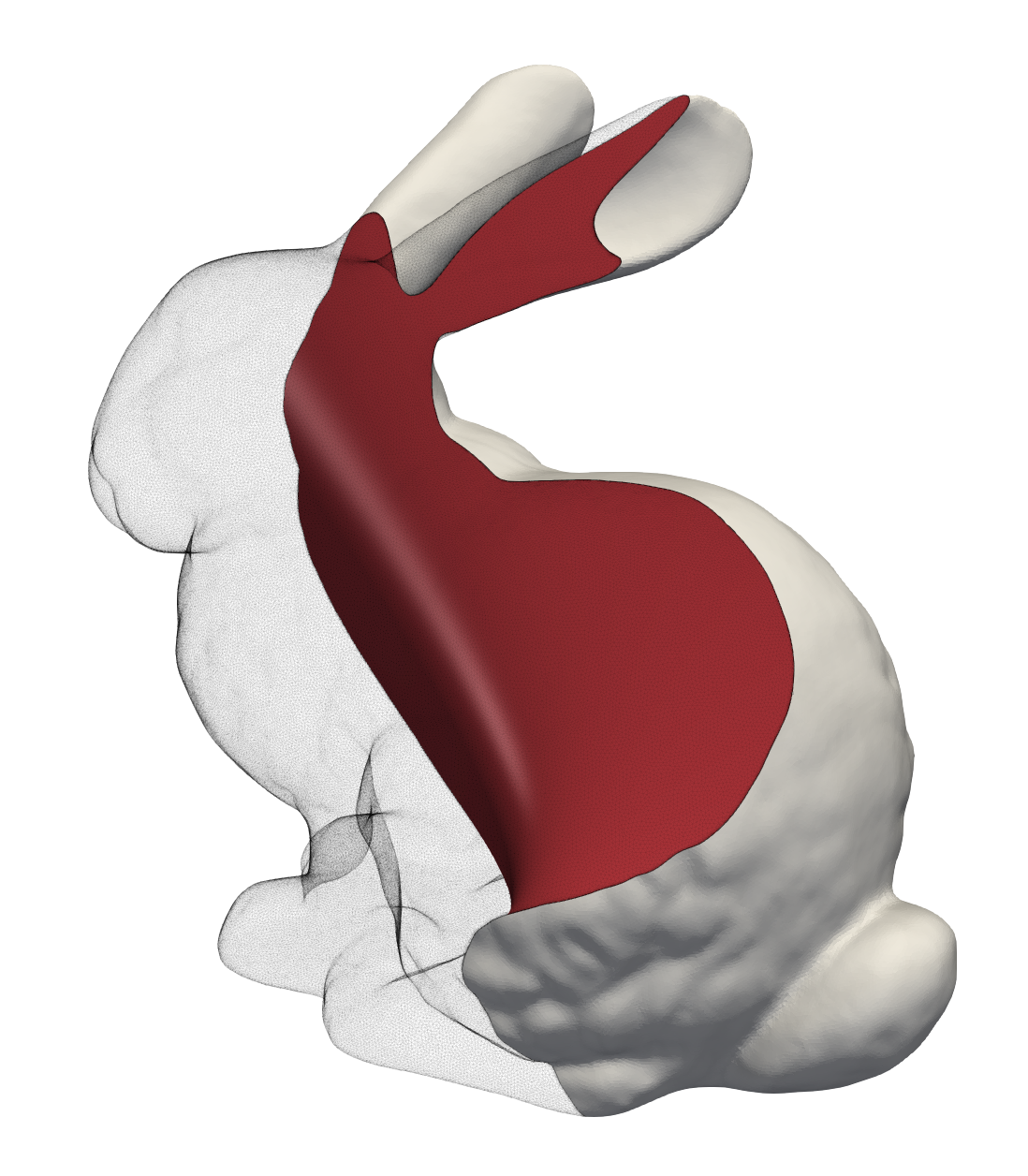}\end{minipage}}
        \subfloat[Clipped bunny (AMR ref.)]{\begin{minipage}{0.39\textwidth}\centering\adjincludegraphics[width=3cm,trim=0 0 0 0,clip=true]{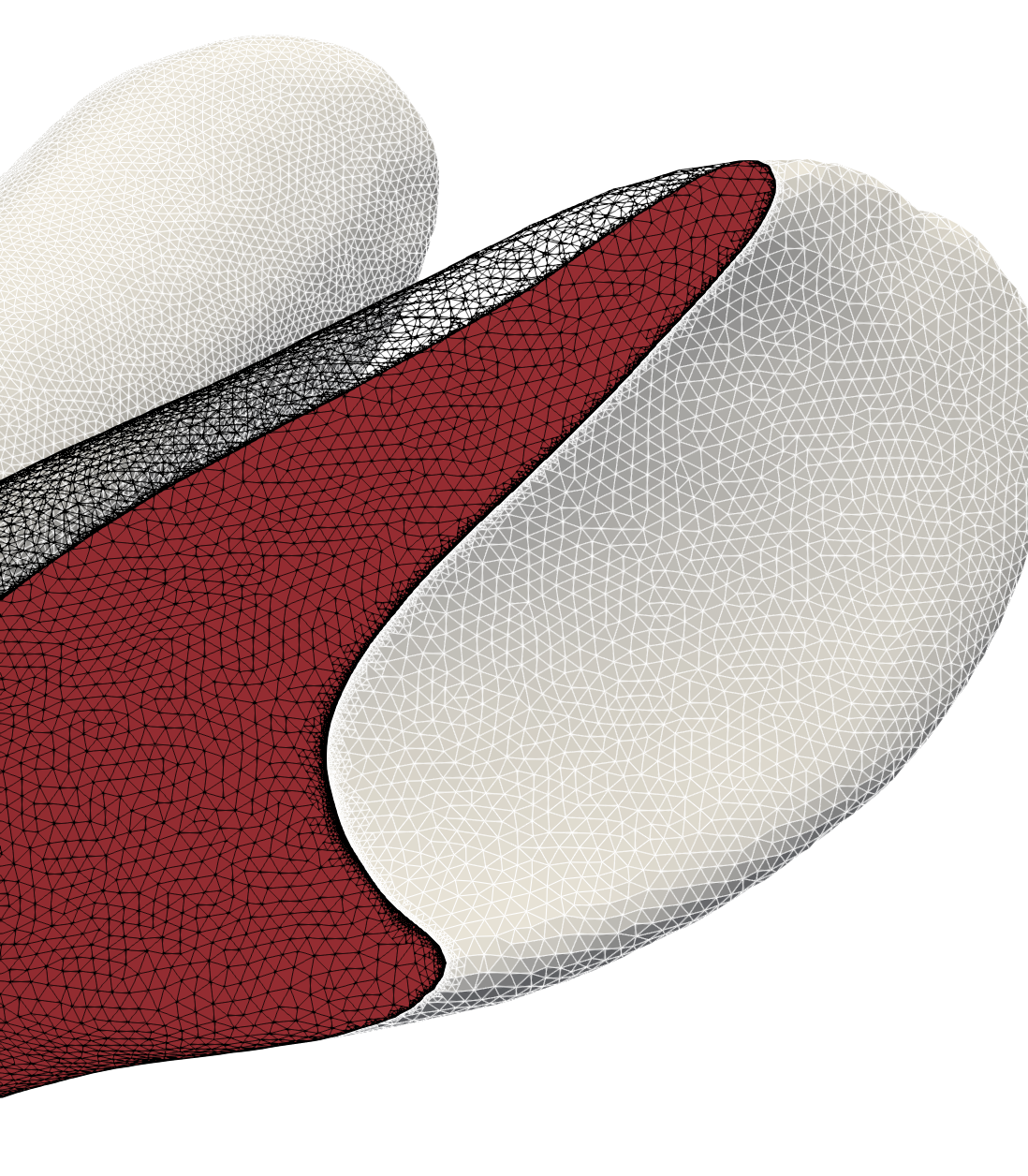}\end{minipage}}
    \end{tabular}
    \caption{Examples of random intersection cases between the five considered polyhedra and a paraboloid. The left column shows the full polyhedra; the central column shows the polyhedra clipped by a paraboloid; the right column shows the edges of the AMR of the polyhedra used for calculating the reference moments.}
    \label{fig:testcases}
\end{figure}

\begin{figure}[tbhp] \centering
    \begin{tabular}{c}
        \subfloat[]{\begin{minipage}{0.33\textwidth}\centering\adjincludegraphics[width=\textwidth,trim=0 0 0 0,clip=true]{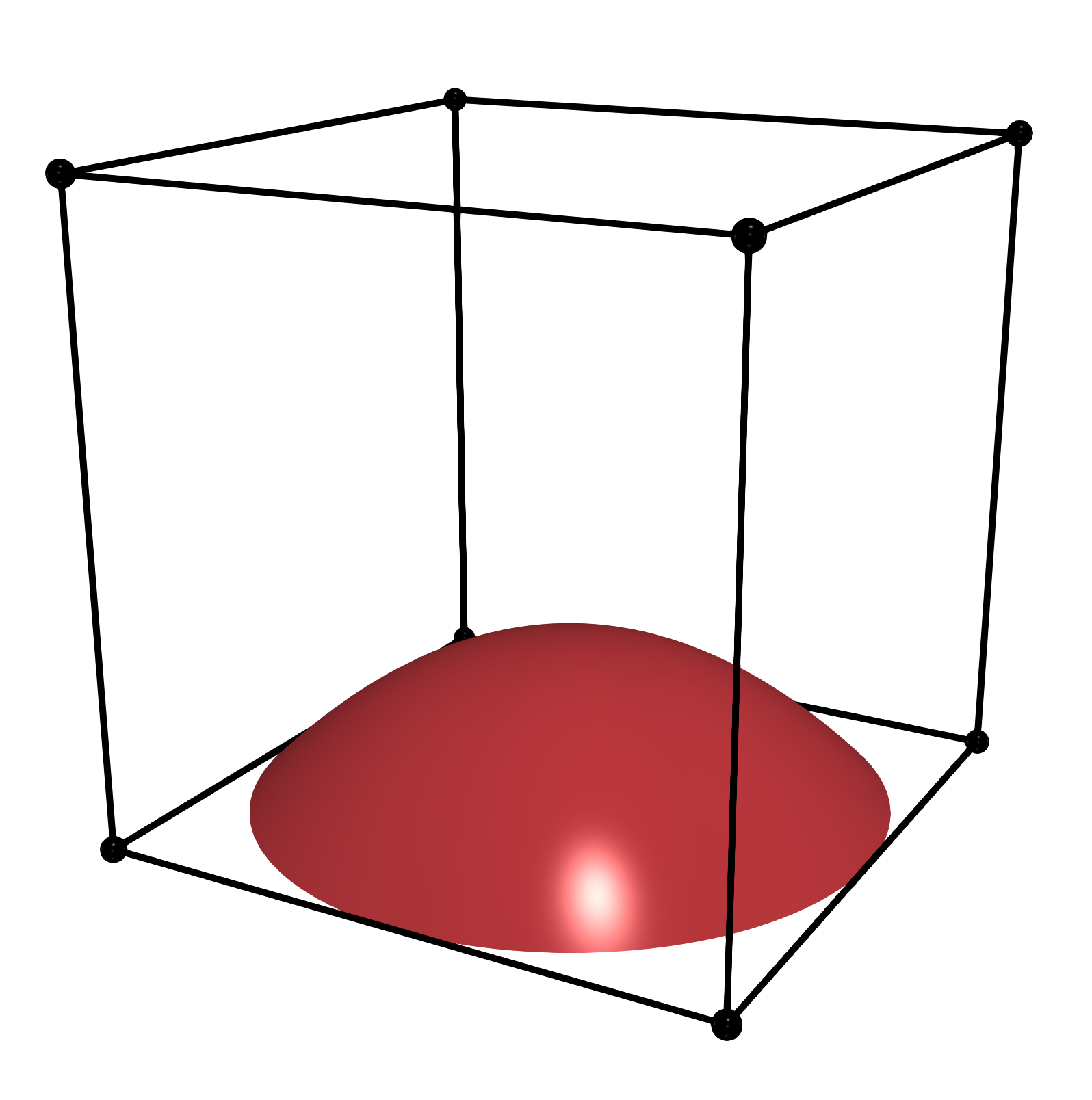}\end{minipage}}
        \subfloat[]{\begin{minipage}{0.33\textwidth}\centering\adjincludegraphics[width=\textwidth,trim=0 0 0 0,clip=true]{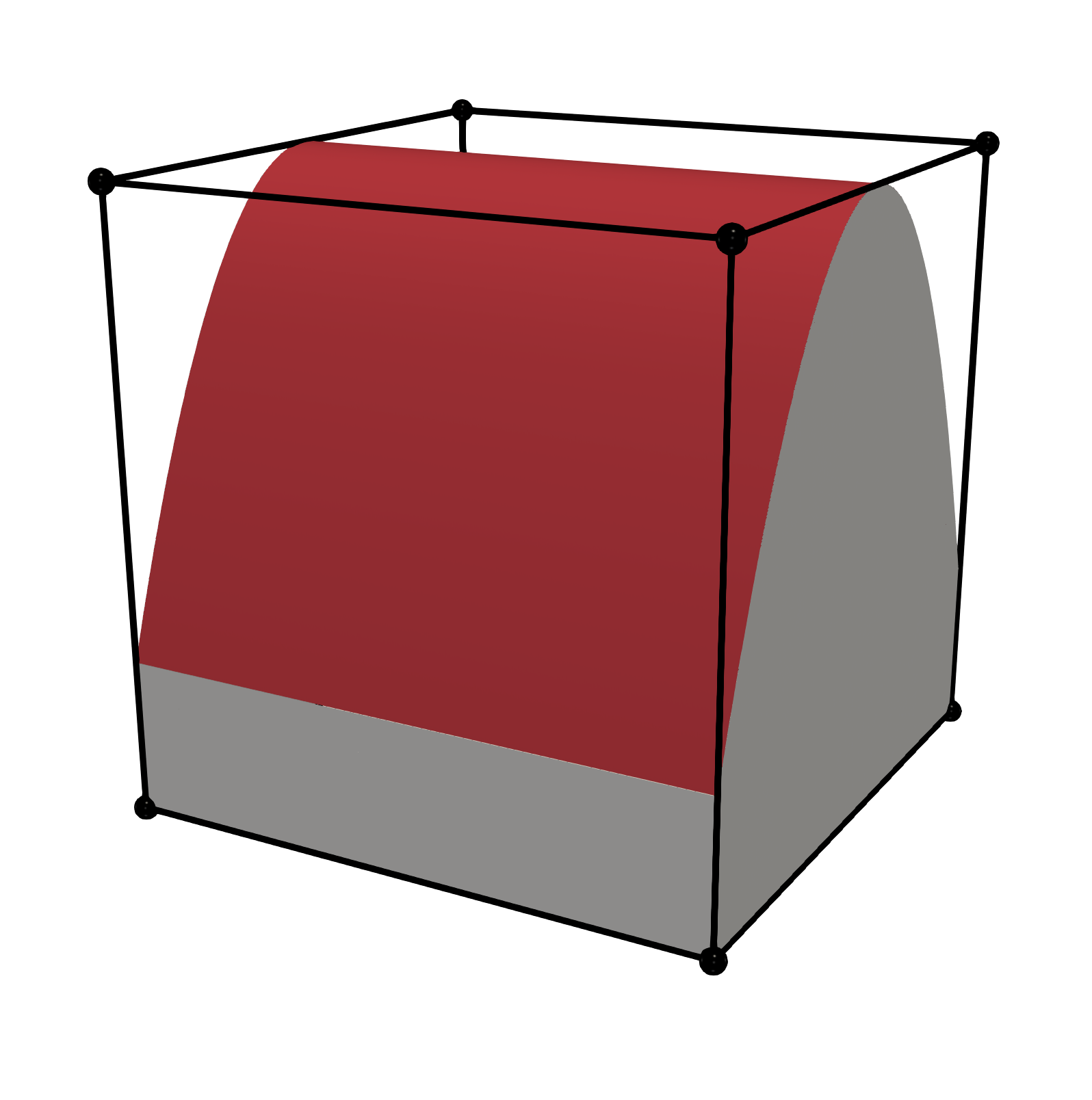}\end{minipage}}
        \subfloat[]{\begin{minipage}{0.33\textwidth}\centering\adjincludegraphics[width=\textwidth,trim=0 0 0 0,clip=true]{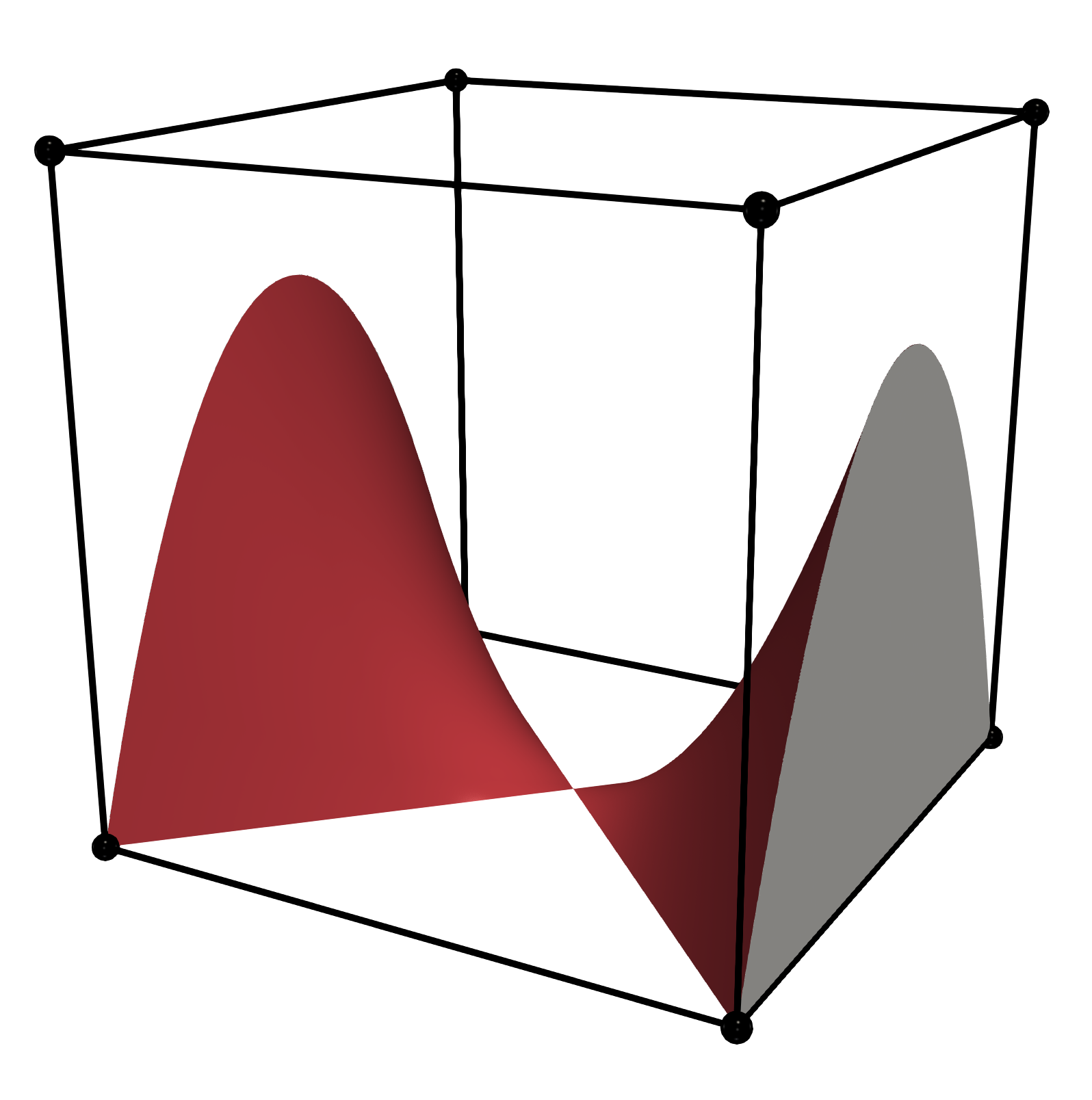}\end{minipage}}
    \end{tabular}
    \caption{Examples of cases with singular intersection configurations and/or ambiguous topologies covered by the graded parameter sweep: (a) the elliptic paraboloid is tangent to four edges of the polyhedron; (b) the parabolic cylinder is tangent to one face of the polyhedron, resulting in a degenerate conic section intersection that is made of two overlapping parallel lines; (c) the intersection of the hyperbolic paraboloid with one face of the polyhedron results in a degenerate hyperbola that is the intersection of two straight lines.}
    \label{fig:degeneratecases}
\end{figure}

A random parameter sweep is first conducted by uniform random sampling of the eight parameters ($t_x, t_y, t_z, \theta_x, \theta_y, \theta_z, \alpha, \beta$) in the parameter space $[-\tfrac{1}{2},\tfrac{1}{2}]^3 \times [-\pi,\pi]^3 \times [-\tfrac{1}{2},\tfrac{1}{2}]^2$, totalling more than $\smash{2\times10^8}$ realizations. A graded parameter sweep is then conducted, in which the eight parameters are chosen in the discrete parameter space $\{-\tfrac{1}{2},-\tfrac{1}{4},0,\tfrac{1}{4},\tfrac{1}{2}\}^3 \times\{-\pi,-\tfrac{\pi}{2},0,\tfrac{\pi}{2},\pi\}^3 \times \{-5,-4,\ldots,4,5\}^2$, resulting in $5^3\times5^3\times11^2 = 1,890,625$ distinct realizations for each geometry\footnote{We do not present the results of a graded parameter sweep on the Stanford bunny, since the ``organic'' nature of this polyhedron renders a graded parameter sweep equivalent to a random one.}. The graded parameter sweep differs from the random one in that it raises many singular intersection configurations (e.g., degenerate conic sections that consist of parallel or intersecting line segments, or conic sections that are parabolas) and/or ambiguous discrete topologies that arise from polyhedron vertices lying on the paraboloid or edges of the polyhedron being tangent to the paraboloid, therefore testing the robustness of our implementation as well as of the procedure described in \Cref{sec:robustness}. For each case of the random and graded parameter sweeps, the reference value of $\smash{\boldsymbol{\mathcal{M}}^{\hat{\mathcal{P}}}}$ for calculating the moment errors is obtained by adaptive mesh refinement (AMR) of the faces of $\smash{\mathcal{P}}$ into triangles, followed by the (exact) integration of $\smash{\boldsymbol{\Phi}_\mathcal{S}}$ and $\smash{\boldsymbol{\Phi}_{\mathcal{F}_i}}$, given in Eq.~\cref{eq:primitives}, over the triangles above and below the paraboloid, respectively. Examples of random intersection cases and their associated AMR are shown in \Cref{fig:testcases}, whereas examples of singular cases raised during the graded parameter sweep are shown in \cref{fig:degeneratecases}. \medskip


The maximum and average moment errors obtained during the random parameter sweep are given in \cref{tab:randomsweep}. For the tetrahedron, cube, dodecahedron, and hollow cube, the average moment error is of the order $\smash{\epsilon_{64} = 2^{-52} \sim \mathcal{O}(10^{-16})}$, whereas the maximum error is about one to two orders of magnitude larger. The average and maximum moment errors for the Stanford bunny, which contains $\mathcal{O}(10^5)$ faces, are each about one order of magnitude larger than for the other polyhedra but still close to machine-zero. The graded parameter sweep, whose results are given in \cref{tab:regularsweep}, exhibits similar moment errors as for the random parameter sweep, even though, by design, it raises many more singular intersection configurations and ambiguous discrete topologies than the random parameter sweep, therefore requiring a more frequent use of the nudging procedure described in \Cref{sec:robustness}.

\begin{table}[tbhp]
    \renewcommand\arraystretch{1.4}
    \scriptsize
    \captionsetup{position=top} 
    \caption{Random parameter sweep results. For each geometry, we provide: the number of tests conducted, the number of recursive levels used for the AMR reference moment calculation, as well as the average and maximum errors in the estimation of the zeroth and first moments.}\label{tab:randomsweep}
    \begin{center}\vspace{-1mm}
        \begin{tabular}{|r|R|R|R|R|R|R|} \hline
            \multicolumn{1}{|c|}{\multirow{2}{*}{\quad\;\textbf{Geometry}\;\quad}} & \multicolumn{1}{c|}{\multirow{2}{*}{\shortstack{\textbf{Number}\\\textbf{of tests}}}} & \multicolumn{1}{c|}{\multirow{2}{*}{\shortstack{\textbf{AMR}\\\textbf{levels}}}} & \multicolumn{2}{c|}{\textbf{Zeroth moment error}} &  \multicolumn{2}{c|}{\textbf{First moments error}} \\ \cline{4-7}
            & & & \multicolumn{1}{c|}{Average} & \multicolumn{1}{c|}{Maximum} &\multicolumn{1}{c|}{Average} & \multicolumn{1}{c|}{Maximum} \\ \hline
    Tetrahedron &  5 \times 10^7 & 17 & 2.3 \times 10^{-16} & 3.8 \times 10^{-15} & 1.6 \times 10^{-16} & 6.9 \times 10^{-14} \\ \hline
    Cube &         5 \times 10^7 & 17 & 2.4 \times 10^{-16} & 2.5 \times 10^{-15} & 8.7 \times 10^{-17} & 2.1 \times 10^{-14} \\ \hline
    Dodecahedron & 5 \times 10^7 & 17 & 3.5 \times 10^{-16} & 1.9 \times 10^{-15} & 6.8 \times 10^{-17} & 6.2 \times 10^{-15} \\ \hline
    Hollow cube &  5 \times 10^7 & 17 & 2.0 \times 10^{-16} & 3.2 \times 10^{-15} & 1.3 \times 10^{-16} & 3.4 \times 10^{-14} \\ \hline
    Stanford bunny &1 \times 10^3 &13 & 5.8 \times 10^{-15} & 4.8 \times 10^{-14} & 2.7 \times 10^{-15} & 3.1 \times 10^{-14} \\ \hline
    \end{tabular}\vspace{-1mm}
    \end{center}
\end{table}
\begin{table}[tbhp]
        \renewcommand\arraystretch{1.4}
        \scriptsize
        \captionsetup{position=top} 
        \caption{Graded parameter sweep results. For each geometry, we provide: the number of tests conducted, the number of recursive levels used for the AMR reference moment calculation, as well as the average and maximum errors in the estimation of the zeroth and first moments.}\label{tab:regularsweep}
        \begin{center}\vspace{-1mm}
        \begin{tabular}{|r|R|R|R|R|R|R|} \hline
        \multicolumn{1}{|c|}{\multirow{2}{*}{\quad\;\textbf{Geometry}\;\quad}} & \multicolumn{1}{c|}{\multirow{2}{*}{\shortstack{\textbf{Number}\\\textbf{of tests}}}} & \multicolumn{1}{c|}{\multirow{2}{*}{\shortstack{\textbf{AMR}\\\textbf{levels}}}} & \multicolumn{2}{c|}{\textbf{Zeroth moment error}} &  \multicolumn{2}{c|}{\textbf{First moments error}} \\ \cline{4-7}
        & & & \multicolumn{1}{c|}{Average} & \multicolumn{1}{c|}{Maximum} &\multicolumn{1}{c|}{Average} & \multicolumn{1}{c|}{Maximum} \\ \hline
        Tetrahedron &  5^6 \times 11^2 & 22 & 2.7 \times 10^{-16} & 4.9 \times 10^{-15} & 2.1 \times 10^{-16} & 3.9 \times 10^{-14} \\ \hline
        Cube &         5^6 \times 11^2 & 22 & 2.6 \times 10^{-16} & 4.7 \times 10^{-15} & 1.0 \times 10^{-16} & 2.1 \times 10^{-14} \\ \hline
        Dodecahedron & 5^6 \times 11^2 & 22 & 3.7 \times 10^{-16} & 1.7 \times 10^{-15} & 8.3 \times 10^{-17} & 5.4 \times 10^{-15} \\ \hline
        Hollow cube &  5^6 \times 11^2 & 22 & 2.0 \times 10^{-16} & 4.2 \times 10^{-15} & 1.6 \times 10^{-16} & 4.8 \times 10^{-14} \\ \hline
        \end{tabular}\vspace{-1mm}
        \end{center}
\end{table}

\subsection{Parameter sweep with nudge}\label{sec:paramnudge} 
In order to further test the robustness of our implementation and of the nudging procedure described in \Cref{sec:robustness}, we consider the same random parameter sweep for the cube geometry as described in \Cref{sec:paramsweep}, but with the addition of a post-sample translation of the polyhedron along~$\mathbf{e}_z$ so as for one of its vertices to lie exactly on the paraboloid. At least one iteration of the nudging procedure of \Cref{sec:robustness} is therefore required for each random case. The results of this random parameter sweep are given in \cref{tab:ambiguouscases}, displaying errors of a similar order of magnitude as for the random parameter sweep presented in \Cref{sec:paramsweep}.


\begin{table}[tbhp]
    \renewcommand\arraystretch{1.4}
    \scriptsize
    \captionsetup{position=top} 
    \caption{Random parameter sweep with one vertex of the polyhedron lying exactly on the paraboloid. We provide: the number of tests conducted, the number of recursive levels used for the AMR reference moment calculation, as well as the average and maximum errors in the estimation of the zeroth and first moments.}\label{tab:ambiguouscases}
    \begin{center}
        \begin{tabular}{|r|R|R|R|R|R|R|} \hline
            \multicolumn{1}{|c|}{\multirow{2}{*}{\quad\;\textbf{Geometry}\;\quad}} & \multicolumn{1}{c|}{\multirow{2}{*}{\shortstack{\textbf{Number}\\\textbf{of tests}}}} & \multicolumn{1}{c|}{\multirow{2}{*}{\shortstack{\textbf{AMR}\\\textbf{levels}}}} & \multicolumn{2}{c|}{\textbf{Zeroth moment error}} &  \multicolumn{2}{c|}{\textbf{First moments error}} \\ \cline{4-7}
            & & & \multicolumn{1}{c|}{Average} & \multicolumn{1}{c|}{Maximum} &\multicolumn{1}{c|}{Average} & \multicolumn{1}{c|}{Maximum} \\ \hline
        Cube (vertex on $\mathcal{S}$) & 5 \times 10^7 & 17 & 2.4 \times 10^{-16} & 7.8 \times 10^{-15} & 1.8 \times 10^{-16} & 4.5 \times 10^{-14} \\ \hline
    \end{tabular}
    \end{center}
\end{table}

\subsection{Timings} \label{sec:timings}
In order to assess the performances of our implementation, the time required for each moment estimation of the random parameter sweeps presented in \Cref{sec:paramsweep,sec:paramnudge} has been measured using OpenMP's \texttt{omp\_get\_wtime()} function \cite{Dagum1998}, which has a precision of 1 nanosecond on the workstation that we used. The characteristics of this workstation are summarized in \cref{tab:workstation}. The \texttt{C++} code implementing the closed-form expressions presented in \Cref{sec:moments} has been compiled with the GNU 10.3.0 suite of compilers \cite{GNU10.3}, using the flags given in \cref{tab:workstation}.

\begin{table}[tbhp]
    \renewcommand\arraystretch{1.4}
    \scriptsize
    \captionsetup{position=top} 
    \caption{Charateristics of the workstation used for the timing results presented in \Cref{sec:timings} (this is the same workstation as used in \cite{Chiodi2022}).}\label{tab:workstation}
    \begin{center}
        \begin{tabular}{|c|l|l|} \hline
            \multirow{9}{*}{\textbf{CPU}}   & vendor\_id & GenuineIntel \\ \cline{2-3}
            & CPU family & 6 \\ \cline{2-3}
            & Model & 158 \\ \cline{2-3}
            & Model name & Intel(R) Core(TM) i7-8700K CPU \@ 3.70 GHz \\ \cline{2-3}
            & Stepping & 10 \\ \cline{2-3}
            & Microcode & 0xca \\ \cline{2-3}
            & Min/max clock CPU frequency & 800 MHz -- 4.70 GHz \\ \cline{2-3}
            & CPU asserted frequency & 4.0 GHz \\ \cline{2-3}
            & Cache size & 12288 KB \\ \hline
            \multirow{3}{*}{\textbf{Compiler}}   & Suite & GNU \\ \cline{2-3}
            & Version & 10.3.0 \\ \cline{2-3}
            & Flags & -O3 -march=native -DNDEBUG -DNDEBUG\_PERF$^\ast$ \\ \hline
        \end{tabular}\vspace{1mm}
    $^\ast$ -DNDEBUG\_PERF is an IRL-specific compiler flag that disables additional debugging assertions \cite{Chiodi2022}.
    \end{center}
\end{table}
The timings are summarized in \cref{tab:timings}, which shows the average moment calculation time for the zeroth moment only, and for both the zeroth and first moments. Overall, the average time for calculating the first moments of a polyhedron clipped by a paraboloid is less than $1\ \mu\mathrm{s}$ per face of the original polyhedron. A direct comparison can be made with the half-space clipping of \cite{Chiodi2022}, which used the same workstation as the current work. It transpires that the clipping of a cube by a paraboloid is on average less than $6$ times more expensive than its clipping by a plane. When the nudging procedure of \Cref{sec:robustness} is required, as is the case for each realization of the modified random parameter sweep presented in \Cref{sec:paramnudge}, the timings of which are shown in the last row of \cref{tab:timings}, an increase in cost by a factor of about $200$ is observed. This is mostly due to the switch to quadruple precision that is operated in conjuction with the random translation and rotation of the polyhedron by a value of $\smash{\epsilon_\text{nudge}}$. It should be noted that such cases occur extremely rarely, unless they are manually engineered as in \Cref{sec:paramnudge}.

\begin{table}[tbhp]
    \renewcommand\arraystretch{1.4}
    \scriptsize
    \captionsetup{position=top} 
    \caption{Timings for the random parameter sweeps presented in \Cref{sec:paramsweep,sec:paramnudge}. The timings measured in \cite{Chiodi2022} for the clipping of a cube by a plane are provided for reference (they were measured on the same workstation as used for the current work).}\label{tab:timings}
    \begin{center}
        \begin{tabular}{|r|R|R|R|R|R|} \hline
            \multicolumn{1}{|c|}{\multirow{3}{*}{\quad\;\textbf{Geometry}\;\quad}} & \multicolumn{1}{c|}{\multirow{3}{*}{\shortstack{\textbf{Number}\\\textbf{of tests}}}} & \multicolumn{4}{c|}{\textbf{Average moment calculation time}} \\ \cline{3-6}
            & & \multicolumn{2}{c|}{Zeroth moment only} & \multicolumn{2}{c|}{Zeroth and first moments} \\ \cline{3-6}
            & & \multicolumn{1}{c|}{$\mu$s/test} & \multicolumn{1}{c|}{$\mu$s/test/face} & \multicolumn{1}{c|}{$\mu$s/test} & \multicolumn{1}{c|}{$\mu$s/test/face} \\ \hline
    Tetrahedron &                           5 \times 10^7 & 1.07 & 0.27 & 1.89 & 0.47 \\ \hline
    Cube &                                  5 \times 10^7 & 1.57 & 0.26 & 2.53 & 0.42 \\ \hline
    Cube (half-space clipping) \cite{Chiodi2022} &  15 \times 10^6 & 0.27 & 0.05 & -- & -- \\ \hline
    Dodecahedron &                          5 \times 10^7 & 2.66 &0.22 & 3.96 & 0.33 \\ \hline
    Hollow cube &                           5 \times 10^7 & 3.27 & 0.27 & 5.12 & 0.43 \\ \hline
    Stanford bunny &                        1 \times 10^3 & 6.22 \times 10^4 & 0.19 & 6.41 \times 10^4 & 0.19 \\ \hline\hline
    Cube (vertex on $\mathcal{S}$) &   5 \times 10^7 & 236 & 39.3 & 541 & 90.2 \\ \hline
\end{tabular}
    \end{center}
\end{table}

\section{Conclusions}
\label{sec:conclusions}
We have derived closed-form expressions for the first moments of a polyhedron clipped by a paraboloid, enabling their robust machine-accurate estimation at a computational cost that is considerably lower than with any other available approach. These expressions have been obtained by consecutive applications of the divergence theorem, transforming the three-dimensional integrals that are the zeroth and first moments of the clipped polyhedron into a sum of one-dimensional integrals. This requires parametrizing the conic section arcs resulting from the intersection of the paraboloid with the polyhedron's faces, which we have chosen to express as rational quadratic B\'ezier curves. The moments of the clipped polyhedron can, as a result, be expressed as the sum of three main contributions that are function of the polyhedron vertices and of the coefficients of the paraboloid. These expressions do not differ based on the type of paraboloid that is considered (elliptic, hyperbolic, or parabolic). Making use of this parametrization, we also show how to express integrals over the curved faces of the clipped polyhedron as the sums of one-dimensional integrals. When ambiguous discrete intersection topologies are detected, e.g., when the paraboloid is tangent to an edge of the polyhedron or intersects the polyhedron at the location of one of its vertices, a nudging procedure is triggered so as to guarantee robust moment estimations. A series of millions of intersection configurations that are randomly chosen, as well as manually engineered so as to raise singular intersection configurations and/or ambigious discrete topologies, have been tested. These showcase an average moment estimation error that is of the order of the machine-zero, and a maximum error that is about one order of magnitude larger. The timing of these moment estimations shows that the clipping of a polyhedron by a paraboloid, with our approach, is on average about 6 times more expensive than its clipping by a plane.

\section*{Reproducibility}
\label{sec:reproducibility}
The code used to produce the results presented in this manuscript is openly available as part of the \href{https://github.com/robert-chiodi/interface-reconstruction-library/tree/paraboloid_cutting}{\texttt{Interface Reconstruction Library~(IRL)}}. A \href{https://github.com/robert-chiodi/interface-reconstruction-library/blob/paraboloid_cutting/docs/markdown/paraboloid-polytope_intersection.md}{\texttt{step-by-step guide}} for reproducing the figures and tables of this manuscript using \texttt{IRL} can be found in the library's documentation. The commit hash corresponding to the version of the code used in this manuscript is \href{https://github.com/robert-chiodi/interface-reconstruction-library/commit/8e77b35d6f6f7888ae1a3527ff36857d30bed8d4}{\texttt{8e77b35}}.

\section*{Acknowledgments}
This project has received funding from the European Union's Horizon 2020 research and innovation programme under the Marie Skłodowska-Curie grant agreement No 101026017 \raisebox{-.125\height}{\includegraphics[height=2ex]{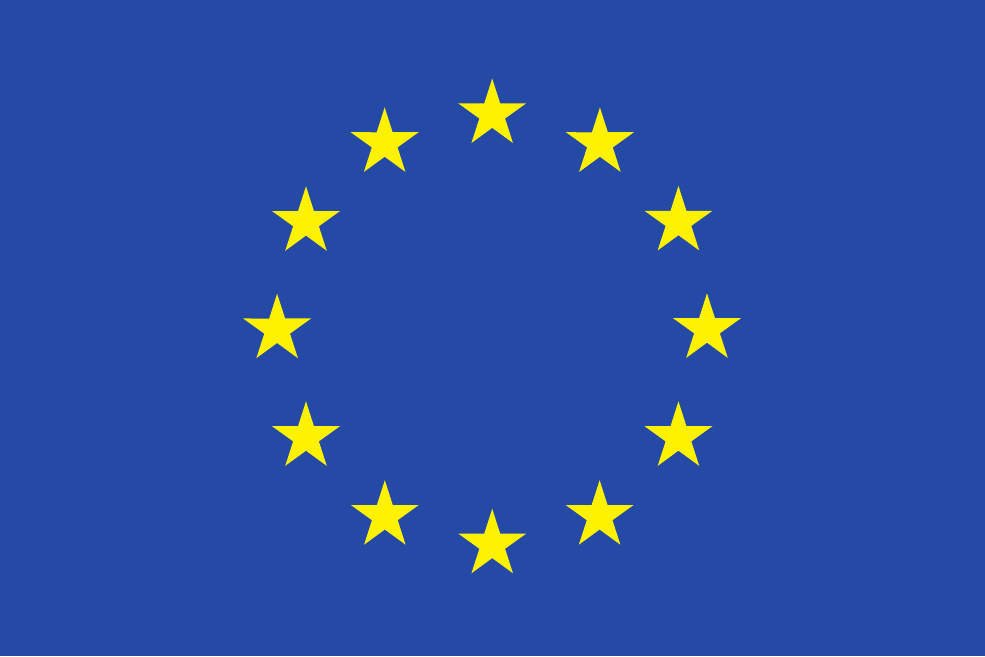}}\,. Robert Chiodi was sponsored by the Office of Naval Research (ONR) as part of the Multidisciplinary University Research Initiatives (MURI) Program, under grant number N00014-16-1-2617. The views and conclusions contained herein are those of the authors only and should not be interpreted as representing those of ONR, the U.S. Navy or the U.S. Government.

\appendix
\section{Third contribution to the moments} 
\begingroup
\allowdisplaybreaks
\label{apdx:M3}
The third contribution to the moments of face $\smash{\mathcal{F}_i}$, $\smash{\boldsymbol{\mathcal{M}}^{\hat{\mathcal{P}}_3}_{i}}$, is obtained by calculating the integral on the right-hand side of Eq.~\cref{eq:defM3} using the definitions in Eqs.~\cref{eq:primitives,eq:defpsiF,eq:defpsiS}, as well as the facts that $\smash{\forall j \in \{1,\ldots,n_{\partial\hat{\mathcal{F}_i}}\}}$,
\begin{align}
    {z}_{i,j} & = \delta_i - \lambda_i {x}_{i,j} - \tau_i {y}_{i,j} \ , \\
    1^{\partial\tilde{\mathcal{S}}}_{i,j} \ne 0 & \Leftrightarrow \left\{
        \begin{array}{ll} 
            {z}_{i,j,0} &= - \alpha {x}_{i,j,0}^2 - \beta {y}_{i,j,0}^2 \\  
            {z}_{i,j,1} &= - \alpha {x}_{i,j,1}^2 - \beta {y}_{i,j,1}^2
        \end{array}\right.  \ ,
\end{align}
and that the control point $\smash{\mathbf{x}_{i,j}^\star}$ belongs to the plane containing the face $\smash{\mathcal{F}_i}$ hence
\begin{align}
            {z}^\star_{i,j} &= \delta_i - \lambda_i {x}_{i,j}^\star - \tau_i {y}_{i,j}^\star \ .
\end{align}
After substitution of these expressions in Eq.~\cref{eq:defM3} and some simplification, $\smash{\boldsymbol{\mathcal{M}}^{\hat{\mathcal{P}}_3}_{i}}$ can be shown to read as in Eq.~\cref{eq:m3k}, with
\begin{equation}
    \boldsymbol{\mathcal{B}}^{(3)} \left(w, \mathbf{x}_a, \mathbf{x}_b, \mathbf{x}_c\right) = \diag\left(\boldsymbol{\mathcal{E}}(w)\right)\left( \boldsymbol{\mathcal{C}} (\mathbf{x}_a, \mathbf{x}_b, \mathbf{x}_c) \boldsymbol{\mathcal{K}} \boldsymbol{\mathcal{D}}(w) \right) \, ,
\end{equation}
where: 
$\boldsymbol{\mathcal{C}}$ is a $4 \times 12$ matrix whose non-zero coefficients contributing to $\mathrm{M}_0^{\hat{\mathcal{P}}}$ are given as
\begin{align}
    {\mathcal{C}}_{1,1}  & = \alpha \left(x_a^2 + 2 x_a x_b + x_b^2\right) + \beta \left(y_a^2 + 2 y_a y_b + y_b^2\right) - 2 z_a - 2 z_b \, ,\\
    {\mathcal{C}}_{1,2}  & = 
    \alpha \left(x_a^2+2 x_a x_b+x_b^2+4 x_a x_c+4 x_b x_c+4 x_c^2\right) \\ & \quad
    +\beta \left(y_a^2+2 y_a y_b+y_b^2+4 y_a y_c+4 y_b y_c+4 y_c^2\right) \nonumber\\ & \quad
    -4 z_a-4 z_b-8 z_c  \, ,\nonumber \\
    {\mathcal{C}}_{1,3}   & = \alpha x_c^2 + \beta y_c^2 - z_c \, ,
\end{align}
whose non-zero coefficients contributing to $\mathbf{e}_x\cdot\mathbf{M}_1^{\hat{\mathcal{P}}}$ are given as
\begin{align}
    {\mathcal{C}}_{2,4} & = 
    \beta \left(-12 {x_a} {y_a} {y_b}+12 {x_a} {y_b}^2+12 {x_b} {y_a}^2-12 {x_b} {y_a} {y_b}\right) \\ & \quad
    -6 {x_a} {z_a}+6 {x_a} {z_b}+6 {x_b} {z_a}-6 {x_b} {z_b}  \, ,\nonumber\\
    {\mathcal{C}}_{2,5} & = 
    \alpha \left(-12 {x_a} {x_c}^2+60 {x_a} {x_c} {x_b}-12 {x_c}^2 {x_b}\right) \\ & \quad
    +\beta \left(28 {x_a} {y_a} {y_c}-8 {x_a} {y_a} {y_b}-4 {x_a} {y_c}^2+20 {x_a} {y_c} {y_b}+8 {x_a} {y_b}^2 \right. \nonumber\\ & \quad\quad\quad
    -28 {x_c} {y_a}^2-8 {x_c} {y_a} {y_c}+20 {x_c} {y_a} {y_b}-8 {x_c} {y_c} {y_b}-28 {x_c} {y_b}^2\nonumber\\ & \quad\quad\quad
    \left.+8 {x_b} {y_a}^2+20 {x_b} {y_a} {y_c}-8 {x_b} {y_a} {y_b}-4 {x_b} {y_c}^2+28 {x_b} {y_c} {y_b}\right) \nonumber\\ & \quad
    +10 {x_a} {z_a}+20 {x_a} {z_c}+14 {x_a} {z_b}-22 {x_c} {z_a}-8 {x_c} {z_c} \nonumber\\ & \quad\quad\quad
    -22 {x_c} {z_b}+14 {x_b} {z_a}+20 {x_b} {z_c}+10 {x_b} {z_b} \, , \nonumber\\
    {\mathcal{C}}_{2,6}  & = 
    \alpha \left(-36 {x_a} {x_c}^2+12 {x_a} {x_c} {x_b}+12 {x_c}^3-36 {x_c}^2 {x_b}\right) \\ & \quad
    +\beta \left(-12 {x_a} {y_c}^2+4 {x_a} {y_c} {y_b}-24 {x_c} {y_a} {y_c}+4 {x_c} {y_a} \right. \nonumber\\ & \quad\quad\quad
    \left.{y_b}+12 {x_c} {y_c}^2-24 {x_c} {y_c} {y_b}+4 {x_b} {y_a} {y_c}-12 {x_b} {y_c}^2\right) \nonumber\\ & \quad
    -10 {x_a} {z_c}+2 {x_a} {z_b}-10 {x_c} {z_a}-12 {x_c} {z_c}-10 {x_c} {z_b} \nonumber\\ & \quad\quad\quad
    +2 {x_b} {z_a}-10 {x_b} {z_c} \, , \nonumber\\
    {\mathcal{C}}_{2,7}  & = 2 \alpha {x_c}^3+2 \beta {x_c} {y_c}^2+2 {x_c} {z_c} \, ,
\end{align}
whose non-zero coefficients contributing to $\mathbf{e}_y\cdot\mathbf{M}_1^{\hat{\mathcal{P}}}$ are given as
\begin{align}
    {\mathcal{C}}_{3,4}  & = \alpha \left(12 {x_a}^2 {y_b}-12 {x_a} {x_b} {y_a}-12 {x_a} {x_b} {y_b}+12 {x_b}^2 {y_a}\right) \\ & \quad
    -6 {y_a} {z_a}+6 {y_a} {z_b}+6 {y_b} {z_a}-6 {y_b} {z_b} \, \nonumber\\
    {\mathcal{C}}_{3,5} & = 
    \alpha \left(-28 {x_a}^2 {y_c}+8 {x_a}^2 {y_b}+28 {x_a} {x_c} {y_a}-8 {x_a} {x_c} {y_c}+20 {x_a} {x_c} {y_b} \right. \\ & \quad\quad\quad
    -8 {x_a} {x_b} {y_a}+20 {x_a} {x_b} {y_c}-8 {x_a} {x_b} {y_b}-4 {x_c}^2 {y_a}-4 {x_c}^2 {y_b} \nonumber\\ & \quad\quad\quad
    \left.+20 {x_c} {x_b} {y_a}-8 {x_c} {x_b} {y_c}+28 {x_c} {x_b} {y_b}+8 {x_b}^2 {y_a}-28 {x_b}^2 {y_c}\right) \nonumber\\ & \quad
    +\beta \left(-12 {y_a} {y_c}^2+60 {y_a} {y_c} {y_b}-12 {y_c}^2 {y_b}\right) \nonumber\\ & \quad\quad\quad
    +10 {y_a} {z_a}+20 {y_a} {z_c}+14 {y_a} {z_b}-22 {y_c} {z_a}-8 {y_c} {z_c} \nonumber\\ & \quad
    -22 {y_c} {z_b}+14 {y_b} {z_a}+20 {y_b} {z_c}+10 {y_b} {z_b} \, , \nonumber\\
    {\mathcal{C}}_{3,6} & = 
    \alpha \left(-24 {x_a} {x_c} {y_c}+4 {x_a} {x_c} {y_b}+4 {x_a} {x_b} {y_c}-12 {x_c}^2 {y_a} \right. \\ & \quad\quad\quad
    \left.+12 {x_c}^2 {y_c}-12 {x_c}^2 {y_b}+4 {x_c} {x_b} {y_a}-24 {x_c} {x_b} {y_c}\right) \nonumber\\ & \quad
    +\beta \left(-36 {y_a} {y_c}^2+12 {y_a} {y_c} {y_b}+12 {y_c}^3-36 {y_c}^2 {y_b}\right) \nonumber\\ & \quad\quad\quad
    -10 {y_a} {z_c}+2 {y_a} {z_b}-10 {y_c} {z_a}-12 {y_c} {z_c}-10 {y_c} {z_b} \nonumber\\ & \quad
    +2 {y_b} {z_a}-10 {y_b} {z_c} \, , \nonumber \\
     {\mathcal{C}}_{3,7} & = 2 \alpha {x_c}^2 {y_c}+2 \beta {y_c}^3+2 {y_c} {z_c} \, ,
\end{align}
and whose non-zero coefficients contributing to $\mathbf{e}_z\cdot\mathbf{M}_1^{\hat{\mathcal{P}}}$ are given as
\begin{align}
    {\mathcal{C}}_{4,8} & =
    +{\alpha\beta} \left(-42 {y_a}^2 {x_a}^2-10 {y_b}^2 {x_a}^2-28 {y_a} {y_b} {x_a}^2-28 {x_b} {y_a}^2 {x_a}\right. \\ & \quad\quad\quad
    \left.-28 {x_b} {y_b}^2 {x_a}-40 {x_b} {y_a} {y_b} {x_a}-10 {x_b}^2 {y_a}^2-42 {x_b}^2 {y_b}^2-28 {x_b}^2 {y_a} {y_b}\right) \nonumber \\ & \quad
    +  \alpha^2 \left(-21 {x_a}^4-28 {x_b} {x_a}^3-30 {x_b}^2 {x_a}^2-28 {x_b}^3 {x_a}-21 {x_b}^4\right) \nonumber \\ & \quad
    +\beta^2 \left(-21 {y_a}^4-28 {y_b} {y_a}^3-30 {y_b}^2 {y_a}^2-28 {y_b}^3 {y_a}-21 {y_b}^4\right) \nonumber \\ & \quad
    40 {z_a}^2+40 {z_b}^2 +48 {z_a} {z_b} \, , \nonumber\\
    {\mathcal{C}}_{4,9} & = 
    {\alpha\beta} \left(-7 {y_c}^2 {x_a}^2-10 {y_b}^2 {x_a}^2+63 {y_a} {y_c} {x_a}^2-21 {y_a} {y_b} {x_a}^2+35 {y_c} {y_b} {x_a}^2 \right. \\ & \quad\quad\quad
    +63 {x_c} {y_a}^2 {x_a}-21 {x_b} {y_a}^2 {x_a}-10 {x_b} {y_c}^2 {x_a}+35 {x_c} {y_b}^2 {x_a}  \nonumber \\ & \quad\quad\quad
    -21 {x_b} {y_b}^2 {x_a}-28 {x_c} {y_a} {y_c} {x_a}+70 {x_b} {y_a} {y_c} {x_a}+70 {x_c} {y_a} {y_b} {x_a}  \nonumber \\ & \quad\quad\quad
    -40 {x_b} {y_a} {y_b} {x_a}-20 {x_c} {y_c} {y_b} {x_a}+70 {x_b} {y_c} {y_b} {x_a}-7 {x_c}^2 {y_a}^2  \nonumber \\ & \quad\quad\quad
    -10 {x_b}^2 {y_a}^2+35 {x_c} {x_b} {y_a}^2-7 {x_b}^2 {y_c}^2-7 {x_c}^2 {y_b}^2+63 {x_c} {x_b} {y_b}^2  \nonumber \\ & \quad\quad\quad
    +35 {x_b}^2 {y_a} {y_c}-20 {x_c} {x_b} {y_a} {y_c}-10 {x_c}^2 {y_a} {y_b}-21 {x_b}^2 {y_a} {y_b}  \nonumber \\ & \quad\quad\quad
    \left.+70 {x_c} {x_b} {y_a} {y_b}+63 {x_b}^2 {y_c} {y_b}-28 {x_c} {x_b} {y_c} {y_b}\right) \nonumber \\ & \quad
    +\alpha^2 \left(63 {x_c} {x_a}^3-21 {x_b} {x_a}^3-21 {x_c}^2 {x_a}^2-30 {x_b}^2 {x_a}^2+105 {x_c} {x_b} {x_a}^2 \right. \nonumber \\ & \quad\quad\quad
    \left.-21 {x_b}^3 {x_a}+105 {x_c} {x_b}^2 {x_a}-30 {x_c}^2 {x_b} {x_a}+63 {x_c} {x_b}^3-21 {x_c}^2 {x_b}^2\right) \nonumber \\ & \quad
    +\beta^2 \left(63 {y_c} {y_a}^3-21 {y_b} {y_a}^3-21 {y_c}^2 {y_a}^2-30 {y_b}^2 {y_a}^2+105 {y_c} {y_b} {y_a}^2 \right. \nonumber \\ & \quad\quad\quad
    \left.-21 {y_b}^3 {y_a}+105 {y_c} {y_b}^2 {y_a}-30 {y_c}^2 {y_b} {y_a}+63 {y_c} {y_b}^3-21 {y_c}^2 {y_b}^2\right) \nonumber \\ & \quad
    -30 {z_a}^2+12 {z_c}^2-30 {z_b}^2-60 {z_a} {z_c}-24 {z_a} {z_b}-60 {z_c} {z_b} \, , \nonumber\\
    {\mathcal{C}}_{4,10} & = 
    {\alpha\beta} \left(-56 {y_c}^2 {x_a}^2-2 {y_b}^2 {x_a}^2+28 {y_c} {y_b} {x_a}^2+84 {x_c} {y_c}^2 {x_a} \right. \\ & \quad\quad\quad
    -92 {x_b} {y_c}^2 {x_a}+28 {x_c} {y_b}^2 {x_a}-224 {x_c} {y_a} {y_c} {x_a}+56 {x_b} {y_a} {y_c} {x_a}\nonumber \\ & \quad\quad\quad
    +56 {x_c} {y_a} {y_b} {x_a}-8 {x_b} {y_a} {y_b} {x_a}-184 {x_c} {y_c} {y_b} {x_a}\nonumber \\ & \quad\quad\quad
    +56 {x_b} {y_c} {y_b} {x_a}-56 {x_c}^2 {y_a}^2-2 {x_b}^2 {y_a}^2+28 {x_c} {x_b} {y_a}^2\nonumber \\ & \quad\quad\quad
    -12 {x_c}^2 {y_c}^2-56 {x_b}^2 {y_c}^2+84 {x_c} {x_b} {y_c}^2-56 {x_c}^2 {y_b}^2\nonumber \\ & \quad\quad\quad
    +84 {x_c}^2 {y_a} {y_c}+28 {x_b}^2 {y_a} {y_c}-184 {x_c} {x_b} {y_a} {y_c} \nonumber \\ & \quad\quad\quad
    \left.-92 {x_c}^2 {y_a} {y_b}+56 {x_c} {x_b} {y_a} {y_b}+84 {x_c}^2 {y_c} {y_b}-224 {x_c} {x_b} {y_c} {y_b}\right) \nonumber \\ & \quad
    +\alpha^2 \left(-6 {x_c}^4+84 {x_a} {x_c}^3+84 {x_b} {x_c}^3-168 {x_a}^2 {x_c}^2-168 {x_b}^2 {x_c}^2 \right. \nonumber \\ & \quad\quad\quad
    \left.-276 {x_a} {x_b} {x_c}^2+84 {x_a} {x_b}^2 {x_c}+84 {x_a}^2 {x_b} {x_c}-6 {x_a}^2 {x_b}^2\right)  \nonumber \\ & \quad
    +\beta^2 \left(-6 {y_c}^4+84 {y_a} {y_c}^3+84 {y_b} {y_c}^3-168 {y_a}^2 {y_c}^2-168 {y_b}^2 {y_c}^2\right. \nonumber \\ & \quad\quad\quad
    \left.-276 {y_a} {y_b} {y_c}^2+84 {y_a} {y_b}^2 {y_c}+84 {y_a}^2 {y_b} {y_c}-6 {y_a}^2 {y_b}^2\right) \nonumber \\ & \quad
    +15 {z_a}^2+24 {z_c}^2+15 {z_b}^2+120 {z_a} {z_c}-6 {z_a} {z_b}+120 {z_c} {z_b} \, , \nonumber \\
    {\mathcal{C}}_{4,11} & = {\alpha\beta} \left(-12 {y_c}^2 {x_c}^2+14 {y_a} {y_c} {x_c}^2-2 {y_a} {y_b} {x_c}^2 \right. \\ & \quad\quad\quad
    +14 {y_c} {y_b} {x_c}^2+14 {x_a} {y_c}^2 {x_c}+14 {x_b} {y_c}^2 {x_c}\nonumber \\ & \quad\quad\quad
    \left.-4 {x_b} {y_a} {y_c} {x_c}-4 {x_a} {y_c} {y_b} {x_c}-2 {x_a} {x_b} {y_c}^2\right)\nonumber \\ & \quad
    +\alpha^2 \left(-6 {x_c}^4+14 {x_a} {x_c}^3+14 {x_b} {x_c}^3-6 {x_a} {x_b} {x_c}^2\right)  \nonumber \\ & \quad
    +\beta^2 \left(-6 {y_c}^4+14 {y_a} {y_c}^3+14 {y_b} {y_c}^3-6 {y_a} {y_b} {y_c}^2\right)\nonumber \\ & \quad
    -7 {z_c}^2-5 {z_a} {z_c}+{z_a} {z_b}-5 {z_c} {z_b} \, , \nonumber\\
    {\mathcal{C}}_{4,12} & = -2 {\alpha\beta} {y_c}^2 {x_c}^2 -\alpha^2 {x_c}^4-\beta^2 {y_c}^4+{z_c}^2\, ;
\end{align}  
$\boldsymbol{\mathcal{K}}$ is the $12\times10$ matrix given as
\begin{equation}
   {\normalsize \boldsymbol{\mathcal{K}}} = {\small\begin{bmatrix}[1.2]
        -\frac{3}{8} & \frac{31}{48} & -\frac{7}{8} & -\frac{1}{16} & 0 & \frac{1}{24} & 0 & 0 & 0 & 0 \\ 
        0 & -\frac{1}{6} & \frac{5}{8} & -\frac{3}{16} & 0 & \frac{1}{24} & 0 & 0 & 0 & 0 \\
        0 & \frac{2}{3} & -3 & \frac{11}{6} &-2 & 0 &0 &0 &0 & 0\\ 
        -\frac{1}{32} & \frac{93}{2240} & 0 & -\frac{163}{3360} & 0 & \frac{5}{168} & 0 & -\frac{1}{140} & 0 & 0 \\
        0 & \frac{1}{70} & -\frac{1}{16} & \frac{29}{1120} & 0 & -\frac{19}{1680} & 0 & \frac{1}{420} & 0 & 0 \\
        0 & -\frac{1}{210} & 0 & \frac{1}{21} & -\frac{1}{8} & \frac{13}{560} & 0 & -\frac{1}{280} & 0 & 0 \\
        0 & \frac{1}{35} & 0 & -\frac{16}{105} & 0 & \frac{56}{105} & -1 & \frac{1}{14} & 0 & 0 \\
        -\frac{1}{128} & \frac{193}{16128} & 0 & -\frac{149}{8064} & 0 & \frac{19}{1120} & 0 & \frac{41}{5040} & 0 & \frac{1}{630} \\
        0 & \frac{4}{945} & -\frac{1}{48} & \frac{65}{6048} & 0 & -\frac{1}{144} & 0 & \frac{11}{3780} & 0 & -\frac{1}{1890} \\
        0 & -\frac{1}{1890} & 0 & \frac{13}{1890} & -\frac{1}{48} & \frac{11}{2016} & 0 & -\frac{5}{3024} & 0 & \frac{1}{3780} \\
        0 & \frac{1}{315} & 0 & -\frac{1}{45} & 0 & \frac{4}{35} & -\frac{1}{4} & \frac{17}{504} & 0 & -\frac{1}{252} \\
        0 & -\frac{1}{63} & 0 & \frac{29}{315} & 0 & -\frac{26}{105} & 0 & \frac{194}{315} & -1 & \frac{1}{18}
    \end{bmatrix}} \, ;
\end{equation}
$\boldsymbol{\mathcal{D}}$ is the vector given as
\begin{align}
    \boldsymbol{\mathcal{D}}(w) & = \begin{bmatrix}[1.2] 0 & w^2 & 0 &w^4 & 0  &w^6 & 0& w^8 &0 & w^{10}  \end{bmatrix}^{\intercal} \\ \nonumber
    & \quad\quad\quad + \Theta(w) \begin{bmatrix}[1.2] w & 0 & w^3 &0 & w^5  & 0 &w^7 & 0& w^9 & 0  \end{bmatrix}^{\intercal} \, ;
\end{align}
and $\boldsymbol{\mathcal{E}}$ is the vector given as
\begin{equation}
    \boldsymbol{\mathcal{E}}(w) = \begin{bmatrix}[1.2] \Lambda(w)^3 & \Lambda(w)^4 & \Lambda(w)^4 & \Lambda(w)^5  \end{bmatrix}^{\intercal} \, ,
\end{equation}
with
\begin{align}
    \Theta(w) & = \left\{ \begin{array}{ll} \text{arctan}\left(\dfrac{1 - w}{\sqrt{1 - w^2}}\right)\dfrac{1}{\sqrt{1-w^2}} & 0 < w < 1 \\ & \\
            \text{arctanh}\left(\dfrac{w - 1}{\sqrt{w^2 - 1}}\right)\dfrac{1}{\sqrt{w^2-1}} & 1 \le w \end{array} \right.  \, , \\
            \Lambda(w) & = \frac{1}{(w - 1) (w + 1)}  \, .
\end{align}
It should be noted that the naive implementation of these expressions may lead to significant round-off errors when $w$ is in the vicinity of $1$. To avoid such problems, we resort to the Taylor series expansion of $\smash{\boldsymbol{\mathcal{B}}^{(3)}}$ around $w = 1$ for its numerical estimation. Using $64$-bit floating-point arithmetics, we have found that the Taylor series expansion of $\smash{\boldsymbol{\mathcal{B}}^{(3)}}$ to order $40$ for $w \in [0.35,1.7]$ is sufficient for producing near machine-zero estimates. This implementation has been used for producing the results presented in \Cref{sec:tests}.
\endgroup

\bibliographystyle{siamplain}
\bibliography{/Users/fabien/Documents/Research/BvWResearch.bib}
\end{document}